\documentstyle[12pt]{article}                  % Use this for Latex (2.09)

\input psfig.sty

\input{amssymb.sty}
\makeatletter
\def\nothing#1{}
\newdimen\earraycolsep
\renewcommand{\thetable}{\arabic{table}}
\renewcommand{\thefigure}{\arabic{figure}}
\renewcommand{\title}[1]{%
 \vspace*{120\p@}%
  {\parindent \z@ \raggedright \reset@font
    \bfseries #1\par
    \nobreak
    \vskip 36\p@
  }}
\def\author#1{{\pretolerance=10000 \raggedright \advance \leftskip by 1in
\noindent #1 \vskip 1pc}}
\def\affiliation#1{{\advance\leftskip by 1in \noindent #1 \vskip -1pc}}
\def\refnote#1{{$^{\hbox{\scriptsize #1}}$}}
\def\affnote#1{\llap{$^{\hbox{\scriptsize #1}}$}}

\renewcommand\section{\@startsection{section}{1}{\z@}{2pc \@plus
      1ex minus .2ex}{1pc \@plus .2ex}{\reset@font
      \normalsize\bfseries\noindent
      {\addtocounter{section}{1}}\Roman{section}\
      {\setcounter{subsection}{0}
      \setcounter{subsubsection}{0}\setcounter{equation}{0}} }}
\renewcommand\subsection{\@startsection{subsection}{2}{\z@}{1pc \@plus 1ex
    minus.2ex}{1pc \@plus .2ex}
    {\reset@font\normalsize\bfseries
    \noindent{\addtocounter{subsection}{1}}%
    {\setcounter{subsubsection}{0}}\Roman{section}.\Roman{subsection}\ }}
\renewcommand\subsubsection{\@startsection{subsubsection}{3}{\parindent}
        {1pc \@plus 1ex minus.2ex}{-0.5em}{\reset@font\normalsize\bfseries%
        {\addtocounter{subsubsection}{1}} \hspace*{.6cm}
        \Roman{section}.\Roman{subsection}.\Roman{subsubsection}
        \hspace*{-7mm}}}

\def\AmS{{\protect\the\textfont2%
        A\kern-.1667em\lower.5ex\hbox{M}\kern-.125emS}}

\def\p@LaTeX{{\family{times}\series{m}\shape{n}\selectfont
L\kern-.36em\raise.3ex\hbox{\scriptsize A}\kern-.15em
T\kern-.1667em\lower.7ex\hbox{E}\kern-.125emX}}

\newlength{\colwidth}

\def\@oddhead{\hfil}
\def\@evenhead{\hfil}
\def\@oddfoot{{\bfseries\hfil\thepage}}
\def\@evenfoot{{\bfseries\thepage\hfil}}
\def\fnum@figure{\footnotesize\raggedright{\bfseries \figurename~\thefigure.}}
\def\fnum@table{\normalsize\raggedright{\bfseries \tablename~\thetable.}}
\long\def\@makecaption#1#2{\vskip 10\p@ {#1 #2\par}}
\long\def\@makefntext#1{\setbox0=\hbox{$\m@th^{\@thefnmark}$}\noindent
\hangindent=\wd0 \box0 #1}
\newbox\@atbox
\long\def\atable#1#2#3{\begin{table}[tbp]\centering\footnotesize
\setbox\@atbox\hbox{#2}
\parbox{\wd\@atbox}{\caption{#1}}\par\smallskip
#2
\par\smallskip\parbox{\wd\@atbox}{\raggedright #3}
\end{table}}

\def\Fb{{\bb F}}

\def\sbsneq{_{\sbs \atop \neq}}
\def\wh{\widehat}
\def\Cb{{\mathbb C}}
\def\Hb{{\mathbb H}}
\def\Nb{{\mathbb N}}
\def\Rb{{\mathbb R}}
\def\Tb{{\mathbb T}}
\def\Zb{{\mathbb Z}}
\def\L{\Lambda}
\def\Ac{{\cal A}}
\def\Bc{{\cal B}}
\def\Ac{{\cal A}}
\def\Hc{{\cal H}}
\def\Lc{{\cal L}}
\def\Nc{{\cal N}}

\def\Sc{{\cal S}}
\def\Uc{{\cal U}}

\def\Ec{{\cal E}}
\def\Hc{{\cal H}}

\def\Kc{{\cal K}}
\def\Lc{{\cal L}}

\def\Sc{{\cal S}}
\def\Uc{{\cal U}}
\newcommand{\wt}{\widetilde}

\def\a{\alpha}
\def\b{\beta}
\def\d{\delta}
\def\lb{\lambda}
\def\g{\gamma}
\def\om{\omega}
\def\s{\sigma}
\def\t{\theta}
\def\ve{\varepsilon}
\def\vp{\varphi}

\def\z{\zeta}

\def\D{\Delta}
\def\G{\Gamma}
\def\Lb{\Lambda}
\def\Om{\Omega}
\def\Si{\Sigma}

\def\ex{\exists}
\def\fl{\forall}
\def\ify{\infty}
\def\lgl{\langle}
\def\nb{\nabla}
\def\op{\oplus}
\def\ot{\otimes}
\def\ov{\overline}
\def\part{\partial}
\def\rgl{\rangle}
\def\sbs{\subset}
\def\semi{>\!\!\!\lhd}
\def\sm{\simeq}
\def\ts{\times}
\def\wdg{\wedge}

\def\ra{\rightarrow}
\def\longra{\longrightarrow}
\def\Ra{\Rightarrow}

\def\text{\hbox}

\def\Aut{\mathop{\rm Aut}\nolimits}
\def\Diff{\mathop{\rm Diff}\nolimits}
\def\Dom{\mathop{\rm Dom}\nolimits}
\def\End{\mathop{\rm End}\nolimits}
\def\id{\mathop{\rm id}\nolimits}
\def\Index{\mathop{\rm Index}\nolimits}

\def\Int{\mathop{\rm Int}\nolimits}
\def\Ker{\mathop{\rm Ker}\nolimits}
\def\Rang{\mathop{\rm Rang}\nolimits}
\def\Re{\mathop{\rm Re}\nolimits}
\def\Res{\mathop{\rm Res}\nolimits}
\def\Sign{\mathop{\rm Sign}\nolimits}

\def\Sup{\mathop{\rm Sup}\nolimits}
\def\Trace{\mathop{\rm Trace}\nolimits}

\def\build#1_#2^#3{\mathrel{
\mathop{\kern 0pt#1}\limits_{#2}^{#3}}}

\def\@nbibitem#1{\noindent \hangindent=2pc \hangafter=1
\refstepcounter{enumi}\hbox to 2pc{\arabic{enumi}.\hfil}%
\immediate\write\@auxout{\string\bibcite{#1}{\arabic{enumi}}}}
\def\numbibliography{%
\section*{REFERENCES}%
\bgroup\footnotesize
\setcounter{enumi}{0}%
\def\newblock{\hskip .11em plus.33em minus.07em}%
\let\bibitem\@nbibitem}
\def\endnumbibliography{\par\egroup}
\makeatother
\def\ra{\rightarrow}
\def\longra{\longrightarrow}
\def\mpo{\mapsto}

\font\tenbb=msbm10
\font\sevenbb=msbm7
\font\fivebb=msbm5
\newfam\bbfam
\textfont\bbfam=\tenbb \scriptfont\bbfam=\sevenbb
\scriptscriptfont\bbfam=\fivebb
\def\bb{\fam\bbfam}

\def\Cb{{\bb C}}
\def\Qb{{\bb Q}}
\def\Rb{{\bb R}}
\def\Sb{{\bb S}}
\def\Tb{{\bb T}}
\def\Zb{{\bb Z}}

\def\Ac{{\cal A}}
\def\Ec{{\cal E}}
\def\Hc{{\cal H}}

\def\vp{\varphi}
\def\Om{\Omega}
\def\om{\omega}
\def\t{\theta}
\def\d{\delta}
\def\ve{\varepsilon}
\def\s{\sigma}
\def\G{\Gamma}
\def\Lb{\Lambda}
\def\g{\gamma}

\def\wdg{\wedge}
\def\fl{\forall}
\def\op{\oplus}
\def\sbs{\subset}

\def\build#1_#2^#3{\mathrel{
\mathop{\kern 0pt#1}\limits_{#2}^{#3}}}

\begin{document}

\title{\centerline{ NONCOMMUTATIVE GEOMETRY}\centerline{ YEAR 2000}}
% TITLE IN CAPITAL LETTERS

\author{\bf Alain CONNES,\refnote{1}}

\affiliation{\affnote{1}  Coll\`ege de France,
3, rue Ulm,
75005 PARIS\\
and\\
I.H.E.S.,
35, route de Chartres,
91440 BURES-sur-YVETTE
}

\vspace{1cm}

\begin{abstract}
Our geometric concepts evolved first through the discovery of NonEuclidean
geometry. The discovery of quantum mechanics in the form of the noncommuting
coordinates on the phase space of atomic systems entails an equally drastic evolution.
We describe a basic construction which extends the familiar duality between ordinary spaces
and commutative algebras to a duality between Quotient spaces and Noncommutative algebras.
The basic tools of the theory, K-theory, Cyclic cohomology, Morita equivalence, Operator theoretic index theorems,
 Hopf algebra symmetry are reviewed. They cover the global aspects of noncommutative spaces, such as the 
transformation $\theta \rightarrow 1/\theta$  for the noncommutative torus $\Tb_{\theta}^2$ which are unseen
in perturbative 
expansions in $\theta$ such as star or Moyal products. We discuss the foundational problem of 
"what is a manifold in NCG" and explain the fundamental role of Poincare duality in K-homology
which is the basic reason for the spectral point of view. This leads us, when specializing to 4-geometries
to a universal algebra called the "Instanton algebra". We describe our joint work with G. Landi which
gives noncommutative spheres $S_{\theta}^4$ from representations of the Instanton algebra.
We show that any compact Riemannian spin manifold whose isometry group has rank $r \geq 2$ admits
isospectral deformations to noncommutative geometries.
We give a survey of several recent developments. First our joint work with H. Moscovici on 
the transverse geometry of foliations which yields a diffeomorphism invariant (rather than the usual covariant
one) geometry on the bundle of metrics on a manifold and a natural extension of cyclic cohomology
to Hopf algebras. Second,  our joint work with D. Kreimer on renormalization and the Riemann-Hilbert
problem. Finally we describe the spectral realization of zeros of zeta and L-functions from the noncommutative space 
of Adele classes on a global field and its relation with the Arthur-Selberg
trace formula in the Langlands program. We end with a tentalizing connection
between the renormalization group and the missing Galois theory at  Archimedian
places.

\end{abstract}

\newpage

\section{Introduction} 
 There are two fundamental sources of 
`bare' facts for the mathematician. These are, on the one hand
 the physical world which is the source of \emph{geometry}, and
 on the other hand the arithmetic of numbers which is the source of \emph{number theory}. 
Any theory concerning either of these subjects can be tested 
 by performing experiments either in the physical world or with numbers. 
That is, there are some real things out there 
 to which we can confront our understanding.

 \noindent If one looks back at the 23 problems of Hilbert then one finds that, 
fortunately, the twentieth century saw very important discoveries which 
nobody could have foreseen by 1900. Two of them (of course by no means 
the only discoveries) involve Hilbert space in a crucial way and will be of particular importance 
for this talk:  The first one is quantum mechanics, 
and the second, equally important in a sense,
 is the extension of class field theory
 to the non-abelian case, 
thanks to the Langlands program.

\noindent  In this lecture  I'll take both of these discoveries as a pretext and point towards the extension
 of our familiar geometrical concepts beyond the classical, commutative case.
 My aim is to discuss the foundation of noncommutative geometry.

  \section{Geometry}    Before I do that, let me remind you, using a simple example, of the power 
of abstraction in mathematics.
 Around 1800, Mathematicians wondered whether it is true that Euclid's fifth axiom
 is actually superfluous. For instance Legendre proved that if you have one triangle
 whose internal angles sum to $\pi$ then that is enough to guarantee ordinary Euclidean geometry.
  However, as we all know Euclid's fifth axiom
 is not superfluous and NonEuclidean Geometry gives a counter-example. 
The simplest model of NonEuclidean Geometry is probably the Klein model.
The points of the geometric space X are the points inside an ellipse, 

$$
\hbox{
\psfig{figure=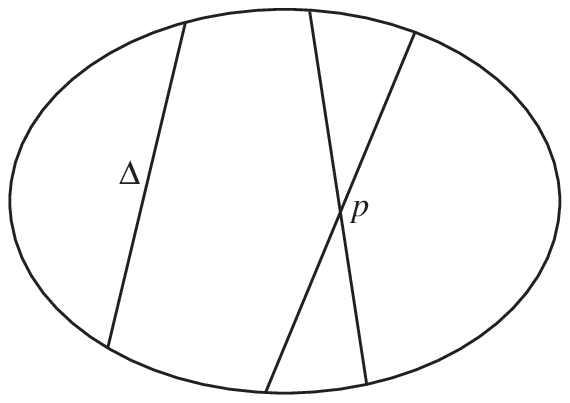}
}
$$
 
 \noindent The lines are the intersections of the ordinary Euclidean lines with X.
 If you take a point $p$, outside the line $\Delta$ then there are distinct lines 
which don't meet $\Delta$ (\emph{i.e.} are parallel to $\Delta$) but meet each other at $p$.

 \noindent  At first this was considered as an esoteric example and Gauss didn't publish his discovery,
 but after some time it became clear that  rather than just being a strange counter-example,
 it was  something with remarkable beauty and power. 
The question then  became ``what is the source of this beauty and power?''
 Often in  mathematics, understanding comes from generalisation, instead of 
 considering the object \emph{per se} what one tries to find are  the concepts 
which embody the power of the object.  

\noindent  A first  generalisation is the \emph{Erlangen} program of Klein 
and the  theory of Lie groups which attributes the beauty of this example 
 to its symmetries, namely the group of projective transformations of the plane 
which preserve  the ellipse. 

 \noindent  The second conceptual generalisation is Riemannian geometry
as explained in Riemann's inaugural  lecture (\cite{[R]})
 in which he reflected on the hypotheses of geometry and 
 introduced two key notions: the concepts of \emph{manifold} and  \emph{line element}.

 \noindent  By a manifold Riemann meant `any space you can think of whose points can vary continuously'. 
For example, a manifold could be a continuous collection of colours, 
the parameter space for some mechanical system or, of course, space.
 In his lecture Riemann explained that it is possible, essentially proceeding by induction, to label the points of such a space
 by a finite collection of real numbers. 

\noindent In Riemannian geometry the distance between two points $x$ and $y$ is given by the following \emph{ansatz}:
\begin{equation}
d(x,y) = \hbox{Inf} \{ \int_\g ds \, | \g \
\text{is a path between $x$ and $y$} \} \label{eq:(1.1)}
\end{equation}

\noindent  Expanding $d(x, y)$ near the diagonal, after raising it to an even power 
to ensure smoothness gives a local formula for ds. The first case 
he considered was the quadratic case (although he explicitly mentionned the quartic case).
 From the Taylor expansion he obtained, in the quadratic case, the well-known formula for the metric,
  \begin{equation}
ds^2 = g_{\mu \nu} \, dx^{\mu} \, dx^{\nu} \, . \label{eq:(1.2)}
\end{equation}

 \noindent Riemann's concept of geometry differs greatly from that of Klein because
 Klein's formulation is based on the idea of rigid motions whereas in Riemannian geometry
 rigid motions are no longer possible because of the variability of the curvature 
and the extraordinary freedom in the choice of the components $g_{\mu\nu}$.

 \noindent The basic notions of ordinary geometry do make sense, for example 
a straight line is given by the geodesic equation, 
\begin{equation}
{d^2 \, x^{\mu} \over dt^2} = -{1\over 2} \, g^{\mu \a}
(g_{\a \nu ,\rho} + g_{\a \rho ,\nu} - g_{\nu \rho
,\a}) \,
{dx^{\nu} \over dt} \ {dx^{\rho} \over dt}              \label{eq:(1.3)}
\end{equation}
but what really vindicated the point of view of Riemann, 
with respect to that of Klein, was another major 
discovery of the twentieth century, General Relativity.

\noindent  One can get a glimpse of this from the following simple fact. 
If we take the Minkowski metric and perturb it to $dx^2 + dy^2 + dz^2 - (1 +
2 V(x,y,z)) dt^2$
using the Newtonian potential $V(x,y,z)$, then the geodesic equation 
can be re-written in the obvious approximation to obtain Newton's law of motion. 
 This makes clear that the variability of the $g_{\mu\nu}$ is precisely 
necessary in order to get a good geometric model of the physical universe. 

\noindent It is  interesting to note that Riemann  was well aware of 
the limits of his own point of view as is  clearly expressed in the last page 
of his inaugural lecture; (\cite{[R]})

\noindent   "Questions about the immeasurably large are idle questions 
for the explanation of Nature. But the situation is quite 
different with questions about the immeasurably small. 
Upon the exactness with which we pursue phenomenon into 
the infinitely small, does our knowledge of their causal 
connections essentially depend. The progress of recent 
centuries in understanding the mechanisms of Nature 
depends almost entirely on the exactness of construction 
which has become possible through the invention of the 
analysis of the infinite and through the simple 
principles discovered by Archimedes, Galileo and Newton, 
which modern physics makes use of. By contrast, in the 
natural sciences where the simple principles for such 
constructions are still lacking, to discover causal 
connections one pursues phenomenon into the spatially 
small, just so far as the microscope permits. Questions 
about the metric relations of Space in the immeasurably 
small are thus not idle ones.

If one assumes that bodies exist independently of 
position, then the curvature is everywhere constant, and 
it then follows from astronomical measurements that it 
cannot be different from zero; or at any rate its 
reciprocal must be an area in comparison with which the 
range of our telescopes can be neglected. But if such an 
independence of bodies from position does not exist, then 
one cannot draw conclusions about metric relations in the 
infinitely small from those in the large; at every point 
the curvature can have arbitrary values in three 
directions, provided only that the total curvature of 
every measurable portion of Space is not perceptibly 
different from zero. Still more complicated relations can 
occur if the line element cannot be represented, as was 
presupposed, by the square root of a differential 
expression of the second degree. Now it seems that the 
empirical notions on which the metric determinations of 
Space are based, the concept of a solid body and that of 
a light ray, lose their validity in the infinitely small; 
it is therefore quite definitely conceivable that the 
metric relations of Space in the infinitely small do not 
conform to the hypotheses of geometry; and in fact one 
ought to assume this as soon as it permits a simpler way 
of explaining phenomena.

The question of the validity of the hypotheses of 
geometry in the infinitely small is connected with the 
question of the basis for the metric relations of space. 
In connection with this question, which may indeed still 
be ranked as part of the study of Space, the above remark 
is applicable, that in a discrete manifold the principle 
of metric relations is already contained in the concept 
of the manifold, but in a continuous one it must come 
from something else. Therefore, either the reality 
underlying Space must form a discrete manifold, or the 
basis for the metric relations must be sought outside it, 
in binding forces acting upon it.

An answer to these questions can be found only by 
starting from that conception of phenomena which has 
hitherto been approved by experience, for which Newton 
laid the foundation, and gradually modifying it under the 
compulsion of facts which cannot be explained by it. 
Investigations like the one just made, which begin from 
general concepts, can serve only to insure that this work 
is not hindered by too restricted concepts, and that 
progress in comprehending the connection of things is not 
obstructed by traditional prejudices.

This leads us away into the domain of another science, 
the realm of physics, into which the nature of the 
present occasion does not allow us to enter".

 \bigskip 
 \section{Quantum mechanics} 

   In fact quantum mechanics showed that indeed the parameter space, or phase space
of the mechanical system given by a single atom fails to be a manifold.
\noindent It is important to convince oneself of this fact and to understand
that this conclusion is indeed dictated by the experimental findings of spectroscopy.
The information we get from the light coming from distant stars is of spectral 
nature, the spectral lines are absorption or emission lines 

$$
\hbox{
\psfig{figure=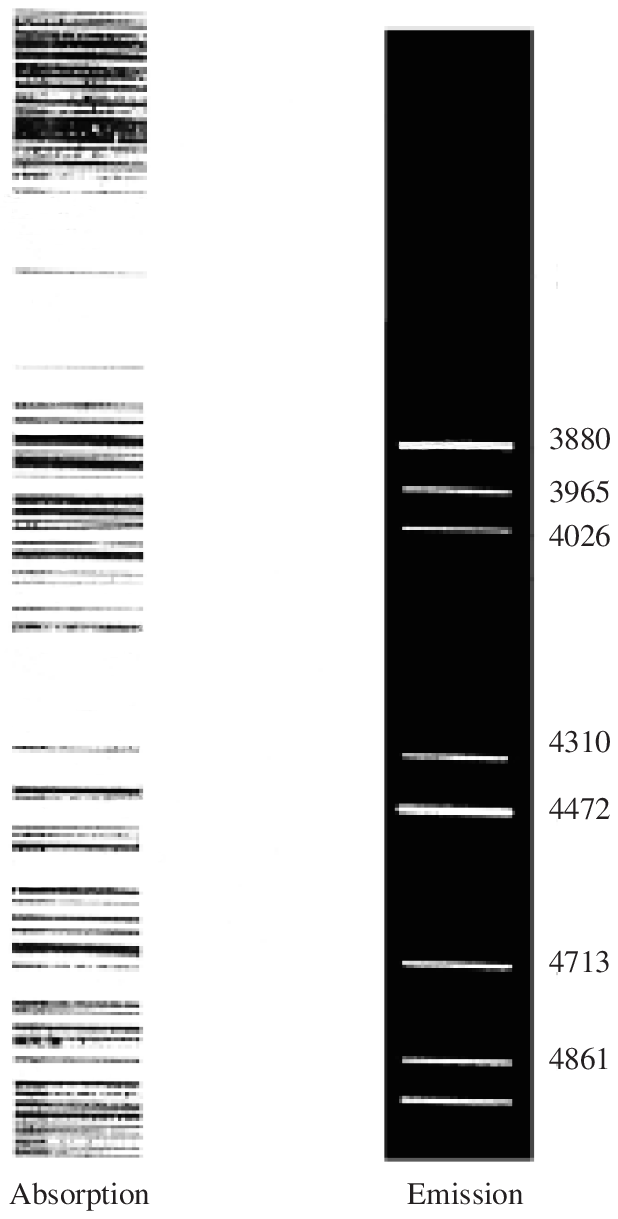}
}
$$

 \noindent One can infer from this spectral information the chemical composition of the star
since the simple elements have recognisable spectra. 
These spectra obey experimentally discovered laws, the most notable being the 
Ritz-Rydberg combination principle.
The principle can be stated as follows; spectral lines are indexed by pairs of objects.
These objects  could be numbers, greek letters, or any kind of labels. The statement of the 
principle then is that certain pairs of spectral lines, when expressed in terms of frequencies,
do add up to give another line in the spectrum.
Moreover, this happens 
 precisely when the labels are of the form $i,j$ and $j,k$.

\noindent  What Heisenberg understood, by analogy with the classical treatment
of the interaction of a mechanical system with the electromagnetic field,
 is that this Ritz-Rydberg combination principle actually dictates an algebraic formula for the  
product of any two observable physical quantities attached to the atomic system.

$$
\hbox{
\psfig{figure=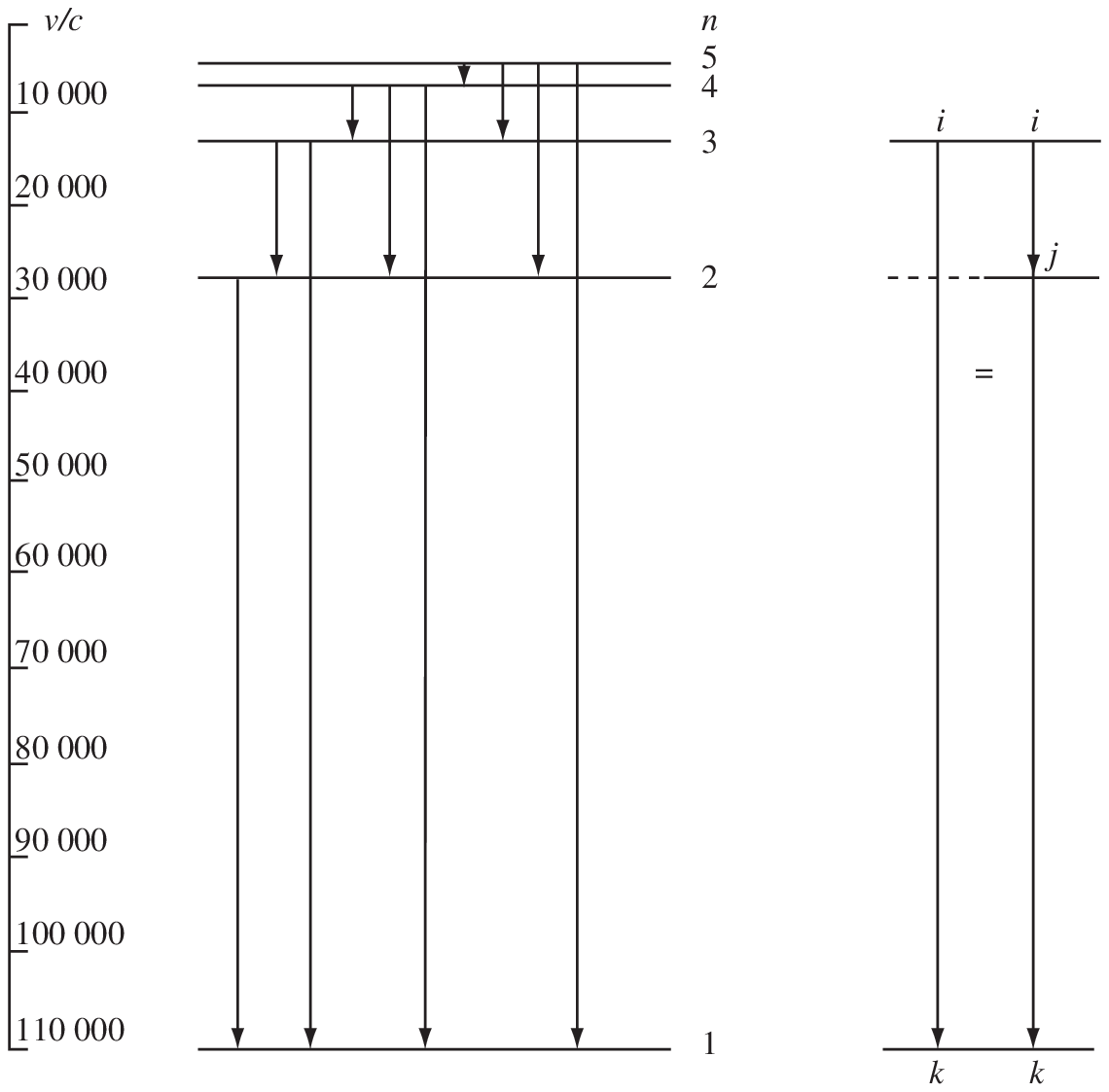}
}
$$

\noindent   Heisenberg wrote down the formula for the product of two observables;
    \begin{equation}
(A \, B)_{(i,k)} \, =  A_{(i,j)} \, B_{(j,k)}             \label{eq:(1.4)}
\end{equation}
\noindent and he noticed of course that this algebra he had found is no longer commutative,
   \begin{equation}
A \, B   \not= B\,A         \label{eq:(1.5)}
\end{equation}

\noindent    Now Heisenberg didn't know about matrices, he just worked it 
 out, but he was told later by Born, Jordan and Dirac that the algebra he had worked out
 was known to mathematicians as the algebra of matrices.

\noindent  Physicists   often tell jokes such as:      
A physicist walks down the main street of a strange town looking 
   for a laundrette. He sees a shop with signs in the window  
  saying `bakery' `grocers' `laundrette', so he enters. 
However,    the shop is owned by a mathematician and
 when the physicist    asks ``when will the washing be ready?'' 
the mathematician    replies ``we don't clean clothes, we just sell signs!''. 

 \noindent  In the case of Heisenberg and also that of Einstein who was  helped out by Riemann,
 this was no joke.

 \noindent However, soon after  Heisenberg's discovery, Schrodinger came up
 with his equation  so physicists happily returned to the study of partial differential  equations,
 and the message of Heisenberg was buried to a great extent.
\noindent   Most of my work has been an attempt to take this discovery of  Heisenberg seriously.
 On reflection, this discovery actually clearly displays the limitation of Riemann's formulation of geometry.
 If we look at the phase space of an atomic system and follow Riemann's procedure to parametrize
its points by finitely many real numbers, we first split the  manifold into the levels 
on which some particular  function is constant, but we then need to iterate this process
and apply it to the level hypersurfaces.
  However, according to Heisenberg this doesn't work because as  
soon as we make the first measurement, we alter the situation drastically.
 The right way to think about this new phenomenon is to think in  terms of 
 a new kind of space in which the coordinates do not commute. 

\noindent The starting point of noncommutative geometry is to take this new notion of space seriously.

  \section{Noncommutative geometry}   The basis of noncommutative geometry is twofold. 

\noindent On the one hand there is a wealth of examples of 
spaces whose coordinate algebra is no longer commutative 
but which have obvious relevance in physics or mathematics. The first examples
 came, as we saw above, from phase space in quantum mechanics but there 
are many others, such as the leaf spaces of foliations, the 
duals of nonabelian groups, the space of Penrose
tilings, the Brillouin zone in solid state physics, the noncommutative tori which 
appear
naturally in string theory and in
M-theory compactification, and the Adele class space
which as we shall see below provides a natural spectral realisation of zeros of zeta 
functions. Finally 
various recent models
 of space-time itself are interesting examples of noncommutative spaces.

\noindent On the other hand the stretching of geometric thinking 
imposed by passing to noncommutative spaces forces one
to rethink about most of our familiar notions.
The difficulty is not to add arbitrarily the adjective
quantum to our geometric words but to
develop far reaching extensions of classical concepts, ranging from the simplest 
which is measure theory, 
to the most sophisticated which is geometry itself.

\noindent Let us first discuss in greater detail the general principles that allow 
to construct
huge classes of  such spaces, it is a vital ingredient indeed 
 since there is no way to build a satisfactory theory without being able to test it on a large 
variety of examples. We have two principles which allow us to  construct examples.

\noindent The first is deformation  from the commutative to the noncommutative 
which allows to explore the neighborhood of the comutative world.

\noindent The second is a new and very  important mathematical principle; the quotient operation. 
Most of  the spaces we are concerned with are not defined by naming every 
one of their points, but by giving a much bigger set and dividing  it by an equivalence relation. 

\noindent   It turns out that there are two ways of extending the
 geometric -  algebraic duality 
\begin{equation}
\hbox{Space}  \leftrightarrow \hbox{Commutative algebra} \label{eq:(0.1)}
\end{equation}
between a space $X$ and the algebra of   functions on that space, when you want to identify two points $a$ and $b$. 
 The first way which gives the usual algebra of functions associated to the quotient is
 to restrict oneself to functions which have the same  value at the two points.

\begin{equation}
{\cal A} = \{ f ; f(a) = f(b) \} \, . \label{eq:(0.2)}
\end{equation}
  The second way is to keep  the two points $a$ and $b$, but to allow them
 to  `speak' to each other by using matrices with off-diagonal  elements.
 It consists, instead of taking the subalgebra
given by 4-\ref{eq:(0.2)}, to adjoin to the algebra of functions on $\{ a,b \}$
the identification of $a$ with $b$. The obtained algebra is the algebra of
two by two matrices
\begin{equation}
{\cal B} = \left\{ f = \left[ \matrix{
f_{aa} &f_{ab} \cr f_{ba} &f_{bb} \cr} \right] \right\} \, . \label{eq:(0.3)}
\end{equation}
 When one computes the spectrum of this algebra it turns out that  it is composed 
of only one point, so the two points $a$ and $b$ have been  identified.
 As we shall see this second method is very powerful and allows 
one to  construct thousands of very interesting examples.
It allows to refine the above
duality of algebraic geometry to,
\begin{equation}
\hbox{Quotient-Space}  \leftrightarrow \hbox{Noncommutative algebra} \label{eq:(0.4)}
\end{equation}
in the situation where the space one is contemplating is obtained by the
operation  of quotient.

 \noindent  At first sight it might seem that, as far as the general theory is concerned, passing from the commutative to the noncommutative
situation would just be a matter of cleverly rewriting in algebraic terms our familiar geometric 
notions without using commutativity anywhere. If noncommutative geometry was just that it would be boring indeed. 
Fortunately, even at the coarsest level which is measure theory, it became clear at the beginning of
the seventies that the noncommutative world is full of beautiful totally unexpected facts which have 
no commutative counterpart whatsoever.
\noindent The prototype of such facts is the following
\begin{equation}
\hbox{ Noncommutative measure spaces evolve with time!}   \label{eq:(0.5)}
\end{equation}
 
\noindent  In other words there is a `god-given' one parameter group of
 automorphisms of the algebra $M$ of measurable coordinates. 
\noindent It is given by the group homomorphism, (\cite{Co_2})
\begin{equation}
\delta : \Rb \rightarrow {\rm Out} (M) = {\rm Aut} (M) / {\rm Int} (M) \, 
\label{eq:(0.6)}
\end{equation}
from the additive group $\Rb$ to the group of automorphism classes of $M$ modulo 
inner 
automorphisms.

\noindent I discovered this fact in 1972 when working on the Tomita-Takesaki theory (\cite{T})
and it convinced me that there are amazing  
 features of noncommutative spaces which have no counterpart in the static commutative  case.

  \section{A basic example} 

\noindent Let us start with a prototype example of quotient space in which the distinction between 
the quotient operations (4-\ref{eq:(0.2)})
and (4-\ref{eq:(0.3)}) appears clearly, and which played a key role in 1980 at the early stage of the 
theory (\cite{C3}). 
 This
example is the following: consider the 2-torus
\begin{equation}
M = \Rb^2 / \Zb^2 \, . \label{eq:(0.5)}
\end{equation}
The space $X$ which we contemplate is the space of solutions of the
differential equation,
\begin{equation}
dx = \theta dy \qquad x,y \in \Rb/\Zb \label{eq:(0.6)}
\end{equation}
where $\theta \in ]0,1[$ is a fixed irrational number.

$$
\hbox{
\psfig{figure=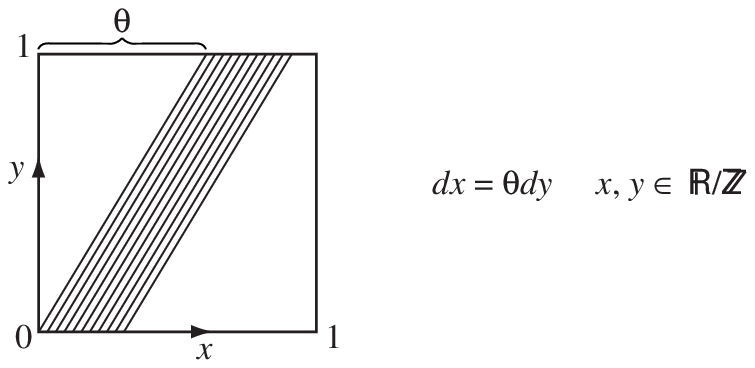}
}
$$

\medskip

\noindent Thus the space we are 
interested in here is just the space of leaves of the foliation defined by the 
differential 
equation 5-\ref{eq:(0.6)}.
We can label such a leaf by a point of the transversal given by $y=0$ which is a 
circle $S^1 
= \Rb /\Zb$, 
 but clearly two 
points of the transversal which differ by an integer multiple of $\theta$ give rise 
to the 
same leaf. Thus 
\begin{equation}
X = S^1 / \theta \Zb \label{eq:(0.7)}
\end{equation}
{\it i.e.} $X$ is the quotient of $S^1$ by the equivalence relation which
identifies any two points on the orbits of the irrational rotation
\begin{equation}
R_{\theta} x = x + \theta \quad {\rm mod 1} \, . \label{eq:(0.8)}
\end{equation}
When we deal with $S^1$ as a space in the various categories (smooth,
topological, measurable) it is perfectly described by the corresponding
algebra of functions,
\begin{equation}
C^{\infty} (S^1) \subset C(S^1) \subset L^{\infty} (S^1) \, . \label{eq:(0.9)}
\end{equation}
When one applies the naive operation (4-\ref{eq:(0.2)}) to pass to the quotient, one
finds, irrespective of which category one works with, the trivial answer
\begin{equation}
\Ac = \Cb \, . \label{eq:(0.10)}
\end{equation}
The operation (4-\ref{eq:(0.3)}) however gives very interesting algebras, by no
means reduced to $\Cb$. Elements of the 
algebra $\Bc$ associated to the transversal $ S^1$ by the operation (4-\ref{eq:(0.3)})
 are just matrices $a(i,j)$ where the 
indices $(i,j)$ are arbitrary pairs of elements $i,j$ of $ S^1$ which belong to the 
same leaf, i.e. give the same element of $X$.
 The algebraic rules are the same as for ordinary matrices. In the above situation 
since the equivalence is given by a group action, the
construction coincides with the crossed product familiar to algebraist from the theory of
central simple algebras.

\noindent An element of $\Bc$ is given by a power series
\begin{equation}
b = \sum_{n \in \Zb} b_n U^n \label{eq:(0.11)}
\end{equation}
where each $b_n$ is an element of the algebra 5-\ref{eq:(0.9)}, while the
multiplication rule is given by
\begin{equation}
U h U^{-1} = h \circ R_{\theta}^{-1} \, . \label{eq:(0.12)}
\end{equation}
Now the algebra 5-\ref{eq:(0.9)} is generated by the function $V$ on $S^1$,
\begin{equation}
V(\alpha) = {\rm exp} (2\pi i \alpha) \qquad \alpha \in S^1 \label{eq:(0.13)}
\end{equation}
and it follows that $\Bc$ admits the generating system $(U,V)$ with
presentation given by the relation
\begin{equation}
VU = \lambda \, U V \qquad \lambda = {\rm exp} 2\pi i \theta \, .
\label{eq:(0.14)}
\end{equation}
Thus, if for instance we work in the smooth category a generic element $b$ of
$\Bc$ is given by a power series
\begin{equation}
b = \sum_{\Zb^2} b_{nm} U^n V^m \ , \ b \in {\cal S} (\Zb^2) \label{eq:(0.15)}
\end{equation}
where ${\cal S} (\Zb^2)$ is the Schwartz space of sequences of rapid decay on
$\Zb^2$.

\noindent  This algebra is by no means trivial and has a very rich
and interesting algebraic structure. It is (canonically up to Morita equivalence) 
associated 
to the
foliation 5-\ref{eq:(0.6)} and the interplay between the geometry of the
foliation and the algebraic structure of $\Bc$ begins by noticing that to a
{\it closed transversal} $T$ of the foliation corresponds canonically a {\it
finite projective module} over $\Bc$.
 Elements 
of the module associated to the transversal $T$ are rectangular matrices, 
$ \xi(i,j)$ where  $(i,j) \in T \times S^1$ while $i$ and $j$ belong to the 
same leaf, i.e. give the same element of $X$. The right action of $a(i,j) \in \Bc$ is 
by matrix multiplication.

\noindent  From the transversal $x=0$, one obtains the following right module over 
$\Bc$. The 
underlying linear
space is the usual Schwartz space,
\begin{equation}
{\cal S} (\Rb) = \{ \xi , \, \xi (s) \in \Cb \quad \forall s \in \Rb \}
\label{eq:(0.16)}
\end{equation}
of smooth functions on the real line all of whose derivatives are of rapid decay.

\noindent  The right module structure is given by the action of the generators $U,V$
\begin{equation}
(\xi U) (s) = \xi (s+\theta) \ , \ (\xi V) (s) = e^{2\pi is} \xi (s) \quad
\forall s \in \Rb \, . \label{eq:(0.17)}
\end{equation}
One of course checks the relation 5-\ref{eq:(0.14)}, and it is a beautiful fact
that as a right module over $\Bc$ the space ${\cal S} (\Rb)$ is {\it finitely
generated} and {\it projective} ({\it i.e.} complements to a free module). It
follows that it has the correct algebraic atributes to deserve the name of
``noncommutative vector bundle'' according to the dictionary,
$$
\matrix{
\hbox{Space} &\hbox{Algebra} \cr
\cr
\hbox{Vector bundle} &\hbox{ \, \,Finite projective module.} \cr
}
$$
The concrete description of the general
finite projective modules over ${\cal A}_{\t}$ is obtained by combining the 
results of \cite{[26], C3, Rieffel}. 
They are classified up to isomorphism by a pair 
of integers $(p,q)$ such that $p+q \t \geq 0$ and  the corresponding modules ${\cal 
H}_{p,q}^{\t}$
are obtained by the above construction from the 
transversals given by closed geodesics of the torus $M$.

\noindent The algebraic counterpart of a vector bundle is its space of smooth 
sections
$C^{\ify} (X,E)$ and one can in particular compute its dimension by computing
the trace of the identity endomorphism of $E$. If one applies this method in
the above noncommutative example, one finds
\begin{equation}
{\rm dim}_{\Bc} ({\cal S}) = \theta \, . \label{eq:(0.18)}
\end{equation}
The appearance of non integral dimension is very exciting and displays a
basic feature of von Neumann algebras of type II. The dimension of a vector
bundle is the only invariant that remains when one looks from the measure
theoretic point of view ({\it i.e.} when one takes the third algebra in
5-\ref{eq:(0.9)}). The von Neumann algebra which describes the quotient space $X$
 from the measure theoretic point of view is the crossed
product,
\begin{equation}
R = L^{\infty} (S^1) \semi_{R_{\theta}} \Zb \label{eq:(0.19)}
\end{equation}
and is the well known hyperfinite factor of type II$_1$. In particular the
classification of finite projective modules $\Ec$ over $R$ is given by a positive
real number, the Murray and von Neumann {\it dimension},
\begin{equation}
{\rm dim}_R (\Ec) \in \Rb_+ \, . \label{eq:(0.20)}
\end{equation}
The next surprise is that even though the {\it dimension} of the above module
is irrational, when we compute the analogue of the first Chern class, {\it i.e.} of
the integral of the curvature of the vector bundle, we obtain an integer.
Indeed the two commuting vector fields which span the tangent space for an
ordinary (commutative) 2-torus correspond algebraically to two commuting
derivations of the algebra of smooth functions. These derivations continue to
make sense when the generators $U$ and $V$ of $C^{\infty} (\Tb^2)$ no longer
commute but satisfy 5-\ref{eq:(0.14)} so that they generate $\Bc = C^{\infty} 
(\Tb_{\theta}^2)$. They 
are given by the same formulas as
in the commutative case,
\begin{equation}
\delta_1 = 2\pi i U \, \frac{\partial}{\partial U} \ , \delta_2 = 2\pi i V \,
\frac{\partial}{\partial V} \label{eq:(0.21)}
\end{equation}
so that $\delta_1 \left( \sum b_{nm} U^n V^m \right) = 2\pi i \sum n b_{nm}
U^n V^m$ and similarly for $\delta_2$. One still has of course
\begin{equation}
\delta_1 \delta_2 = \delta_2 \delta_1 \label{eq:(0.22)}
\end{equation}
and the $\delta_j$ are still derivations of the algebra $\Bc= C^{\infty} 
(\Tb_{\theta}^2)$,
\begin{equation}
\delta_j (bb') = \delta_j (b) b' + b \delta_j (b') \quad \forall b , b' \in
\Bc \, . \label{eq:(0.23)}
\end{equation}
The analogues of the notions of connection and curvature of vector bundles
are straightforward to obtain (\cite{C3}) since a connection is just given by the
associated covariant differentiation $\nabla$ on the space of smooth
sections. Thus here it is given by a pair of linear operators,
\begin{equation}
\nabla_j : {\cal S} (\Rb) \rightarrow {\cal S} (\Rb) \label{eq:(0.24)}
\end{equation}
such that
\begin{equation}
\nabla_j (\xi b) = (\nabla_j \xi)b + \xi \delta_j (b) \quad \forall \xi \in
{\cal S} \ , b \in \Bc \, . \label{eq:(0.25)}
\end{equation}
One checks that, as in the usual case, the trace of the curvature $\Omega =
\nabla_1 \nabla_2 - \nabla_2 \nabla_1$, is independent of the choice of the
connection. Now the remarkable fact here is that (up to the correct powers of
$2\pi i$) the total curvature of ${\cal S}$ is an integer. In fact for the
following choice of connection the curvature $\Omega$ is constant, equal to
$\frac{1}{\theta}$ so that the irrational number $\theta$ disappears in the
total curvature, $\theta \times \frac{1}{\theta}$
\begin{equation}
(\nabla_1 \xi) (s) = - \frac{2\pi i s}{\theta} \, \xi (s) \ \,\,\,(\nabla_2 \xi)
(s) = \xi' (s) \, . \label{eq:(0.26)}
\end{equation}
With this integrality, one could get the wrong impression that the algebra
$\Bc= C^{\infty} (\Tb_{\theta}^2)$ looks very similar to the algebra $C^{\infty} 
(\Tb^2)$ of 
smooth
functions on the 2-torus. A striking difference is obtained by looking at the 
range of Morse functions. The range of a Morse function on $\Tb^2$ is of course a connected interval.
 For the above noncommutative torus $\Tb_{\theta}^2$ the range of a Morse function is the spectrum of a real valued
function such as
\begin{equation}
h = U + U^* + \mu (V+V^*) \label{eq:(0.27)}
\end{equation}
and it can be a Cantor set, {\it i.e.} have infinitely many disconnected pieces. This
shows that the one dimensional pictures of our space $\Tb_{\theta}^2$ are
truly different from what they are in the commutative case. The above noncommutative
torus $\Tb_{\theta}^2$ is the simplest example of noncommutative manifold, it arises 
naturally
not only from foliations but also from the Brillouin zone in the Quantum Hall
effect as understood by J. Bellissard, and in M-theory as we shall see next.
In the Quantum Hall
effect, the above integrality of the total curvature corresponds to the observed integrality of the 
Hall  conductivity 

$$
\hbox{
\psfig{figure=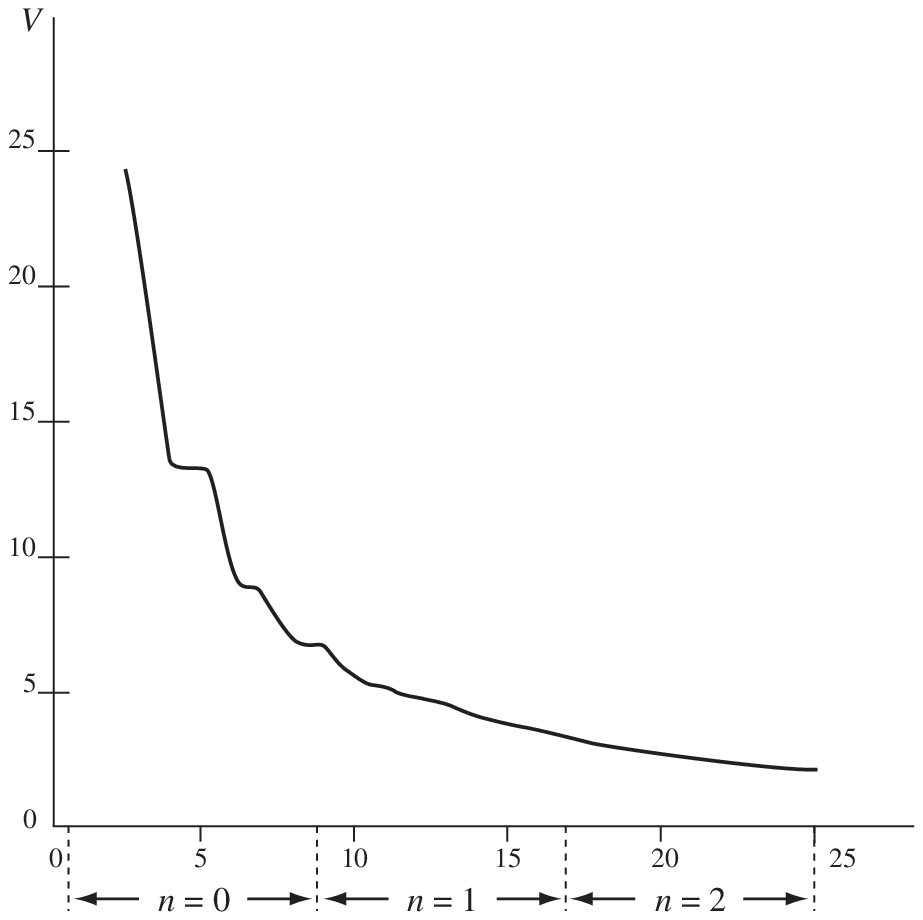}
}
$$

\noindent The analogue of the Yang-Mills action functional and the classification of 
Yang-Mills connections on the 
noncommutative tori was developped in \cite {Co-R}, with the primary goal of 
finding 
a "manifold shadow" for these noncommutative spaces.
These moduli spaces turned out indeed to fit this purpose perfectly, allowing for 
instance to find the usual Riemannian space of gauge equivalence classes of 
Yang-Mills connections as an invariant of the noncommutative metric. 

\noindent The next 
surprise came from the natural occurence (as an unexpected guest) of 
both the noncommutative tori and the components of the Yang-Mills connections in 
the 
classification of the BPS states in M-theory \cite{CDS}.

\noindent  In the matrix formulation of M-theory the basic equations to obtain periodicity 
of two of the basic coordinates $X_i$ turn out to be the 
following, 
\begin{equation}
U_i X_j U_i^{-1} =X_j + a \delta_i^j , i= 1,2  \label{eq:(0.136)}
\end{equation}
where the $U_i$ are unitary gauge transformations.

 \noindent  The multiplicative commutator $U_1 U_2 U_1^{-1}U_2^{-1}$ is then central
and in the irreducible case its scalar value $ \lambda =\exp 2\pi i \t $ brings in 
the algebra of coordinates on the noncommutative torus.
The $X_j$ are then the components of the Yang-Mills connections. It is quite 
remarkable that the same picture emerged from the
other information one has about M-theory concerning its relation with 11 
dimensional 
supergravity and that string theory dualities could be interpreted
using Morita equivalence. The latter relates the 
values of $\theta$ on an orbit of $SL (2,\Zb)$ and simply illustrates that the leaf-space of the 
original foliation is independent of which transversal is used to parametrize it. This
 type of relation between for instance $\theta$ and $1/{\theta}$ would be invisible in a purely deformation theoretic 
perturbative 
expansion like the one given by the Moyal product.

\noindent   Nekrasov and Schwarz \cite{N-S} showed that 
Yang-Mills gauge theory on
 noncommutative $\Rb^4$ gives a conceptual understanding of the nonzero B-field 
desingularization of the 
moduli space of instantons obtained by perturbing the ADHM equations.

\noindent In \cite{Witten}, Seiberg and Witten exhibited the unexpected relation 
between the 
standard gauge theory and the noncommutative one, and clarified 
the limit in which the entire string dynamics is described by a gauge theory
on a noncommutative space.

\smallskip

\noindent  One should understand from the very start
that foliations provide an inexhaustible source of interesting examples of noncommutative spaces. 
\noindent In the above example of $\Tb_{\theta}^2$ we could make use of
the special vector fields on the torus in order to obtain the analogues of
elementary notions of differential geometry. It is quite important 
 to develop the general
theory independently of these special features and this is what we shall do in section VII.
We shall start by the noncommutative analogues of topology and vector bundles which are necessary preliminary steps.

\bigskip
\section{Topology}\noindent The development of the topological ideas was prompted by the work of 
Israel 
Gel'fand, whose C* algebras give the required framework for noncommutative 
topology. 
The two 
main driving forces were the Novikov conjecture on homotopy invariance of higher 
signatures 
of ordinary manifolds as well as the Atiyah-Singer Index theorem. It has led, 
through the 
work of  Atiyah, Singer, Brown, Douglas, Fillmore, Miscenko and Kasparov  
\cite{[AT]} \cite{Singer} \cite{B-D-F} \cite{[21]} \cite{[18]} to
the  recognition 
that not only the 
Atiyah-Hirzebruch K-theory but more importantly the dual K-homology admit Hilbert 
space 
techniques and functional analysis as their natural framework. The cycles in the 
K-homology 
group $K_*(X)$ of a compact space X are indeed given by Fredholm representations 
of 
the C* 
algebra A of continuous functions on X. The central tool is the Kasparov bivariant 
K-theory. 
A basic example of C* algebra to which the theory applies is the group ring of a 
discrete 
group and this makes it clear that restricting oneself to commutative algebras is an undesirable 
assumption.

\noindent For a $C^*$ algebra $A$, let $K_0 (A)$, $K_1 (A)$ be its $K$ theory 
groups. 
Thus $K_0 (A)$ is the algebraic $K_0$ theory of the ring $A$ and $K_1 (A)$ 
is the algebraic $K_0$ theory of the ring $A \ot C_0 (\Rb) = C_0 (\Rb , 
A)$. If $A \ra B$ is a morphism of $C^*$ algebras, then there are induced 
homomorphisms of abelian groups $K_i (A) \ra K_i (B)$. Bott periodicity 
provides a six term $K$ theory exact sequence for each exact sequence $0 
\ra J \ra A \ra B \ra 0$ of $C^*$ algebras and excision shows that the $K$ groups 
involved in 
the exact sequence only depend on the respective $C^*$ algebras. As an exercice to 
appreciate 
the power of this abstract tool one should for instance use the six term $K$ 
theory 
exact sequence
 to give a short proof of the 
Jordan curve theorem.

\noindent Discrete groups, Lie groups, group actions and foliations give rise 
through 
their convolution algebra to a canonical $C^*$ algebra, and hence to $K$ 
theory groups. The analytical meaning of these $K$ theory groups is clear 
as a receptacle for indices of elliptic operators. However, these groups 
are difficult to compute. For instance, in the case of semi-simple Lie 
groups the free abelian group with one generator for each irreducible 
discrete series representation is contained in $K_0 \, C_r^* G$ where $C_r^* G$ 
is the reduced $C^*$ algebra of $G$. Thus an explicit determination of the 
$K$ theory in this case in particular involves an enumeration of the 
discrete series.

\noindent We introduced with P. Baum \cite{[BC]} a geometrically
defined $K$  theory which 
specializes to discrete groups, Lie groups, group actions, and foliations. 
Its main features are its computability and the simplicity of its 
definition. In the case of semi-simple Lie groups it elucidates the role of 
the homogeneous space $G/K$ ($K$ the maximal compact subgroup of $G$) in 
the Atiyah-Schmid geometric construction of the discrete series \cite{[4]}. Using 
elliptic operators we constructed a natural map $\mu$ from our geometrically 
defined $K$ theory groups to the above analytic ({\it i.e.} $C^*$ algebra) $K$ 
theory groups. Much progress has been made in the past years to 
determine the range of validity of the isomorphism between the geometrically 
defined $K$ theory groups and the above analytic ({\it i.e.} $C^*$ algebra) $K$ 
theory groups.  We refer to the three Bourbaki seminars 
\cite{BBB}, \cite{BBBB}, 
\cite{BBBBB} for an update on this topic and for a precise account
of the various  contributions.
Among the most important contributions are those of Kasparov and Higson who showed 
that the 
conjectured isomorphism holds for all
amenable groups, thus proving the Novikov conjecture for all 
amenable groups and the Kadison conjecture (i.e. the absence of nontrivial idempotents
in the reduced $ C^*$-algebra) for all torsion free amenable groups. The 
conjectured isomorphism also holds for real semi-simple Lie groups thanks in particular 
to the 
work of A. Wassermann.
Moreover the recent work of V. Lafforgue
crossed the barrier of property T, showing that it holds for cocompact subgroups of 
rank one 
Lie groups and also of $SL(3, \Rb)$ or of p-adic Lie groups.
 He also gave the first general conceptual proof of the isomorphism for real or 
p-adic 
semi-simple Lie groups (and as a corollary a direct K-theoretic proof of the construction 
of all discrete series representations by Dirac-induction).
The proof of the isomorphism is certainly accessible for all connected locally 
compact 
groups.
 The proof by G. Yu of the analogue (due to J. Roe) of the conjecture in the context 
of 
coarse geometry for metric spaces which are uniformly embeddable in hilbert space, 
and the 
work of G. Skandalis J. L. Tu, J. Roe and N. Higson on the groupoid case
got very striking consequences such as the injectivity of the map $\mu$ for exact 
$C^*_r(\Gamma)$ due to  Kaminker,  Guentner and Ozawa, but recent progress
 due to Gromov, Higson, Lafforgue and Skandalis gives counterexamples to the general 
conjecture for locally compact groupoids
 for 
the simple reason that the functor $G \rightarrow K_0(C^*_r(G))$ is not half exact, 
unlike 
the functor given by the geometric group.
This makes the general problem of computing $K(C^*_r(G))$  really interesting. It 
shows that 
besides determining the large class of locally compact groups
 for which the original conjecture is valid, one should understand how to 
take homological algebra into account to deal with the correct general formulation. 

\bigskip
\section{Differential Topology}

\noindent The development of differential geometric ideas,
including de Rham homology, connections and curvature
of vector bundles, etc... took place during the eighties
thanks to cyclic cohomology which came from two different horizons  
(\cite{C2} \cite{Co$_{15}$} \cite{Co$_{17}$} \cite{Co$_{18}$} \cite{Ts$_1$}).

\noindent In the commutative case, for a compact space $X$, 
we have at our disposal in $K$-theory a tool of great relevance, the Chern character
\begin{equation}
{\rm ch} : K^* (X) \ra  H^*(X,\Qb)   \label{eq:(7.1)}
\end{equation}
which relates the $K$-theory of $X$ to the cohomology of $X$. 
When $X$ is a smooth manifold the Chern character may be 
calculated explicitly by the differential calculus of forms, 
currents, connections and curvature. More precisely, 
given a smooth vector bundle $E$ over $X$, 
or equivalently the finite projective module, $\Ec = 
C^{\infty} (X,E)$ over $\Ac = C^{\infty} (X)$ of smooth 
sections of $E$, the Chern character of $E$
\begin{equation}
{\rm ch} (E) \in H^* (X,\Rb)   \label{eq:(7.2)}
\end{equation}
is represented by the closed differential form:
\begin{equation}
{\rm ch} (E) = \hbox{trace} \, (\exp (\nabla^2 / 2\pi i)) \label{eq:(7.3)}
\end{equation}
for any connection $\nabla$ on the vector bundle $E$. Any 
closed de Rham current $C$ on the manifold $X$ determines a 
map $\vp_C$ from $K^* (X)$ to $\Cb$ by the equality
\begin{equation}
\vp_C (E) = \langle C,{\rm ch} (E) \rangle \label{eq:(7.4)}
\end{equation}
where the pairing between currents and differential forms is 
the usual one.

\noindent One obtains in this way numerical invariants of $K$-theory 
classes whose knowledge for arbitrary closed currents $C$ is 
equivalent to that of ${\rm ch} (E)$.

\noindent The noncommutative torus gave a striking example where it was 
obviously worthwhile to adapt the above construction of 
differential geometry to the noncommutative framework 
(\cite{C3}). As an easy preliminary step towards cyclic 
cohomology one can reformulate the essential ingredient of 
the construction without direct reference to derivations in 
the following way (\cite{Co$_{18}$}).

\medskip

\noindent By a cycle of dimension 
$n$ we mean a triple $(\Om , d , \int )$ where $(\Om , d )$ is a
graded differential algebra, 
and $\int : \Om^n \ra \Cb$ is a closed graded trace on $\Om$.

\smallskip

\noindent  Let $\Ac$ be an algebra over $\Cb$. Then a cycle over 
$\Ac$ is given by a cycle $(\Om , d , \int )$ and a 
homomorphism $\rho : \Ac \ra \Om^0$.

\noindent Thus a {\it cycle} over an algebra $\Ac$ is a way to embed 
$\Ac$ as a subalgebra of a differential graded algebra (DGA). We 
shall see in f) below the role of the graded trace.

\noindent The usual notions of connection and curvature extend in a 
straightforward manner to this context (\cite{Co$_{18}$}).

\noindent  Let $\Ac \build 
\longra_{}^{\rho} \Om$ be a cycle over $\Ac$, and $\Ec$ a 
finite projective module over $\Ac$. Then a connection 
$\nabla$ on $\Ec$ is a linear map $\nabla : \Ec \ra \Ec 
\otimes_{\Ac} \Om^1$ such that
\begin{equation}
\nabla (\xi x) = (\nabla \xi) x + \xi \otimes d\rho (x) \ , 
\quad \forall \, \xi \in \Ec \ , \quad x \in \Ac \, . \label{eq:(7.6)}
\end{equation}

\medskip

\noindent Here $\Ec$ is a {\it right} module over $\Ac$ and $\Om^1$ is 
considered as a bimodule over $\Ac$ using the homomorphism 
$\rho : \Ac \ra \Om^0$ and the ring structure of $\Om^*$. Let 
us list a number of easy properties (\cite{Co$_{18}$}):

\medskip

\noindent  a)  Let $e \in {\rm 
End}_{\Ac} (\Ec)$ be an idempotent and $\nabla$ a connection 
on $\Ec$; then $\xi \mpo (e \otimes 1) \nabla  \xi$ is a 
connection on $e \Ec$.

\smallskip

\noindent b) Any finite projective module $\Ec$ admits a 
connection.

\smallskip

\noindent c)The space of connections is an affine space over the 
vector space
\begin{equation}
{\rm Hom}_{\Ac} (\Ec , \Ec \otimes_{\Ac} \Om^1) \, .\label{eq:(7.7)}
\end{equation}

\noindent d) Any connection $\nabla$ extends uniquely to a linear 
map of $\widetilde{\Ec} = \Ec \otimes_{\Ac} \Om$ into itself 
such that
\begin{equation}
\nabla (\xi \otimes \om) = (\nabla \xi) \om + \xi \otimes d 
\om \ , \quad \forall \, \xi \in \Ec \ , \quad \om \in \Om \, 
.\label{eq:(7.8)}
\end{equation}

\noindent  e) The map $\t = \nabla^2$ of $\widetilde{\Ec}$ to 
$\widetilde{\Ec}$ is an endomorphism: $\t \in {\rm End}_{\Om} 
(\widetilde{\Ec})$ and with $\d(T) = \nabla T - (-1)^{degT} T \nabla$, one has 
$\d^2 (T) = \t T - T \t$ for all $T 
\in {\rm End}_{\Om} (\widetilde{\Ec})$.

\smallskip

\noindent  f) For $n$ even, $n = 2m$, the equality
 \begin{equation}
\langle [\Ec] , [\tau] \rangle = 
\frac{1}{m!} \, \int \, \t^m,\label{eq:(7.9)}
\end{equation}
 defines an additive 
map from the $K$-group $K_0 (\Ac)$ to the scalars.

\medskip
\noindent Of course one can reformulate f) by dualizing the 
closed graded trace $\int$, i.e. by considering the homology 
of the quotient $\Om / [\Om , \Om]$ (\cite{[Kar]}) and one might be tempted 
at first sight to assert that a noncommutative algebra often 
comes naturally equipped with a natural embedding in a DGA
which should suffice for the Chern character. 
This however would be rather naive and would overlook for 
instance the role of {\it integral} cycles for which the 
above additive map only affects {\it integer} values.

\noindent  The 
starting point of cyclic cohomology is the ability to compare 
different cycles on the same algebra. In fact the invariant 
of $K$-theory defined in f) by a given cycle only depends on 
the multilinear form
 \begin{equation}
\vp (a^0 , \ldots , a^n) = \int \rho (a^0) \, d (\rho (a^1)) 
\, d (\rho (a^2)) \ldots d (\rho (a^n)) \qquad \forall \, a^j 
\in \Ac \label{eq:(7.10)}
\end{equation}
(called the character of the cycle) and the functionals thus 
obtained are exactly those multilinear forms on $\Ac$ such 
that

\smallskip

$\vp$ is {\it cyclic} i.e.
 \begin{equation}
\vp (a^0 , a^1 , \ldots , a^n) = (-1)^n \, \vp (a^1 , a^2 , 
\ldots , a^0) \qquad \forall \, a_j \in \Ac \, ,\label{eq:(7.11)}
\end{equation}

$b \vp = 0$ where
 \begin{equation}
(b \vp) (a^0 , \ldots , a^{n+1}) =  
\sum_0^n (-1)^j \, \vp (a^0 , \ldots , a^j a^{j+1} , \ldots , 
a^{n+1}) + (-1)^{n+1} \, \vp (a^{n+1} a^0 , a^1 , \ldots , 
a^n) \, .
 \label{eq:(7.12)}
\end{equation}
This second condition means that $\vp$ is a Hochschild 
cocycle. In particular such a $\vp$ admits a Hochschild class
 \begin{equation}
I(\vp) \in H^n (\Ac , \Ac^*) \label{eq:(7.13)}
\end{equation}
for the Hochschild cohomology of $\Ac$ with coefficients in 
the bimodule $\Ac^*$ of linear forms on $\Ac$.

\noindent The $n$-dimensional {\it cyclic cohomology} of $\Ac$ is 
simply the cohomology $HC^n (\Ac)$ of the {\it subcomplex} of 
the Hochschild complex given by cochains which are {\it 
cyclic} i.e. fulfill \ref{eq:(7.11)}. One has an obvious ``forgetful'' map
\begin{equation}
HC^n (\Ac) \build \longra_{}^{I} H^n (\Ac , \Ac^*) \label{eq:(7.14)}
\end{equation}
but the real story starts with the following long exact 
sequence which allows in many cases to compute cyclic 
cohomology from the $B$ operator acting on Hochschild 
cohomology:

\medskip

\noindent {\bf Theorem 1.} {\it The following triangle is 
exact:}
$$
\matrix{
&H^* (\Ac , \Ac^*) \cr
^B \swarrow &&\nwarrow ^I \cr\cr
HC^* (\Ac) &\build \longra_{}^{S} &HC^* (\Ac) \cr
}
$$
\medskip 

\noindent The operator $S$ is obtained by tensoring cycles by the canonical
2-dimensional generator of the cyclic cohomology of $\Cb$.

\noindent The operator $B$ is explicitly defined at the cochain 
level by the equality
$$
\displaylines{
B = AB_0 \, , \ B_0 \, \vp (a^0 , \ldots ,a^{n-1}) = \vp 
(1,a^0 , \ldots , a^{n-1}) - (-1)^n \, \vp (a^0 , \ldots 
,a^{n-1} , 1) \cr
(A\psi) (a^0 , \ldots ,a^{n-1}) = \sum_{0}^{n-1} 
(-1)^{(n-1)j} \psi (a^j , a^{j+1} , \ldots , a^{j-1}) \, . 
\cr
}
$$
Its conceptual origin lies in the notion of cobordism of 
cycles which allows to compare different inclusion of $\Ac$ 
in DGA as follows. By a {\it chain} of dimension $n+1$ we 
shall mean a quadruple $(\Om , \partial \Om , d , \int )$ where 
$\Om$ and $\partial \Om$ are differential graded algebras of 
dimensions $n+1$ and $n$ with a given surjective morphism $r 
: \Om \ra \partial \Om$ of degree $0$, and where $\int : 
\Om^{n+1} \ra \Cb$ is a graded trace such that
\begin{equation}
\int d \om = 0 \ , \quad \forall \, \om \in \Om^n \ \hbox{such 
that} \ r (\om) = 0 \, .\label{eq:(7.16)}
\end{equation}
By the {\it boundary} of such a chain we mean the cycle 
$(\partial \Om ,d, \int')$ where for $\om' \in (\partial 
\Om)^n$ one takes $\int' \om' = \int d \om$ for any $\om \in 
\Om^n$ with $r(\om) = \om'$. One easily checks, using the 
surjectivity of $r$, that $\int'$ is a graded trace on 
$\partial \Om$ and is closed by construction.

\noindent We shall say that two cycles $\Ac \build \longra_{}^{\rho} \Om$ 
and $\Ac \build 
\longra_{}^{\rho'} \Om'$  over $\Ac$ are {\it cobordant}  if there 
exists a chain $\Om''$ with boundary $\Om \oplus 
\widetilde{\Om'}$ (where $\widetilde{\Om'}$ is obtained from 
$\Om'$ by changing the sign of $\int$) and a homomorphism 
$\rho'' : \Ac \ra \Om''$ such that $r \circ \rho'' = (\rho , 
\rho')$.

\noindent The conceptual role of the operator $B$ is clarified by the 
following result,

\medskip
\medskip

\noindent {\bf Theorem 2.}  {\it Two cycles over $\Ac$ are 
cobordant if and only if their characters $\tau_1 , \tau_2 
\in HC^n (\Ac)$ differ by an element of the image of $B$, 
where}
$$
B : H^{n+1} (\Ac , \Ac^*) \ra H C^n (\Ac) \, .
$$

\smallskip

\noindent  The operators $b,B$ given as above by
$$
\displaylines{
(b\vp) (a^0 , \ldots , a^{n+1}) = \cr
\sum_0^n (-1)^j \, \vp (a^0 , \ldots , a^j a^{j+1} , \ldots , 
a^{n+1}) + (-1)^{n+1} \, \vp (a^{n+1} a^0 , a^1 , \ldots , a^n)
}
$$
$$
\displaylines{
B = AB_0 \, , \ B_0 \, \vp (a^0 , \ldots , a^{n-1}) = \vp (1, 
a^0 , \ldots , a^{n-1}) - (-1)^n \, \vp (a^0 , \ldots , 
a^{n-1}, 1) \cr
(A\psi) (a^0 , \ldots , a^{n-1}) = \sum_{0}^{n-1} (-1)^{(n-1)j} 
\, \psi (a^j , a^{j+1} , \ldots , a^{j-1})
}
$$
satisfy $b^2 = B^2 = 0$ and $bB = -Bb$ and periodic cyclic 
cohomology which is the inductive limit of the $H C^n (\Ac)$ under
 the periodicity map $S$ admits an equivalent description as the cohomology 
of the $(b,B)$ bicomplex.

\noindent With these notations one has the 
following formula for the Chern character of the class of an idempotent $e$, up to 
normalization one has
\begin{equation}
 Ch_n (e) = \, (e-1/2) \ot e \ot e \ot ...  \ot e,
\end{equation}
where $\ot$ appears 2n times in the right hand side of the equation.

\noindent  Both the Hochschild and Cyclic 
cohomologies of the algebra $\Ac = C^{\infty} (V)$ of smooth 
functions on a manifold $V$ were computed in \cite{Co$_{17}$} and \cite{Co$_{18}$}.

\smallskip

\noindent Let $V$ be a smooth compact manifold and $\Ac$ the locally 
convex topological algebra $C^{\infty} (V)$. Then the following map $\vp 
\ra C_{\vp}$ is a canonical isomorphism of the continuous 
Hochschild cohomology group $H^k (\Ac , \Ac^*)$ with the space 
of $k$-dimensional de Rham currents on $V$:
$$
\langle C_{\vp} , f^0 \, d \, f^1 \wdg \ldots \wdg d \, f^k 
\rangle = \frac{1}{k!} \sum_{\s \in S_k} \ve (\s) \, \vp (f^0 , 
f^{\s (1)} , \ldots , f^{\s (k)})
$$
$\fl \, f^0 , \ldots , f^k \in C^{\infty} (V)$.

\smallskip

\noindent Under the isomorphism $C$ the operator $I \circ B : H^k (\Ac , 
\Ac^*) \ra H^{k-1} (\Ac , \Ac^*)$ is ($k$ times) the de Rham 
boundary $b$ for currents.

\medskip

\noindent {\bf Theorem 3.} {\it Let $\Ac$ be the locally convex 
topological algebra $C^{\infty} (V)$. Then

\smallskip

{\rm 1)} For each $k$, $HC^k (\Ac)$ is canonically isomorphic 
to the direct sum
$$
{\rm Ker} \, b \op H_{k-2} (V,\Cb) \op H_{k-4} (V,\Cb) \op 
\cdots
$$
where $H_q (V,\Cb)$ is the usual de Rham homology of $V$ and 
$b$ the de Rham boundary.

\smallskip

{\rm 2)} The periodic cyclic cohomology of $C^{\infty} (V)$ is 
canonically isomorphic to the de Rham homology $H_* (V,\Cb)$, 
with filtration by dimension.}

\medskip

\noindent As soon as we pass to the noncommutative case, more subtle 
phenomena arise. Thus for instance the filtration of the 
periodic cyclic homology (dual to periodic cyclic cohomology) 
together with the lattice $K_0 (\Ac) \sbs HC_{\rm ev} (\Ac)$, 
for $\Ac = C^{\infty} (\Tb_{\t}^2)$, gives an even analogue of 
the Jacobian of an elliptic curve. More precisely the 
filtration of $HC_{\rm ev}$ yields a canonical foliation of the 
torus $HC_{\rm ev} / K_0$ and one can show that the foliation 
algebra associated as above to the canonical transversal 
segment $[0,1]$ is isomorphic to $C^{\infty} (\Tb_{\t}^2)$. 

\smallskip

\noindent A simple example of cyclic cocycle on a nonabelian group ring 
is provided by the following formula. Any {\it group cocycle} 
$c \in H^* (B\G) = H^* (\G)$ gives rise to a cyclic cocycle 
$\vp_c$ on the algebra $\Ac = \Cb \G$
$$
\vp_c (g_0 , g_1 , \ldots , g_n) = \left\{ \matrix{
0 \hfill &\hbox{if} &g_0 \ldots g_n \not= 1 \hfill \cr
c(g_1 , \ldots , g_n) \hfill &\hbox{if} &g_0 \ldots g_n = 1 
\hfill \cr
} \right.
$$
where $c \in Z^n (\G , \Cb)$ is suitably normalized, and the 
formula is extended by linearity to $\Cb \G$. The cyclic 
cohomology of group rings is given by,

\medskip

\noindent {\bf Theorem 4.} \cite{Bu} {\it Let $\G$ be a 
discrete group, $\Ac = \Cb \G$ its group ring.

\smallskip

{\rm a)} The Hochschild cohomology $H^* (\Ac , \Ac^*)$ is 
canonically isomorphic to the cohomology $H^* ((B\G)^{\Sb^1} , 
\Cb)$ of the free loop space of the classifying space of $\G$.

\smallskip

{\rm b)} The cyclic cohomology $HC^* (\Ac)$ is canonically 
isomorphic to the $\Sb^1$-equivariant cohomology $H_{\Sb^1}^* 
((B\G)^{\Sb^1} , \Cb)$.}

\medskip

\noindent The role of the free loop space in this theorem is not 
accidental and is clarified in general by the equality
$$
B \Lb = BS^1
$$
of the classifying space $B\Lb$ of the {\it cyclic category} 
with the classifying space of the compact group $S^1$.
We refer to appendix XVIII for this point.

\smallskip

\noindent As we saw in section V the integral curvature of vector bundles 
on $\Tb_{\t}^2$ was surprisingly giving an integer, in spite of 
the irrationality of $\t$. The conceptual understanding of this 
type of integrality result lies in the existence of a natural 
lattice of {\it integral cycles} which we now describe.

\medskip

\noindent {\bf Definition.} {\it Let $\Ac$ be an algebra, a 
Fredholm module over $\Ac$ is given by:

\smallskip

{\rm 1)} a representation of $\Ac$ in a Hilbert space $\Hc$;

\smallskip

{\rm 2)} an operator $F = F^*$, $F^2 = 1$, on $\Hc$ such that
$$
[F , a] \ \hbox{is a compact operator for any} \ a \in \Ac \, .
$$
}

\smallskip

\noindent Such a Fredholm module will be called {\it odd}. An {\it even} 
Fredholm module is given by an odd Fredholm module $(\Hc , F)$ 
as above together with a $\Zb / 2$ grading $\g$, $\g = \g^*$, 
$\g^2 = 1$ of the Hilbert space $\Hc$ such that:

\smallskip

a) $\g a = a \g$ $\fl \, a \in \Ac$

\smallskip

b) $\g F = -F \g$.

\smallskip

\noindent The above definition is, up to trivial changes, the same as 
Atiyah's definition \cite{[AT]} of abstract elliptic operators, and the 
same as Kasparov's definition \cite{[18]} for the cycles in $K$-homology, 
$KK(A,\Cb)$, when $A$ is a $C^*$-algebra.

\smallskip

\noindent The main point is that a Fredholm module over an algebra $\Ac$ 
gives rise in a very simple manner to a DGA containing $\Ac$. 
One simply defines $\Om^k$ as the linear span of operators of 
the form,
$$
\om = a^0 \, [F , a^1] \ldots [F , a^k] \qquad a^j \in \Ac
$$
and the differential is given by
$$
d\om = F \om - (-1)^k \, \om F \qquad \fl \, \om \in \Om^k \, .
$$
One easily checks that the ordinary product of operators gives 
an algebra structure, $\Om^k \, \Om^{\ell} \sbs \Om^{k+\ell} $ 
and that $d^2 = 0$ owing to $F^2 = 1$.

\smallskip

\noindent Moreover if one assumes that the size of the differential $da = 
[F,a]$ is controlled, i.e. that
$$
\vert da \vert^{n+1} \quad \hbox{is trace class},
$$
 then one obtains a natural closed graded trace of degree $n$ 
by the formula,
$$
\int \om = {\rm Trace} \, (\om)
$$
(with the supertrace ${\rm Trace} \, (\g \om)$ in the even 
case, see \cite{[Co]} for details).

\noindent Hence the original 
Fredholm module gives rise to a {\it cycle} over $\Ac$. Such 
cycles have the remarkable {\it integrality} property that 
when we pair them with the $K$ theory of $\Ac$ we only get {\it 
integers} as follows from an elementary index formula 
(\cite{[Co]}).

\noindent We let $Ch_*(\Hc , F) \in HC^n(\Ac)$ be the character of the cycle associated
 to a Fredholm module $(\Hc , F)$ over $\Ac$.
This formula defines the Chern character in $K$-homology.

\noindent Cyclic cohomology got many applications \cite{L}, it led for instance
to the proof of the Novikov conjecture for hyperbolic
groups \cite{[C-M1]}. Basically,
by extending the Chern-Weil characteristic classes
to the general framework it allows for many concrete
computations of differential geometric nature on noncommutative spaces.
It also showed the depth of the relation between the classification of 
factors 
and the 
geometry of foliations.

\noindent Von Neumann algebras arise very naturally in geometry 
from foliated manifolds $(V,F)$. The 
von Neumann algebra $L^{\infty} (V,F)$ of a foliated manifold is easy to describe, 
its 
elements are random operators $T = 
(T_f)$, i.e. bounded measurable families of operators $T_f$ parametrized  
by the leaves $f$ of the foliation. For each leaf $f$ the operator $T_f$ acts in 
the Hilbert space $L^2 (f)$ of square integrable densities on the 
manifold $f$. Two random operators are identified if they are equal for 
almost all leaves $f$ (i.e. a set of leaves whose union in $V$ is 
negligible). The algebraic operations of sum and product are given by,
\begin{equation}
(T_1 + T_2)_f = (T_1)_f + (T_2)_f \, , \ (T_1 \, T_2)_f = (T_1)_f \, (T_2)_f 
\,  ,  \label{eq:(38)}
\end{equation} 
i.e. are effected pointwise.

\noindent All types of factors occur from this geometric construction and the 
continuous 
dimensions of Murray and von-Neumann play an essential role in the longitudinal 
index 
theorem. 

\noindent  Using cyclic cohomology together with the 
following 
simple fact,
\begin{equation}
\matrix{
&\hbox{``A connected group can only act trivially on a homotopy} \cr 
&\hbox{invariant cohomology theory'',} \hfill \cr
} \label{eq:(42)}
\end{equation} 
one proves (cf. \cite{[Co1]}) that for any codimension one foliation $F$ of  
a compact manifold $V$ with non vanishing Godbillon-Vey class one has,
\begin{equation}
{\rm Mod} (M) \ \hbox{has finite covolume in} \ \Rb_+^* \, ,  \label{eq:(43)}
\end{equation} 
where ${\rm Mod} (M)$ is the flow of weights of $M=L^{\infty} (V,F)$.

\noindent In the recent years J. Cuntz and D. Quillen (\cite{CQ} \cite{C-Q$_2$} 
\cite{C-Q$_4$} ) have developed a powerful new approach 
to cyclic cohomology which allowed them to prove excision in full generality.
\bigskip

\bigskip
\section{Calculus and Infinitesimals}

 \noindent The central notion of noncommutative geometry comes from
the identification of the noncommutative analogue of the two basic concepts in 
Riemann's
formulation of Geometry, namely those of manifold and
of infinitesimal line element. Both of these noncommutative analogues are of 
spectral nature and combine 
to give rise to the notion of spectral triple and spectral manifold, which will be 
described below.
\noindent  We shall first describe an operator theoretic framework for the 
calculus 
of infinitesimals which will 
provide a natural home for the line element $ds$.

\noindent  I first have to make a little excursion, and I want it as
naive as possible.  I want  to turn back to an extremely naive question 
about
what is an infinitesimal.  Let me first explain one answer that was proposed for
this intuitive idea  of infinitesimal and let me explain why this answer is
not satisfactory and then give another answer  which hopefully is satisfactory.
So,
I remember quite a long time ago to have seen an answer which was proposed
by non standard analysis. The book I was reading \cite{[B-W]} was
starting from the following problem:

 \noindent You play
a
game of throwing darts at some target called
$\Omega$ 

$$
\hbox{
\psfig{figure=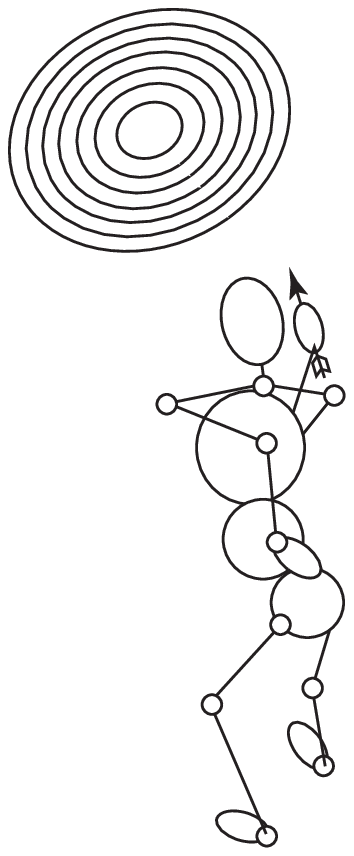}
}
$$

\noindent and the question which is asked is: what is the probability
$dp
(x)$ that actually when you send the dart you land exactly at a given point $x \in 
\Omega$?
Then the following argument was given: certainly this probability
$dp(x)$ is smaller than $1/2$ because you can cut the target into two equal 
halves, only one of which contains $x$. For the same reason $dp(x)$
is smaller than $1/4$, and so on and so forth. So
what
you find out is  that $dp(x)$ is smaller than any positive real number $\epsilon$. 
 On the other hand, if you give the answer that $dp(x)$ is $0$,
this
is not really satisfactory, because whenever you send the dart it will land
somewhere. So now, if you ask a mathematician about this naive question, he might 
very well
answer: well, $dp(x)$ is a 2-form, or it's a measure, or something like
that. But then you can try to ask him more precise questions, for
instance "what is the exponential of  $-\frac 1{ dp(x)}$ ".
And
then it will be hard for him to give a satisfactory answer, because you know that the
Taylor expansion of the function $f(y)=e^{-\frac 1{ y}}$ is zero at $y=0$. Now the book I was reading claimed to
give
an answer, and it was what is called  a non standard number. So I
worked on this theory for some time, learning some logics, until eventually
I realized there was a very bad obstruction preventing one to get concrete answers. 
It is
the following: it's a little lemma that one can easily prove, that if you
are
given a non standard number you can canonically produce  a subset
of
the interval which is not Lebesgue measurable.  Now we know from logic (from results 
of
Paul Cohen and Solovay) that it will forever be impossible to
produce
explicitely a subset of the real numbers, of the interval $[0,1]$, say, that is not
Lebesgue measurable. So, what this says is that for instance in this
example,
nobody will actually be able  to name a non standard number. A nonstandard number is 
some 
sort of chimera
which is impossible to grasp and certainly not a concrete object. In
fact when you look at nonstandard analysis you find out that except for the use of 
ultraproducts, which is very efficient, it just shifts the order in logic by one
step; it's not doing much more. Now, what  I want to explain is  that to
the above
 naive question there is a very beautiful and simple answer which
is provided by quantum mechanics. This answer will be obtained just by going
through the usual dictionary of quantum mechanics, but looking at it more
closely. So, let us thus look at the first two lines of the following dictionary 
which 
translates
classical notions into the language of operators in the Hilbert space $\Hc$:
\[
\begin{array}{cc}
\hbox{Complex variable} &\hbox{Operator in} \
\Hc \\
 & \\
\hbox{Real variable} &\hbox{Selfadjoint operator} \\
 & \\
\hbox{Infinitesimal} &\hbox{Compact operator} \\
 & \\
\hbox{Infinitesimal of order} \ \a &\hbox{Compact operator with characteristic 
values} \\
&\mu_n \ \hbox{satisfying } \ \mu_n = O(n^{-\a}) \ , \
n\ra \ify \\
\hbox{Integral of an infinitesimal }
&{\int \!\!\!\!\! {{\scriptstyle -}}} \, T = \
\hbox{ Coefficient of logarithmic} \\
 \hbox{of order 1}
&\hbox{divergence  in the trace of } \ T \, . \\
\end{array}
\]

\smallskip

\noindent The first two lines of the dictionary are familiar from quantum
mechanics. The range  of a complex variable corresponds to the
{\it spectrum} of an operator. The holomorphic functional calculus gives a
meaning to $f(T)$ for all holomorphic functions $f$ on the spectrum of
$T$. It is only holomorphic functions which operate in this generality
which reflects the difference between  complex and real analysis.
When $T=T^*$ is selfadjoint then  $f(T)$ has a meaning for all
Borel functions $f$.

\noindent The size of the  infinitesimal $T \in \Kc$
is governed by the order of decay of the sequence of characteristic values 
$\mu_n = \mu_n (T)$  as $n \ra \ify$. In particular, for all
real positive $\a$ the following condition defines infinitesimals of order
$\a$:
\begin{equation}
\mu_n (T) = O (n^{-\a}) \qquad \hbox{when} \ n\ra
\ify                                    \label{eq:(2.10)}
\end{equation}
(i.e. there exists   $C>0$ such that  $\mu_n (T) \leq C
n^{-\a} \quad \fl \, n\geq 1$). Infinitesimals of order $\a$
also  form a two--sided ideal and moreover,
\begin{equation}
T_j \ \hbox{of order} \ \a_j \Ra T_1 T_2 \
\hbox{of order} \ \a_1 + \a_2 \, .              \label{eq:(2.11)}
\end{equation}

\smallskip

\noindent  Hence, apart from commutativity, intuitive
properties of the infinitesimal calculus are fulfilled.

\noindent  Since the size of an infinitesimal is measured by the
sequence $ \ \mu_n
\downarrow 0$  it might seem that one does not need the operator
formalism at all, and that it would be enough to replace the ideal $\Kc$ in
$\Lc (\Hc)$ by the ideal $c_0 (\Nb)$ of sequences converging to zero  in
the algebra $\ell^{\ify} (\Nb)$ of bounded sequences. A variable would just 
be a bounded sequence, and an infinitesimal a sequence $\mu_n, \mu_n \ra 0$.
  However, this commutative version does not allow for the
existence of variables with range a continuum since all elements of
$\ell^{\ify} (\Nb)$ have a point spectrum and a discrete spectral measure.
Only {\it noncommutativity} of $\Lc (\Hc)$ allows for the coexistence of variables 
with Lebesgue spectrum together with
infinitesimal variables. As we shall see shortly, it is precisely this lack of 
commutativity
between the line element and the coordinates on a space that will provide the 
measurement of 
distances.

\smallskip

\noindent The integral is obtained by the following analysis, mainly due to 
Dixmier 
(\cite{[Dx]}), of the logarithmic
divergence of the partial traces
\begin{equation}
\Trace_N (T) = \sum_{0}^{N-1} \mu_n (T) \ , \ T\geq 0 \,
.                                                       \label{eq:(2.22)}
\end{equation}
In fact, it is useful to define $\Trace_{\Lb} (T)$ for any positive real
$\Lb > 0$ by piecewise affine interpolation for noninteger $\Lb$. 
\smallskip

\noindent  Define for all order 1 operators   $T \geq 0$
\begin{equation}
\tau_{\Lb} (T) = {1\over \log \Lb} \, \int_e^{\Lb} \,
{\Trace_{\mu} (T) \over \log \mu} \ {d\mu \over \mu}
                                                        \label{eq:(2.26)}
\end{equation}
which is the Cesaro mean  of the function ${\Trace_{\mu} (T) \over \log
\mu}$  over  the scaling group $\Rb_+^*$.

\noindent For  $T \geq 0$, an infinitesimal of order 1, one has
\begin{equation}
\Trace_{\Lb} (T) \leq C \, \log \Lb                     \label{eq:(2.27)}
\end{equation}
so that $\tau_{\Lb} (T)$ is bounded.
The essential  property is the following  {\it
asymptotic additivity} of the coefficient $\tau_{\Lb} (T)$ of
the logarithmic divergence (\ref{eq:(2.27)}):
\begin{equation}
\vert \tau_{\Lb} (T_1 +T_2) - \tau_{\Lb} (T_1) -
\tau_{\Lb} (T_2) \vert \leq 3C \ {\log (\log \Lb) \over
\log \Lb}                                               \label{eq:(2.28)}
\end{equation}
for  $T_j \geq 0$.

\smallskip

An easy consequence of (\ref{eq:(2.28)}) is that any  limit point
$\tau$ of the nonlinear functionals $\tau_{\Lb}$ for ${\Lb}\ra \infty$
defines a positive and linear trace on the two--sided ideal of infinitesimals of 
order $1$,

\noindent In practice the choice of the limit point $\tau$ is irrelevant
because in all important examples $T$ is a {\it measurable }
operator, i.e.:
\begin{equation}
\tau_{\Lb} (T) \ \hbox{converges when } \ \Lb \ra
\ify \, .                                               \label{eq:(2.30)}
\end{equation}
Thus  the value $\tau (T)$ is independent of the choice of the limit point $\tau$ 
and is 
denoted
\begin{equation}
{\int \!\!\!\!\!\! -} \ T \, .                          \label{eq:(2.31)}
\end{equation}
 The first interesting example is provided by pseudodifferential
operators $T$ on a differentiable manifold $M$. When $T$ is of order 1
in the above sense, it is measurable and ${\int \!\!\!\!\!
-} T$ is the non-commutative residue of $T$ (\cite{[Wo]}). 
It has a local expression in terms of the distribution
kernel $k(x,y)$, $x,y \in M$. For $T$ of order $1$ the kernel $k(x,y)$
diverges logarithmically near the diagonal,
 \begin{equation}
k(x,y) = - a(x) \log \vert x-y \vert + 0(1) \ (\hbox{for} \ \  y \ra x)
                                                        \label{eq:(2.32)}
\end{equation}
 where   $ a(x)$ is a 1--density independent of the choice of Riemannian
distance $\vert x-y \vert$.
Then one has  (up to  normalization),
\begin{equation}
{\int \!\!\!\!\!\! -} \ T = \int_M a(x).                \label{eq:(2.33)}
\end{equation}
The right hand side of this formula makes sense for all
pseudodifferential operators (cf. \cite{[Wo]}) since one
can see that the kernel of such an operator is asymptotically of the form
\begin{equation}
k(x,y) = \sum a_k (x,x-y) - a(x) \log \vert x-y \vert +
0(1)                                                    \label{eq:(2.34)}
\end{equation}
where $a_k (x,\xi)$ is homogeneous of  degree $-k$ in $\xi$,
and the  1--density $a(x)$ is defined  intrinsically.

\smallskip

\noindent The same principle  of  extension of
${\int \!\!\!\!\! -}$ to infinitesimals of order
$<1$ works for hypoelliptic operators
and more generally as we shall see below, for  spectral triples
whose  dimension spectrum is simple.

 \noindent  We can now
go back to our initial naive question about the target and the darts, we find that 
quantum 
mechanics gives us  an
obvious infinitesimal which answers the question: it is the inverse of the Dirichlet 
Laplacian
for
the domain $\Omega$. Thus there is now a clear meaning for the exponential of
$\frac {-1}{dp}$, that's the well known heat kernel which is an infinitesimal
 of arbitrarily large order as we expected from the Taylor expansion. 

\noindent From the H. Weyl theorem on the asymptotic
behavior of eigenvalues of $\D$ it follows that $dp$ is of order 1, and that given a 
function 
$f$ on $\Om$ the product $f
\, dp$ is measurable, while
 \begin{equation}
{\int \!\!\!\!\!\! -} \, f \, dp = \int_{\Om} f(x_1 ,x_2)
\, dx_1 \wdg dx_2                                       \label{eq:(2.36)}
\end{equation}
gives the ordinary integral of $f$  with respect to the measure given by the area of 
the 
target. 

\bigskip
  
\section{ Spectral triples}

\noindent In this section we shall come back to the two basic notions 
introduced by Riemann in the classical framework, those of {\it 
manifold} and of {\it line element}. We shall see that both of these 
notions adapt remarkably well to the noncommutative framework and 
this will lead us to the notion of spectral manifold which 
noncommutative geometry is based on.

\smallskip

\noindent In ordinary geometry of course you can give a manifold by a cooking 
recipe, by charts and local diffeomorphisms, and one could be 
tempted to propose an analogous cooking recipe in the 
noncommutative case. This is pretty much what is achieved by the 
general construction of the algebras of foliations and it is a good 
test of any general idea that it should at least cover that large 
class of examples.

\smallskip

\noindent But at a more conceptual level, it was recognized long ago by 
geometors that the main quality of the homotopy type of an oriented 
manifold is to satisfy Poincar\'e duality not only in ordinary 
homology but also in $K$-homology. Poincar\'e duality in ordinary 
homology is not sufficient to describe homotopy type of manifolds 
\cite{Mi-S} but D.~Sullivan \cite{Sull} showed (in the simply connected 
PL case of dimension $\geq 5$ ignoring 2-torsion) that it is 
sufficient to replace ordinary homology by $KO$-homology.
Moreover the Chern character of the $KO$-homology fundamental class
contains all the rational information on the Pontrjagin classes.

\smallskip

\noindent The characteristic property of {\it differentiable manifolds} which 
is carried over to the noncommutative case is {\it Poincar\'e 
duality} in $KO$-homology \cite{Sull}.

\smallskip

\noindent Moreover, as we saw above in the discussion of Fredholm modules, 
$K$-homology admits a fairly simple definition in terms of Hilbert 
space and Fredholm representations of algebras.

\smallskip

\noindent For an ordinary manifold the choice of the fundamental cycle in 
$K$-homology is a refinement of the choice of orientation of the 
manifold and in its simplest form is a choice of Spin-structure. Of 
course the role of a spin structure is to allow for the 
construction of the corresponding Dirac operator which gives a 
corresponding Fredholm representation of the algebra of smooth 
functions.

\smallskip

\noindent What is rewarding is that this will not only guide us towards the 
notion of noncommutative manifold but also to a formula, of 
operator theoretic nature, for the line element $ds$.

\smallskip

\noindent The infinitesimal unit of length``$ds$'' should be an infinitesimal in the sense of 
section VIII and one way to get an intuitive understanding of the 
formula for $ds$ is to consider Feynman diagrams which physicist 
use currently in the computations of quantum field theory. Let us 
contemplate the diagram 

$$
\hbox{
\psfig{figure=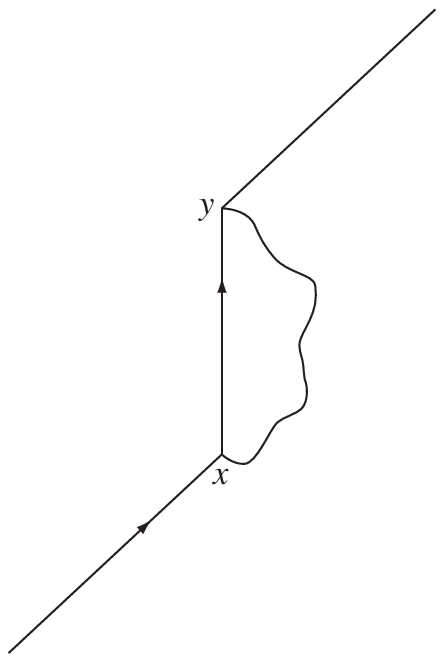}
}
$$

 \noindent which is involved in the 
computation of the self-energy of an electron in QED. The two 
points $x$ and $y$ of space-time at which the photon (the wiggly 
line) is emitted and reabsorbed are very close by and our ansatz 
for $ds$ will be at the intuitive level,
 
 \begin{equation}
ds = \times\!\!\!\!-\!\!\!-\!\!\!-\!\!\!\!\!\times \, . \label{eq:(9.1)}
\end{equation}
The right hand side has good meaning in physics, it is called the 
Fermion propagator and is given by
 \begin{equation}
\times\!\!\!\!-\!\!\!-\!\!\!-\!\!\!\!\!\times = D^{-1} \label{eq:(9.2)}
\end{equation}
where $D$ is the Dirac operator.

\smallskip

\noindent We thus arrive at the following basic ansatz,
 \begin{equation}
ds = D^{-1} \, .\label{eq:(9.3)}
\end{equation}
In some sense it is simpler than the ansatz giving $ds^2$ as 
$g_{\mu \nu} \, dx^{\mu} \, dx^{\nu}$, the point being that the 
spin structure allows really to extract the square root of $ds^2$ 
 (as is well known Dirac found the corresponding operator as 
a differential square root of a Laplacian).

\smallskip

\noindent The first thing we need to do is to check that we are still able to 
measure distances with our ``unit of length'' $ds$. In fact we saw 
in the discussion of the quantized calculus that variables with 
continuous range cant commute with ``infinitesimals'' such as $ds$ 
and it is thus not very surprising that this lack of commutativity 
allows to compute, in the classical Riemannian case, the geodesic 
distance $d(x,y)$ between two points. The precise formula is
 \begin{equation}
d(x,y) = {\rm Sup} \, \{ \vert f(x) - f(y) \vert \, ; \ f \in \Ac 
\, , \ \Vert [D,f] \Vert \leq 1 \} \label{eq:(9.4)}
\end{equation}
where $D = ds^{-1}$ as above and $\Ac$ is the algebra of smooth 
functions. Note that if $ds$ has the dimension of a length $L$, then $D$ 
has dimension $L^{-1}$ and the above expression for $d(x,y)$ also has the 
dimension of a length.

\smallskip

\noindent Thus we see in the classical geometric case that both the 
fundamental cycle in $K$-homology and the metric are encoded in the 
{\it spectral triple} $(\Ac , \Hc , D)$ where $\Ac$ is the algebra 
of functions acting in the Hilbert space $\Hc$ of spinors, while 
$D$ is the Dirac operator.

\noindent To get familiar with this notion one should check that
 we recover the volume form of the 
Riemannian
metric by the equality (valid up to a normalization constant \cite{[Co]})
\begin{equation}
{\int \!\!\!\!\!\! -}f \,\vert ds\vert^n =  \,
\int_{M_n} f \, \sqrt{g} \ d^n x                        \label{eq:(9.5)}
\end{equation}
but the first interesting point is that besides this coherence with the usual 
computations
there are new simple questions we can ask now such as "what is the two-dimensional 
measure
of a four manifold" in other words "what is its area ?". Thus one should compute
\begin{equation}
{\int \!\!\!\!\!\! -} \, ds^2                         \label{eq:(9.6)}
\end{equation}
It is obvious from invariant theory that this should be proportional to the 
Hilbert--Einstein action but doing the direct computation is a worthwile exercice
(cf. \cite{[Kas]} \cite{[K-W]}), the exact result being
\begin{equation}
{\int \!\!\!\!\!\! -} \, ds^2 = {-1 \over 48\pi^2} \,
\int_{M_4} r \, \sqrt{g} \ d^4 x                        \label{eq:(9.7)}
\end{equation}
where as above $dv=\sqrt{g} \ d^4 x$ is the volume form, 
$ds=D^{-1}$  the length element, {\it i.e.} the inverse of the Dirac operator and $r$ is 
the scalar 
curvature.

\smallskip

\noindent In the general framework of Noncommutative Geometry the confluence 
of the Hilbert space incarnation of the two notions
of  metric and fundamental class for a manifold led very naturally
to define a  geometric space as given by a {\it spectral triple:}
\begin{equation}
(\Ac ,\Hc ,D)                                           \label{eq:(9.8)}
\end{equation}
 where $\cal{A}$ is a concrete algebra of coordinates 
represented on a Hilbert  space $\cal{H}$ and the operator
 $D$ is the inverse of the line  element. 
\begin{equation}
ds=1/D.\label{eq:(9.9)}
\end{equation}
This definition is entirely spectral; the elements of the algebra 
 are operators, the points, if they exist, come from
 the joint  spectrum of operators and the line element is an operator. 

\noindent The basic properties of such spectral triples are easy to formulate 
and do not make any reference to the commutativity of the algebra 
$\Ac$. They are 
\begin{equation}
[D,a] \ \hbox{is bounded for any} \ a \in \Ac \, ,\label{eq:(9.10)}
\end{equation}
\begin{equation}
D = D^*  \ \hbox{and} \ (D+\lambda)^{-1} \ \hbox{is a compact 
operator} \ \forall \, \lambda \not\in \Cb \, .\label{eq:(9.11)}
\end{equation}
(Of course $D$ is an {\it unbounded} operator).

\smallskip

\noindent There is no difficulty to adapt the above formula for the distance 
in the general noncommutative case, one uses the same, the points $x$ and $y$ 
being replaced by arbitrary states $\vp$ and $\psi$ on the algebra 
$\Ac$.
Recall that a state is a normalized positive  linear form on $\Ac$ such that $\vp (1) = 1$,
\begin{equation}
\vp : \bar{\Ac} \ra \Cb \ , \ \vp (a^* a) \geq 0 \ ,
\quad \fl \, a \in \bar{\Ac} \ , \ \vp (1) = 1 \, .
                                                        \label{eq:(9.12)}
\end{equation}
The distance between two states is given by, 
\begin{equation}
d(\vp ,\psi) = \Sup \, \{ \vert \vp (a) - \psi (a) \vert \
; \ a\in \Ac \ , \ \Vert [D,a]\Vert \leq 1 \} \, .
                                                        \label{eq:(9.13)}
\end{equation}

\noindent The significance of $D$ is two-fold. On the one hand
it defines the metric by the above equation, on the other hand its
homotopy class represents the K-homology fundamental class of the
space under consideration.

\noindent It is crucial to understand from the start the tension between the conditions 9-\ref{eq:(9.10)} 
and  9-\ref{eq:(9.11)}. The first condition would be trivially fulfilled if $D$ were 
bounded but condition  9-\ref{eq:(9.11)} shows that it is unbounded. To understand 
this tension let us work out a very simple case. We let the algebra 
$\Ac$ be generated by a single unitary operator $U$. Let us show 
that if the index pairing between $U$ and $D$, i.e. the index of 
$PUP$ where $P$ is the orthogonal projection on the positive 
eigenspace of $D$, {\it does not vanish} then the number $N(E)$ of 
eigenvalues of $D$ whose absolute value is less than $E$ grows at 
least like $E$ when $E \ra \infty$. This means that in the above 
circumstance $ds = D^{-1}$ is of order one or less.

\smallskip

To prove this we choose a smooth function $f \in C_c^{\infty} 
(\Rb)$ identically one near 0, even and with Support $(f) \sbs 
[-1,1]$. We then let $R(\ve) = f (\ve D)$. One first shows 
(\cite{[Co]}) that the operator norm of the commutator $[R(\ve) , 
U]$ tends to 0 like $\ve$. It then follows that the trace norm 
satisfies
\begin{equation}
\Vert [R(\ve) , U] \Vert_1 \leq C \, \ve \, N(1/\ve) \label{eq:(9.14)}
\end{equation}
as one sees using the control of the rank of $R(\ve)$ from 
$N(1/\ve)$. The index pairing is given by $-\frac{1}{2} \, {\rm 
Trace} \, (U^* [F,U])$ where $F$ is the sign of $D$ and one has,
\begin{equation}
{\rm Trace} \, (U^* [F,U]) = \lim_{\ve \ra 0} {\rm Trace} \, (U^* 
[F,U] \, R(\ve)) = \lim_{\ve \ra 0} {\rm Trace} \, (U^* \, F 
[U,R(\ve)]) \, 
. \label{eq:(9.15)}
\end{equation}
Thus the limit being non zero we get a lower bound on the trace 
norm of $[U,R(\ve)]$ and hence on $\ve \, N \left( \frac{1}{\ve} 
\right)$ which shows that $N(E)$ grows at least like $E$ when $E 
\ra \infty$.

\noindent This shows that $ds$ cannot be too small (it cannot be of order $\a 
> 1$). In fact when $ds$ is of order 1 one has the following index 
formula,
\begin{equation}
{\rm Index} \, (PUP) = -\frac{1}{2} \, \int\!\!\!\!\!\!- \ U^{-1} 
[D,U] \, \vert ds \vert \, .\label{eq:(9.16)}
\end{equation}
The simplest case 
in which the index pairing between $D$ and $U$ does not vanish, 
with $ds$ of order 1, is obtained by requiring the further 
condition,
\begin{equation}
U^{-1} [D,U] = 1 \, .\label{eq:(9.17)}
\end{equation}
It is a simple exercise to compute the geometry on $S^1 = {\rm 
Spectrum} \, (U)$ given by an irreducible representation of 
condition \ref{eq:(9.17)}. One obtains the standard circle with length $2\pi$.

\noindent The above index formula is a special case of a general result (\cite{[Co]}) 
which computes the $n$-dimensional Hochschild class of the Chern 
character of a spectral triple of dimension $n$. 

\medskip

\noindent {\bf Theorem 5.} {\it Let $(\Hc , F)$ be a Fredholm 
module over an involutive algebra $\Ac$. Let $D$ be an unbounded 
selfadjoint operator in $\Hc$ such that $D^{-1}$ is of order $1/n$
, ${\rm Sign} \, D = F$, and such that for any $a 
\in 
\Ac$ the operators $a$ and $[D,a]$ are in the domain of all powers 
of the derivations $\d$, given by $\d (x) = [\vert D \vert , x]$.
Let $\tau_n \in HC^n (\Ac)$ be the Chern character of 
$(\Hc , F)$.}
\smallskip

 \noindent  {\it For every $n$-dimensional Hochschild cycle $c \in Z_n (\Ac 
, \Ac)$, $c=\sum \, a^0 \, \ot a^1 \ldots \ot a^n$, one has 
$ \langle \tau_n , c 
\rangle ={\int \!\!\!\!\!\! -} \, \sum \, a^0 \, [D,a^1] \ldots [D,a^n] \, \vert D \vert^{-n}$.} 
\medskip

\noindent  We refer to \cite{[Co]} for precise normalization and to \cite{Gracia+V} for the detailed proof.
By construction, this formula  is scale invariant, 
i.e. it remains unchanged if we replace $D$ by $\lb D$ for $\lb \in 
\Rb_+^*$. 
The operators $T_c$ of the form
\begin{equation}
T_c = \sum \, a^0 \, [D,a^1] \ldots [D,a^n] \, \vert D \vert^{-n}
\end{equation}
 are {\it measurable} in the sense of section VIII.

\smallskip

\noindent  The 
 long exact sequence of cyclic cohomology 
(Section VII) shows that the Hochschild class of $\tau_n$ is the 
obstruction to a better summability of $(\Hc , F)$, indeed $\tau_n$ belongs to 
the image $S(HC^{n-2} (\Ac))$ (which is the case if the degree of summability
can be improved by $2$) if and only if the Hochschild 
cohomology class $I(\tau_n) \in H^n (\Ac, \Ac^*)$ is equal to $0$.

\medskip

\noindent In particular, the above theorem implies nonvanishing of residues when the 
cohomological dimension of ${\rm ch}_* (\Hc , F)$ is not lower than 
$n$:

\medskip

\noindent {\bf Corollary.} {\it With the hypothesis of Theorem 5
and if the Hochschild class of ${\rm ch}_* (\Hc , F)$ pairs 
nontrivially with $H_n (\Ac, \Ac)$ one has}
\begin{equation}
{\int \!\!\!\!\!\! -} \vert D \vert^{-n} \ne 0 \, .
\end{equation}

\medskip

\noindent In other words the residue of the function $\zeta (s) = {\rm Trace} 
\, (\vert D \vert^{-s})$ at $s=n$ cannot vanish.

\smallskip

\noindent In higher dimension, the Hochschild class of the character suffices 
to determine the index pairing with the $K$-theory class of an 
idempotent $e$ provided the lower dimensional components of ${\rm 
ch} (e)$ vanish. As we saw above these components are given, up to 
normalization by,
\begin{equation}
{\rm ch}_n (e) = \left( e -  \frac{1}{2} \right) \ot e \ot \cdots 
\ot e
\end{equation}
(with $2n$ tensor signs) and as such cannot vanish. But both Hochschild and cyclic 
cohomology are Morita invariant, which implies that the class of 
${\rm ch} (e)$ in the normalized $(b,B)$ bicomplex (in homology) does not 
change when we project each of its components ${\rm ch}_n (e)$ on 
the commutant of a matrix algebra $M_q (\Cb) \sbs \Ac$. The formula 
for this projection $\langle {\rm ch}_n (e) \rangle$ in terms of 
the matrix components $e_{ij}$,
\begin{equation}
e = [e_{ij}] \ , \qquad e_{ij} \in M_q (\Cb)' \cap \Ac
\end{equation}
is the following,
\begin{equation}
\langle {\rm ch}_n (e) \rangle = \sum \left( e_{i_0 i_1} - 
\frac{1}{2} \, \d_{i_0 i_1} \right) \ot e_{i_1 i_2} \ot e_{i_2 i_3} 
\ot \cdots \ot e_{i_{2n} i_0}
\end{equation}
and there are very interesting situations in which all the lower 
components $\langle {\rm ch}_{j} (e) \rangle$ actually vanish,
\begin{equation}
\langle {\rm ch}_{j} (e) \rangle = 0 \qquad j < m \, .
\end{equation}
For $m=1$ for instance we can take $q=2$ and the condition $\langle 
{\rm ch}_0 (e) \rangle = 0$ means that $e$ is of the form,
\begin{equation}
e = \left[ \matrix{
t &z \cr
z^* &(1-t) \cr
} \right] \, .
\end{equation}
(The equation $e^2 = e$ then means that $t^2 + z^* \, z = t$, $tz + 
z(1-t) = z$, $z^* \, t + (1-t) \, z^* = z^*$, $z^* \, z + (1-t)^2 = 
(1-t)$ which shows that the algebra generated by the components 
$z$, $z^*$, $t$ of $e$ is abelian).

\smallskip

\noindent It then follows automatically that $\langle {\rm ch}_1 (e) \rangle$ 
is a Hochschild cycle and hence by theorem 5, that if $ds = D^{-1}$ 
is of order $\frac{1}{2}$ the index pairing is given by,
\begin{equation}
\hbox{Index} D^+_e  =  - \, \int\!\!\!\!\!\!- \ \g \left(e - 
\frac{1}{2}\right) \, [D,e]^2 \, ds^2 \, .\label{eq:(9.20)}
\end{equation}
Exactly as above this shows that $ds$ cannot be of order $\a > 
\frac{1}{2}$ if the index pairing is non zero, and we also get the 
analogue of equation 9-\ref{eq:(9.17)} in the form,
\begin{equation}
\left\langle \left( e - \frac{1}{2} \right) \, [D,e]^2 
\right\rangle = \g
\end{equation}
where $\langle \ \rangle$ is simply the projection on the commutant 
of $M_2 (\Cb)$ in $\Lc (\Hc)$.

\smallskip

\noindent This equation together with \ref{eq:(9.20)} implies that the area $\int\!\!\!\!\!- 
\ ds^2$ is an integer since it is given by a Fredholm index. One 
can show that the algebra $\Ac$ generated by the components of $e$ 
is $C(S^2)$ the algebra of continuous functions on $S^2$ and that 
any Riemannian metric $g$ on $S^2$ with fixed volume form gives a 
solution to the above equations.

\smallskip

\noindent There is a converse to that result (\cite{Co3}) but it requires 
the further hypothesis that $D$ is of order one:
\begin{equation}
[[D,e_{ij}] , e_{k\ell}] = 0
\end{equation}
where the $e_{ij}$ are the components of the idempotent $e$, i.e. 
are the generators of the algebra.

\smallskip

\noindent This order one condition is the counterpart in our operator 
theoretic setting of the ``quadratic'' nature of Riemann's equation 
$ds^2 = g_{\mu \nu} \, dx^{\mu} \, dx^{\nu}$. It is easier to 
formulate in terms of the square root which we extracted using the 
spin structure. We shall come later to the correct formulation of 
the order one condition when the algebra $\Ac$ is noncommutative.

\smallskip

\noindent To end this section let us move on to the four dimensional case, 
i.e. $n=2$. We take $q=4$, i.e. we deal with $M_4 (\Cb)$.

\noindent We first determine the $C^*$ algebra generated by $M_4 (\Cb)$ and a 
projection $e = e^*$ such that $\left\lgl e - \frac{1}{2} \right\rgl = 0$ as 
above and whose two by two matrix expression is of the form,
\begin{equation}
[ e^{ij} ] = \left[ \matrix{q_{11} &q_{12} \cr q_{21} &q_{22} \cr} \right] 
\label{eq58}
\end{equation}
where each $q_{ij}$ is a $2 \ts 2$ matrix of the form,
\begin{equation}
q = \left[ \matrix{\a &\b \cr -\b^* &\a^* \cr} \right] \, . \label{eq59}
\end{equation}
Since $e = e^*$, both $q_{11}$ and $q_{22}$ are selfadjoint, moreover since 
$\left\lgl e - \frac{1}{2} \right\rgl = 0$, we can 
find $t = t^*$ such that,
\begin{equation}
q_{11} = \left[ \matrix{t &0 \cr 0 &t \cr} \right] \, , \ q_{22} = \left[ 
\matrix{(1-t) &0 \cr 0 &(1-t) \cr} \right] \, . \label{eq60}
\end{equation}
We let $q_{12} = \left[ \matrix{\a &\b \cr -\b^* &\a^* \cr} \right]$, we then 
get from $e = e^*$,
\begin{equation}
q_{21} = \left[ \matrix{\a^* &-\b \cr \b^* &\a \cr} \right] \, . \label{eq61}
\end{equation}
We thus see that the commutant $\Ac$ of $M_4 (\Cb)$ is generated by $t,\a,\b$ 
and we first need to find the relations imposed by the equality $e^2 = e$.

\noindent In terms of $e = \left[ \matrix{ t &q \cr q^* &1-t \cr} \right]$, 
the equation $e^2 = e$ means that $t^2 - t + q q^* = 0$, $t^2 - t + q^* q = 
0$ and $[t,q] = 0$. This shows that $t$ commutes with $\a$, $\b$, $\a^*$ and 
$\b^*$ and since $qq^* = q^* q$ is a diagonal matrix
\begin{equation}
\a \a^* = \a^* \a \, , \ \a \b = \b \a \, , \ \a^* \b = \b \a^* \, , \ \b 
\b^* = \b^* \b \label{eq62}
\end{equation}
so that the $C^*$ algebra $\Ac$ is abelian, with the only further relation, 
(besides $t = t^*$),
\begin{equation}
\a \a^* + \b \b^* + t^2 - t = 0 \, . \label{eq63}
\end{equation}
This is enough to check that,
\begin{equation}
\Ac = C (S^4) \label{eq64}
\end{equation}
where $S^4$ appears naturally as quaternionic projective space,
\begin{equation}
S^4 = P_1 (\Hb) \, . \label{eq65}
\end{equation}
The original $C^*$ algebra is thus,
\begin{equation}
B = C (S^4) \ot M_4 (\Cb) \, . \label{eq66}
\end{equation}
We shall now check that the two dimensional component $\left\lgl Ch_1(e) \right\rgl$
automatically vanishes as an element of the (normalized) (b,B)-bicomplex so that,
\begin{equation}
\left\lgl Ch_{n}(e) \right\rgl = 0 \, , \ 
n = 0,1  . \label{eq67}
\end{equation}
With $q = \left[ \matrix{ \a &\b \cr -\b^* &\a^* \cr} \right]$, we get,
\begin{eqnarray}
&\lgl Ch_1(e) \rgl &= \Biggl\lgl \left( t - \frac{1}{2} \right) 
\, (dq \, dq^* - dq^* \, dq) \\
&&+ \, q \, (dq^* \, dt - dt \, dq^*) + q^* \, (dt \, dq - dq \, dt ) 
\Biggl\rgl \nonumber
\label{eq69}
\end{eqnarray}
where the expectation in the right hand side is relative to $M_2 (\Cb)$ and we use 
the notation $dx$ instead of the tensor notation.

\noindent The diagonal elements of $\om = dq \, dq^*$ are 
$$
\om_{11} = d \a \, d\a^* + d\b \, d\b^* \, , \ \om_{22} = d\b^* \, d\b + 
d\a^* \, d\a
$$
while for $\om' = dq^* \, dq$ we get, 
$$
\om'_{11} = d\a^* \, d\a + d\b \, d\b^* \, , \ \om'_{22} = d\b^* \, d\b + d\a 
\, d\a^* \, .
$$
It follows that, since $t$ is diagonal,
\begin{equation}
\left\lgl \left( t - \frac{1}{2} \right) \, (dq \, dq^* - dq^* \, dq) 
\right\rgl = 0 \, . \label{eq70}
\end{equation}
The diagonal elements of $q \, dq^* \, dt = \rho$ are
$$
\rho_{11} = \a \, d\a^* \, dt + \b \, d\b^* \, dt \, , \ \rho_{22} = \b^* \, 
d\b \, dt + \a^* \, d\a \, dt
$$
while for $\rho' = q^* \, dq \, dt$ they are
$$
\rho'_{11} = \a^* \, d\a \, dt + \b \, d\b^* \, dt \, , \ \rho'_{22} = \b^* \, 
d\b \, dt + \a \, d\a^* \, dt \, .
$$
Similarly for $\s = q \, dt \, dq^*$ and $\s' = q^* \, dt \, dq$ one gets the 
required cancellations so that,
\begin{equation}
\left\lgl Ch_1(e) \right\rgl = 0 \, , 
\label{eq71}
\end{equation}
It follows thus that $\left\lgl Ch_2(e) \right\rgl $ is a Hochschild cycle and that 
for any $ds= D^{-1}$ of order $\frac{1}{4}$ commuting with $M_4 (\Cb)$, the index pairing
of $ D$ with $e$ is
\begin{equation}
\hbox{Index} D^+_e  =\, \int\!\!\!\!\!\!- \ \g \left(e - 
\frac{1}{2}\right) \, [D,e]^4 \, ds^4 \, .\label{eq:(9.21)}
\end{equation}
Exactly as above this shows that $ds$ cannot be of order $\a > 
\frac{1}{4}$ if the index pairing is non zero, and we also get the 
analogue of equation 9-\ref{eq:(9.17)} in the form,
\begin{equation}
\left\langle \left( e - \frac{1}{2} \right) \, [D,e]^4 
\right\rangle = \g \label{eq:(9:36)}
\end{equation}
where $\langle \ \rangle$ is simply the projection on the commutant 
of $M_4 (\Cb)$ in $\Lc (\Hc)$.

\smallskip

\noindent This equation together with (41) implies the integrality of the 4-dimensional volume,
\begin{equation}
\int\!\!\!\!\!- 
\ ds^4 \in \Nb, \label{eq:(9.40)}
\end{equation}
 since it is given by a Fredholm index.

\noindent  One 
can show that the algebra $\Ac$ generated by the components of $e$ 
is $C(S^4)$ the algebra of continuous functions on $S^4$ and that 
any Riemannian metric $g$ on $S^4$  gives a 
solution to the above equations, provided its volume form is,

\begin{equation}
v = \frac{1}{1-2t} \ d\a \wdg d \, \ov{\a} \wdg d\b \wdg d \, \ov{\b} \, . 
\label{eq72}
\end{equation}
As in the two dimensional case there is a converse, assuming the order one condition on $D$.
  
\noindent The next question is how is $D$  to be chosen from within the homotopy class
which characterizes its $K$-homology class? 
There are two  answers to this question.   The first uses the naive idea of a formal metric, 
\begin{equation}
G = \sum_{\mu , \nu = 1}^{d}   dx^{\mu} g_{\mu \nu} (dx^{\nu})^* \in \Om_+^2 
({\cal 
A}) \, , \label{eq:(9.41)}
\end{equation}
and the choice of $D$ is performed by minimizing the action functional, 
\begin{equation}
A = \sum_{\mu , \nu = 1}^{d} {\int \!\!\!\!\! -} [ D,x^{\mu}] g_{\mu \nu} ([ D,x^{\nu}])^* \vert D^{-4} \vert
 \, , \label{eq:(9.41)}
\end{equation} 
among the $D$'s which fulfill equation (42) holding $G$ fixed.

\noindent    The minimum is then given by the Dirac operator associated to the unique Riemannian metric 
with volume form $v$ in the conformal
class of $g_{\mu \nu} dx^{\mu}dx^{\nu}$.

\noindent The second way to select $D$ from within its $K$-homology class is to use
an action functional with the largest possible invariance group which is the unitary group
of Hilbert space. The corresponding action is then spectral and only depends upon the eigenvalues
of $D$. The simplest such action is of the form, \cite{[C-C]}
\begin{equation}
S(D) = \hbox{Trace}(f(D)).
\end{equation}
where $f$ is an even function vanishing at $\infty$.
If we take for $f$ a step function equal to 1 in $[-\Lb,\Lb]$, the value of $S(D)$ is,
\begin{equation}
N(\Lb) = \# \ \hbox{eigenvalues of $D$ in} \ [-\Lb,\Lb] .
\end{equation}
This step function $N(\Lb)$ is the superposition of two terms,
$$
 N(\Lb)= \lgl N(\Lb) \rgl +N_{\rm osc} (\Lb).
$$
 The oscillatory part $N_{\rm osc} (\Lb)$ 
is the same as for a random matrix, governed by the statistic
dictated by the symmetries of the system and does not concern us here.
The average part $ \lgl N(\Lb) \rgl $ is computed by a semiclassical approximation 
and the leading term in the asymptotic expansion is,
\begin{equation}
\frac{\Lb^4}{2}  \int\!\!\!\!\!- 
\ ds^4 
\end{equation}
which by (43)  is independent of the choice of $D$ in its $K$-homology class.

\noindent If we restrict ourselves to solutions 
given by ordinary Riemannian metrics the next term in the asymptotic expansion is the 
Hilbert--Einstein action
functional for the Riemannian metric,
\begin{equation}
  \,{-\Lb^2 \over 96\pi^2} \,\int_{S_4} r \, \sqrt{g} \ d^4 x 
\end{equation}
Other nonzero terms in
the asymptotic expansion  are  cosmological, Weyl gravity
and  topological terms.

\bigskip
\section{Noncommutative 4-manifolds and the Instanton algebra}

\noindent In this section, based on our collaboration with G. Landi (\cite{CLa}),
 we shall show that the basic equation for an instanton in 
dimension 4, namely
\begin{equation}
e=e^2=e^* \label{eq1}
\end{equation}
and
\begin{equation}
\langle {\rm ch}_0 (e) \rangle = 0 \ , \quad \langle {\rm ch}_1 (e) \rangle = 0 
\label{eq2}
\end{equation}
(where ${\rm ch}_n$ are the components of the Chern character,
\begin{equation}
{\rm ch}_n (e) = \left( e - \frac{1}{2} \right) \ot e \ot \ldots \ot e 
\label{eq3}
\end{equation}
and $\langle \ \rangle$ is the projection onto the commutant of a $4 \times 4$ 
matrix algebra) do admit noncommutative solutions. In other words the algebra 
generated by the 16 components of the $4 \times 4$ matrix,
\begin{equation}
e=[e_{ij}] \label{eq4}
\end{equation}
will be noncommutative.

\smallskip

\noindent In fact this prompts us to introduce, a priori, the algebra $\Ac$ with 16 
generators $e_{ij}$ and whose presentation is given by the relations (1) and 
(2). The relation $\langle {\rm ch}_0 (e) \rangle = 0$ just means that
\begin{equation}
e_{11} + e_{22} + e_{33} + e_{44} = 2 \label{eq5}
\end{equation}
and the equation $e=e^*$ defines the involution in $\Ac$. The relation $e^2 = e$ 
is easy to comprehend as a quadratic relation between the generators.

\smallskip

\noindent The relation $\langle {\rm ch}_1 (e) \rangle = 0$ is more delicate to understand 
since it involves tensors and the simplest way to think about it is to represent 
the $e_{ij}$ as operators in Hilbert space $\Hc$. What we ask then is that,
\begin{equation}
\sum \left( e_{ij} - \frac{1}{2} \ \d_{ij} \right) \ot \widetilde{e}_{jk} \ot 
\widetilde{e}_{ki} = 0 \label{eq6}
\end{equation}
where the $\sim$ means that we take the class modulo the scalar multiples of 1.

\smallskip

\noindent This allows to define what is a unitary representation $\pi$ of the algebra $\Ac$ 
and we can endow its elements, i.e. polynomials in the noncommuting generators 
$e_{ij}$, with the $C^*$-norm,
\begin{equation}
\Vert x \Vert = \sup_{\pi}^{} \, \Vert \pi (x) \Vert \label{eq7}
\end{equation}
where $\pi$ ranges through all unitary representations. It is easy to show that 
for $x \in \Ac$ the supremum is finite since in any unitary representation, the 
$e_{ij}$ satisfy,
\begin{equation}
\Vert \pi (e_{ij}) \Vert \leq 1 \label{eq8}
\end{equation}
as matrix elements of a selfadjoint idempotent.

\medskip

\noindent {\bf Definition.} {\it We let $C ({\rm Gr})$ be the $C^*$ completion 
of $\Ac$ and $C^{\infty} ({\rm Gr})$ the smooth closure of $\Ac$ in $C({\rm 
Gr})$.}

\medskip

\noindent The letter ${\rm Gr}$ stands for the Grassmanian but our construction has little 
to do with the known ``noncommutative Grassmanians''. The really non-trivial condition is 
the cubic condition \ref{eq6}.
In fact as we saw above 
the same construction in dimension 2 does give a {\it commutative} answer namely 
$P_1 (\Cb)$.

\smallskip

\noindent One should observe from the outset that the compact Lie group $SU(4)$ acts by 
automorphisms,
\begin{equation}
PSU(4) \subset {\rm Aut} \, (C^{\infty} ({\rm Gr})) \label{eq9}
\end{equation}
by the following operation,
\begin{equation}
e \ra U \, e \, U^* \label{eq10}
\end{equation}
where $U \in SU(4)$ is viewed as a $4 \times 4$ matrix and $e = [e_{ij}]$ is as 
above.

\smallskip

\noindent What we saw in section IX is that there is a surjection,
\begin{equation}
C ({\rm Gr}) \ra C(S^4) \label{eq11}
\end{equation}
while the corresponding symmetry group breaks down to $SO(4)$, the isometry 
group of the 3-sphere from which $S^4$ is obtained by suspension. We shall now show that the algebra $C({\rm Gr})$ is 
noncommutative by constructing explicit surjections,
\begin{equation}
C ({\rm Gr}) \ra C(S_{\theta}^4) \label{eq12}
\end{equation}
whose form is dictated by natural deformations of the 4-sphere similar in 
spirit to the above deformation of $\Tb^2$ to $\Tb_{\theta}^2$.

\noindent We first determine the $C^*$ algebra generated by $M_4 (\Cb)$ and a 
projection $e = e^*$ such that $\left\lgl e - \frac{1}{2} \right\rgl = 0$ as 
above and whose two by two matrix expression is of the form,
\begin{equation}
[ e^{ij} ] = \left[ \matrix{q_{11} &q_{12} \cr q_{21} &q_{22} \cr} \right] 
\label{eq58}
\end{equation}
where each $q_{ij}$ is a $2 \ts 2$ matrix of the form,
\begin{equation}
q = \left[ \matrix{\a &\b \cr -\lambda \b^* &\a^* \cr} \right] \, . \label{eq59}
\end{equation}
where $ \lambda = {\rm exp} 2\pi i \theta$  is a complex number of modulus one, different from -1 for convenience.
Since $e = e^*$, both $q_{11}$ and $q_{22}$ are selfadjoint, moreover since 
$\left\lgl e - \frac{1}{2} \right\rgl = 0$, we can 
find $t = t^*$ such that,
\begin{equation}
q_{11} = \left[ \matrix{t &0 \cr 0 &t \cr} \right] \, , \ q_{22} = \left[ 
\matrix{(1-t) &0 \cr 0 &(1-t) \cr} \right] \, . \label{eq60}
\end{equation}
We let $q_{12} = \left[ \matrix{\a &\b \cr -\lambda \b^* &\a^* \cr} \right]$, we then 
get from $e = e^*$,
\begin{equation}
q_{21} = \left[ \matrix{\a^* &- \bar{\lambda} \b \cr \b^* &\a \cr} \right] \, . \label{eq61}
\end{equation}
We thus see that the commutant $\Bc_{\theta}$ of $M_4 (\Cb)$ is generated by $t,\a,\b$ 
and we first need to find the relations imposed by the equality $e^2 = e$.

\noindent In terms of $e = \left[ \matrix{ t &q \cr q^* &1-t \cr} \right]$, 
the equation $e^2 = e$ means that $t^2 - t + q q^* = 0$, $t^2 - t + q^* q = 
0$ and $[t,q] = 0$. This shows that $t$ commutes with $\a$, $\b$, $\a^*$ and 
$\b^*$ and since $qq^* = q^* q$ is a diagonal matrix
\begin{equation}
\a \a^* = \a^* \a \, , \ \a \b =\lambda  \b \a \, , \ \a^* \b =\bar{\lambda} \b \a^* \, , \ \b 
\b^* = \b^* \b \label{eq62}
\end{equation}
so that the $C^*$ algebra $\Bc_{\theta}$ is not abelian for $ \lambda$ different from 1.
The only further relation is, 
(besides $t = t^*$),
\begin{equation}
\a \a^* + \b \b^* + t^2 - t = 0 \, . \label{eq63}
\end{equation}
We denote by $S^4_{\theta}$ the corresponding noncommutative space, so that $C(S^4_{\theta})=\Bc_{\theta}$.
It is by construction the suspension of the noncommutative 3-sphere $S^3_{\theta}$ whose coordinate algebra is 
generated by $\a$ and $\b$ as above for the special value $t=1/2$.
This noncommutative 3-sphere is related by analytic continuation of the parameter $q$ to the quantum 
group $SU(2)_q$ but the usual theory requires $q$ to be real whereas we need a complex number of 
modulus one which spoils the unitarity of the coproduct.

\noindent We shall now check that the two dimensional component $\left\lgl Ch_1(e) \right\rgl$
automatically vanishes as an element of the (normalized) (b,B)-bicomplex.
\begin{equation}
\left\lgl Ch_{n}(e) \right\rgl = 0 \, , \ 
n = 0,1  . \label{eq67}
\end{equation}
With $q = \left[ \matrix{ \a &\b \cr -\lambda \b^* &\a^* \cr} \right]$, we get,
\begin{eqnarray}
&\lgl Ch_1(e) \rgl &= \Biggl\lgl \left( t - \frac{1}{2} \right) 
\, (dq \, dq^* - dq^* \, dq) \\
&&+ \, q \, (dq^* \, dt - dt \, dq^*) + q^* \, (dt \, dq - dq \, dt ) 
\Biggl\rgl \nonumber
\label{eq69}
\end{eqnarray}
where the expectation in the right hand side is relative to $M_2 (\Cb)$ and we use 
the notation $dx$ instead of the tensor notation.

\noindent The diagonal elements of $\om = dq \, dq^*$ are computed as above,
$$
\om_{11} = d \a \, d\a^* + d\b \, d\b^* \, , \ \om_{22} = d\b^* \, d\b + 
d\a^* \, d\a
$$
while for $\om' = dq^* \, dq$ we get, 
$$
\om'_{11} = d\a^* \, d\a + d\b \, d\b^* \, , \ \om'_{22} = d\b^* \, d\b + d\a 
\, d\a^* \, .
$$
It follows that, since $t$ is diagonal,
\begin{equation}
\left\lgl \left( t - \frac{1}{2} \right) \, (dq \, dq^* - dq^* \, dq) 
\right\rgl = 0 \, . \label{eq70}
\end{equation}
The diagonal elements of $q \, dq^* \, dt = \rho$ are
$$
\rho_{11} = \a \, d\a^* \, dt + \b \, d\b^* \, dt \, , \ \rho_{22} = \b^* \, 
d\b \, dt + \a^* \, d\a \, dt
$$
while for $\rho' = q^* \, dq \, dt$ they are
$$
\rho'_{11} = \a^* \, d\a \, dt + \b \, d\b^* \, dt \, , \ \rho'_{22} = \b^* \, 
d\b \, dt + \a \, d\a^* \, dt \, .
$$
Similarly for $\s = q \, dt \, dq^*$ and $\s' = q^* \, dt \, dq$ one gets the 
required cancellations so that,
\begin{equation}
\left\lgl Ch_1(e) \right\rgl = 0 \, , 
\label{eq71}
\end{equation}
It follows thus that $\left\lgl Ch_2(e) \right\rgl $ is a Hochschild cycle and that 
for any $ds= D^{-1}$ of order $\frac{1}{4}$ commuting with $M_4 (\Cb)$, the index pairing
of $ D$ with $e$ is
\begin{equation}
\hbox{Index} D^+_e  =  \, \int\!\!\!\!\!\!- \ \g \left(e - 
\frac{1}{2}\right) \, [D,e]^4 \, ds^4 \, .\label{eq:(9.21)}
\end{equation}
Exactly as above this shows that $ds$ cannot be of order $\a > 
\frac{1}{4}$ if the index pairing is non zero, and we also get the 
analogue of equation 9-\ref{eq:(9.17)} in the form,
\begin{equation}
\left\langle \left( e - \frac{1}{2} \right) \, [D,e]^4 
\right\rangle = \g \label{eq:(9:36)}
\end{equation}
where $\langle \ \rangle$ is simply the projection on the commutant 
of $M_4 (\Cb)$ in $\Lc (\Hc)$.

\smallskip

\noindent This equation together with (23) implies the integrality of the 4-dimensional volume,
\begin{equation}
\int\!\!\!\!\!- 
\ ds^4 \in \Nb, \label{eq:(9.40)}
\end{equation}
 since it is given by a Fredholm index.
We shall refer to \cite{CLa} for the explicit construction of solutions of (24).
It should be clear to the reader that this amply justifies the clarification to which we turn 
next, of the notion of manifold in Noncommutative Geometry.

\bigskip
\section{Noncommutative Spectral Manifolds}

\noindent In our discussion in section IX of the K-homology fundamental class of a 
manifold we skipped over the nuance between K-homology and KO-homology.
This nuance turns out to be essential in the noncommutative case. Thus to describe 
the fundamental class of a noncommuative space by a spectral triple  $(\Ac ,\Hc ,D)$, 
will require an additional "real structure" on the Hilbert space $\Hc$ given by an antilinear
isometry $J$.
\noindent The anti-linear isometry $J$ is given in Riemannian geometry  by the
charge conjugation operator and in the noncommutative case by the Tomita-Takesaki 
antilinear 
conjugation operator \cite{T}.

\noindent The action of $\Ac$ satisfies the commutation rule,
$[a,b^0] = 0 \quad \fl \, a,b \in \Ac$ where
\begin{equation}
b^0 = J b^* J^{-1} \qquad \fl b \in \Ac
\end{equation}
so $\Hc$ becomes an  $\Ac$-bimodule using the representation of $\Ac \ot \Ac^0$, where $\Ac^0$ is
the opposite algebra, given by,
\begin{equation}
a \ot b^0 \rightarrow aJ b^* J^{-1} \qquad \fl a, b \in \Ac
\end{equation}
\noindent This allows to overcome the main difficulty of the noncommutative case which is that the 
diagonal in the square of the space 
 no longer corresponds to an algebra homomorphism (the map $ x \ot y \ra  \, xy$ 
is 
no longer an algebra homomorphism),

\noindent The {\it fundamental class} of a noncommutative space
is a class $\mu$ in the  $KR$--homology  of the algebra $\Ac \ot \Ac^0$
equipped with the  involution
\begin{equation}
\tau (x \ot y^0) = y^* \ot (x^*)^0 \qquad \fl \, x,y \in
\Ac                                                     \label{eq:(1.15)}
\end{equation}
where $\Ac^0$ denotes the algebra opposite to $\Ac$.  
The  $KR$-homology cycle representing $\mu$ is given by a spectral
triple, as above, equipped with an anti-linear
isometry $J$ on  $\Hc$ which implements the involution $\tau$,
\begin{equation}
J w J^{-1} = \tau (w) \qquad \fl \, w \in \Ac \ot \Ac^0
\, ,                                                    \label{eq:(1.16)}
\end{equation}
$KR$-homology (\cite{[18]} \cite{[At]}) is periodic with period $8$ and the  
dimension modulo $8$ is specified by the
following commutation rules. One has $J^2 = \ve$, $JD = \ve' DJ$, $J\g = \ve'' \g 
J$
where  $\ve ,\ve' ,\ve'' \in \{ -1,1\}$ and with $n$ the dimension modulo 8,

\bigskip

\begin{center}
\begin{tabular}
{|c| r r r r r r r r|} \hline
{\bf n }&0 &1 &2 &3 &4 &5 &6 &7 \\ \hline
\hline
$\ve$  &1 & 1&-1&-1&-1&-1& 1&1 \\
$\ve'$ &1 &-1&1 &1 &1 &-1& 1&1 \\
$\ve''$&1 &{}&-1&{}&1 &{}&-1&{} \\  \hline
\end{tabular}
\end{center}

\bigskip

 \noindent The class $\mu$ specifies only the stable homotopy class of 
the
spectral triple $(\Ac ,\Hc ,D)$ equipped with the isometry $J$ (and $\Zb
/2$--grading $\g$ if $n$ is even). The non-triviality of this homotopy
class shows up in the intersection form
\begin{equation}
K_* (\Ac) \ts K_* (\Ac) \ra \Zb                             \label{eq:(4.104)}
\end{equation}
which is  obtained from the Fredholm index of $D$ with   coefficients
in $K_* (\Ac \ot \Ac^0)$. Note that it is  defined without using
the diagonal map $m:\Ac \ot \Ac \ra \Ac$, which is not a homomorphism
in the noncommutative case. This form is quadratic
or symplectic according to the value of $n$ modulo $8$.

 \smallskip

\noindent The  Kasparov intersection product \cite{[18]} allows to formulate
the Poincar\'e duality in terms of the invertibility  of $\mu$,
\begin{equation}
\ex \, \b \in KR_n (\Ac^0 \ot \Ac) \ , \ \b \ot_{\Ac}
\mu = \id_{\Ac^0} \ , \ \mu \ot_{\Ac^0} \b = \id_{\Ac}
\, .                             \label{eq:(4.103)}
\end{equation}

\noindent It implies the isomorphism $K_* (\Ac) \
\stackrel{\cap \mu} \longra K^* (\Ac)$. 

\smallskip

\noindent  The condition that D is an operator of order one becomes
\begin{equation}
[[D,a],b^0] = 0 \qquad \fl \, a,b \in \Ac \, .                             \label{eq:(4.102)}
\end{equation}
(Notice that  since $a$ and  $b^0$ commute this condition
is equivalent to  $[[D,a^0],b]=0 \quad \fl \, a,b \in \Ac$.)

\smallskip

\noindent  One can show that the von Neumann algebra $\Ac''$ generated by $\Ac$ in
$\Hc$ is automatically finite and hyperfinite and there is a complete
list of such algebras up to isomorphism.  The algebra $\Ac$
is stable under smooth  functional calculus in its norm closure $A =
\bar{\Ac}$ so that $K_j (\Ac) \sm K_j (A)$, i.e. $K_j (\Ac)$ depends
only on the underlying topology (defined by the $C^*$ algebra $A$).
The integer $\chi = \lgl \mu ,\b \rgl \in \Zb$ gives the Euler characteristic in 
the 
form
\begin{equation}
\chi = \Rang K_0 (\Ac) - \Rang K_1 (\Ac)                               \label{eq:(4.101)}
\end{equation}
and the general operator theoretic index formula of section 13 below, gives a local 
formula 
for $\chi$.

\noindent We gave in \cite{Co3} the necessary and sufficient conditions that a spectral triple
(with real structure $J$) should fulfill in order to come from an ordinary compact Riemannian spin manifold.
These conditions extend in a straightforward manner  to the noncommutative case (\cite{Co3}).
To appreciate the richness of examples which fulfill them we shall just quote the following result (\cite{CLa}),

\bigskip
\noindent {\bf Theorem 6.}
{\it Let $M$ be a compact Riemannian spin manifold. Then if the isometry 
group of $M$ has rank $r \geq 2$, $M$ admits a non-trivial one parameter
isospectral deformation to noncommutative geometries $M_{\theta}$.} 

\bigskip
\noindent The group $\Aut^+ (\Ac)$ of automorphisms $\alpha$ 
of the involutive algebra $\Ac$, which are implemented by a unitary 
operator $U$ in $\Hc$ commuting with $J$,
\begin{equation}
\alpha (x) = U \, x \, U^{-1} \qquad \fl \, x \in
\Ac \ ,                                 \label{eq:(4.100)}
\end{equation}
plays the role of the group $\Diff^+ (M)$ of
diffeomorphisms preserving the K-homology fundamental class for a manifold $M$.

\smallskip

\noindent In the general noncommutative case, parallel to
the normal subgroup  $\Int \Ac \sbs \Aut \Ac$  of inner automorphisms of
$\Ac$,
\begin{equation}
\a (f) = ufu^* \qquad \fl \, f \in \Ac
                                                        \label{eq:(4.8)}
\end{equation}
where $u$ is a unitary element of $\Ac$ (i.e.
$uu^* = u^* u =1$), there exists a natural foliation
of the space of spectral geometries on
$\Ac$  by equivalence classes of
{\it inner deformations}
of a given geometry. To understand how they arise we need to understand 
how to transfer a given spectral geometry to a Morita
equivalent algebra. Given a spectral triple $(\Ac,\Hc ,D)$ and the Morita 
equivalence \cite{[Ri1]} between $\Ac$
and an algebra $\Bc$ where
\begin{equation}
\Bc = \End_{\Ac} (\Ec)
                                                        \label{eq:(4.16)}
\end{equation}
where  $\Ec$ is a finite,  projective, hermitian
right  $\Ac$--module, one gets a
spectral triple on  $\Bc$ by  the choice of  a {\it
hermitian connection} on  $\Ec$. Such a connection  $\nb$
is a linear map $\nb : \Ec \ra \Ec \ot_{\Ac} \Om_D^1$
satisfying the rules (\cite{[Co]})
\begin{equation}
\nb (\xi a) = (\nb \xi) a + \xi \ot da \qquad \fl \, \xi
\in \Ec \ , \ a\in \Ac
                                                        \label{eq:(4.17)}
\end{equation}
\begin{equation}
(\xi , \nb \eta) - (\nb \xi ,\eta) = d(\xi ,\eta) \qquad
\fl \, \xi ,\eta \in \Ec
                                                        \label{eq:(4.18)}
\end{equation}
where  $da = [D,a]$ and where $\Om_D^1 \sbs \Lc (\Hc)$ is
the $\Ac$--bimodule of operators of the form
\begin{equation}
A = \Si \, a_i [D,b_i] \ , \ a_i , b_i \in \Ac \, .
                                                        \label{eq:(4.10)}
\end{equation}

\smallskip

\noindent Any algebra $\Ac$ is Morita equivalent to itself (with $\Ec =
\Ac$) and  when one applies the above construction in the above context one gets 
the 
inner
deformations of the spectral geometry.

\noindent Such a deformation is obtained by
the following  formula (with suitable signs depending on the dimension mod 8) 
without modifying  neither the
representation of $\Ac$ in $\Hc$ nor the anti-linear isometry $J$
\begin{equation}
D\ra D+A+JAJ^{-1}
                                                        \label{eq:(4.9)}
\end{equation}
where $A=A^*$ is an arbitrary selfadjoint operator of the form \ref{eq:(4.10)}.
The action of the group
$\Int (\Ac)$ on the spectral geometries 
is simply the following gauge transformation of $A$
\begin{equation}
\g_u (A) = u[D,u^*] + uAu^* \, .
                                                        \label{eq:(4.11)}
\end{equation}
The required unitary equivalence is
implemented by the following representation of
the unitary group of
$\Ac$ in  $\Hc$,
\begin{equation}
u \ra uJuJ^{-1} = u(u^*)^0 \, .
                                                        \label{eq:(4.12)}
\end{equation}
The transformation (15) is the identity 
in the usual Riemannian case. To get a nontrivial example it suffices
to consider the product of a Riemannian triple by the 
unique spectral
geometry on the finite-dimensional algebra $\Ac_F = M_N (\Cb)$ of $N\ts N$
matrices on $\Cb$, $N\geq 2$. One then has $\Ac = C^{\ify} (M)
\ot \Ac_F$, $\Int (\Ac) = C^{\ify} (M,PSU(N))$ and inner deformations of
the geometry are parameterized by the gauge potentials for the gauge
theory of the group $SU(N)$.  The space of pure states of the algebra
$\Ac$, $P(\Ac)$, is the product $P = M\ts P_{N-1} (\Cb)$ and the metric
on $P(\Ac)$ determined by the formula (9.13) depends on the
gauge potential $A$. It coincides with the Carnot metric \cite{[G]} on $P$
defined by the horizontal distribution given by the connection
associated to $A$.  The group $\Aut (\Ac)$ of
automorphisms of $\Ac$ is the following semi--direct product
\begin{equation}
\Aut (\Ac) = \Uc \semi
\Diff^+ (M)
                                                        \label{eq:(4.13)}
\end{equation}
of the local gauge transformation group $\Int (\Ac)$ by the group of
diffeomorphisms.
\bigskip

 \section{Test with space-time}

What we have done so far is to stretch the usual framework of ordinary geometry 
beyond its
commutative
restrictions (set theoretic restrictions) and of course now it's not
perhaps a
bad idea to test it with what we know about physics and to try to find a better model 
of 
space-time 
within this new framework. The
best
way is to start with the hard core information one has from physics and that
can be summarized by a Lagrangian. This Lagrangian is the Einstein
Lagrangian
plus the standard model Lagrangian. I am not going to write it down, it's a
very complicated expression since just the  standard  model  Lagrangian comprises 
five types 
of 
terms. But one can start understanding something
by looking at the symmetry group of this Lagrangian. Now, if it were just  the 
Einstein 
theory, 
the symmetry group of the Lagrangian
 would just be, by the equivalence
principle, the diffeomorphism group of the space-time manifold. But
 because of the standard model piece the symmetry group of this Lagrangian
 is not just the diffeomorphism group, because 
the gauge theory has another huge  symmetry group which is the
group
of maps from the manifold to the small gauge group, namely $U_1\times
SU_2\times SU_3$ as far as we know. Thus, the symmetry group $G$ of the full
Lagrangian is neither the diffeomorphism group nor the group of gauge
transformations of second kind nor their product, but it is their semi-direct
product. It is exactly like what happens with the Poincar\'e group where you
have
translations and Lorentz transformations, so it is the semi-direct product of
these two subgroups. Now we can ask a very simple question: would there be some
space $X$ so that this group $G$  would be equal to $\Diff (X)$? If such a
space
would exist, then we would have some chance  to actually geometrize
completely
the theory, namely  to be able to say that it's pure gravity on the space
$X$.
Now, if you look for the space $X$ among ordinary manifolds, you have no chance
 since by a result of John Mather the diffeomorphism  group of a (connected) manifold 
is a 
simple
group. A simple group cannot have a nontrivial normal subgroup, so you
cannot have this structure of semi-direct product. 

\noindent However, we can use our dictionary,
and in
this dictionary if we browse through it, we find that what corresponds to
diffeomorphisms for a non commutative space  is just the group $\Aut^+ (\Ac)$ of
automorphisms
of the algebra of coordinates $\Ac$, which preserve the fundamental class in 
$K$-homology, as described above in section XI.

\noindent Now there is a
beautiful fact which is that when an algebra is not commutative, then among its
automorphisms there are very trivial ones, there are automorphisms  which
are there for free, I mean the  inner ones, which associate to an element $x$ of
the
algebra the element $uxu^{-1}$. Of course $uxu^{-1}$ is not, in general equal to $x$ 
because 
the
algebra is not commutative, and these automorphisms form a normal subgroup
of
the group of automorphisms. Thus you see that the group $\Aut^+ (\Ac)$ 
 has the same type of structure, namely it has a normal subgroup of internal
automorphisms and it has a quotient. Now it turns out that there is one very
natural non commutative algebra $\Ac$ whose group of
internal automorphisms corresponds to the group of gauge
transformations and the quotient $\Aut^+ (\Ac) / \Int (\Ac)$ corresponds exactly to 
diffeomorphisms \cite{Schu}. 
It is amusing that the physics vocabulary is actually the same as
the
mathematical  vocabulary. Namely in physics you talk about internal symmetries
and in mathematics you talk about inner automorphisms, you could call them
internal automorphisms. Now the
corresponding
space is a product $ M\ts F$ of an ordinary manifold $M$ by a finite noncommutative
space $F$. The corresponding algebra $\Ac_F$ is the direct sum of the algebras
$\Cb$ , $ \ \Hb$ (the quaternions), and $M_3 (\Cb)$ of  $3\ts 3$ complex
matrices.

\noindent The algebra $\Ac_F$ corresponds to a {\it finite} space where
the standard model fermions and the Yukawa parameters (masses of fermions
and mixing matrix of Kobayashi Maskawa) determine the spectral
geometry in the following manner. The Hilbert space is finite-dimensional
and admits the set of elementary fermions as a basis.
For example for the first generation of quarks, this  set  is
\begin{equation}
u_L , u_R , d_L , d_R,\bar{u}_L , \bar{u}_R , \bar{d}_L , \bar{d}_R
\, .
                                                        \label{eq:(5.4)}
\end{equation}
The algebra $\Ac_F$ admits a natural representation in $\Hc_F$
(see \cite{[Coo2]}) and the Yukawa coupling matrix $Y$ determines the operator $D$.

The detailed structure of
$Y$ (and in particular the fact that color is not broken) allows to check the
axioms of noncommutative geometry.
\smallskip

\noindent The next step consists in the computation of  internal deformations 
\begin{equation}
D\ra D+A+JAJ^{-1}
                                                        \label{eq:(4.9)}
\end{equation}
(cf. section XI),
 of the product geometry $M\ts F$ where $M$ is
a $4$--dimensional Riemannian spin manifold. The computation gives the
standard
model gauge bosons $\g , W^{\pm} ,Z$, the  eight gluons and
the  Higgs fields $\vp$ with
accurate  quantum numbers.

\noindent Now the next question that comes about is how do we
recover the original action functional which contained both the Einstein-Hilbert term 
as well as the standard model ? The answer is very simple: the Fermionic part of 
this action is there from the start and one recovers the bosonic part as follows.
Both the Hilbert--Einstein action
functional for the Riemannian metric, the Yang--Mills action for
the  vector potentials, the self interaction and the minimal coupling for the Higgs 
fields 
all
 appear with the correct signs in the
asymptotic expansion for large $\Lb$ of the number $N(\Lb)$ of
eigenvalues of $D$ which are $\leq \Lb$ (cf. \cite{[C-C]}), 
\begin{equation}
N(\Lb) = \# \ \hbox{eigenvalues of $D$ in} \ [-\Lb,\Lb] .
\end{equation}
Exactly as above, this step function $N(\Lb)$ is the superposition of two terms,
$$
 N(\Lb)= \lgl N(\Lb) \rgl +N_{\rm osc} (\Lb).
$$
 The oscillatory part $N_{\rm osc} (\Lb)$ 
is the same as for a random matrix, governed by the statistic
dictated by the symmetries of the system and does not concern us here.
The average part $ \lgl N(\Lb) \rgl $ is computed by a semiclassical approximation 
from local expressions
 involving the familiar heat equation expansion and delivers the correct terms. 
 We showed above in section IX,  that if one studies natural 
presentations 
of the 
algebra generated by $\Ac$ and $D$ one naturally gets 
 only metrics with a fixed volume form so that the bothering 
cosmological term does not enter in the variational
equations associated to the spectral action $ \lgl N(\Lb) \rgl$. It is tempting to 
speculate that the phenomenological Lagrangian
 of physics, combining matter and gravity appears from the solution of an 
extremely 
simple 
operator theoretic equation along the lines described above in sections IX and X. 
  
\vglue 1cm

\bigskip
\section{Operator theoretic Index Formula}

The power of the general theory comes from 
deeper general theorems such as the local computation of the analogue of
Pontrjagin classes: {\it i.e.} of the components of the cyclic
cocycle  which is the Chern character of the K-homology class of $D$ and
which
make sense in general. This result allows, using the 
infinitesimal 
calculus, 
to go from local to global in the general framework of spectral triples $(\Ac 
,\Hc,D)$. 
\smallskip

 \noindent The  Fredholm index of the operator  $D$
determines (in the odd case) an  additive map  $K_1 (\Ac) \
\stackrel{\vp}{\ra} \Zb$ given by the equality
\begin{equation}
\vp ([u]) = \Index \, (PuP) \ , \ u\in GL_1 (\Ac)
                                                \label{eq:(3.3)}
\end{equation}
where $P$ is the projector  $P = {1+F \over 2}$, $F = \Sign \, (D)$.

\smallskip
 
 \noindent This map is computed by the pairing
of  $K_1 (\Ac)$ with the following  cyclic cocycle
\begin{equation}
\tau (a^0 ,\ldots ,a^n) = \Trace \, (a^0 [F,a^1] \ldots
[F,a^n]) \qquad \fl \, a^j \in \Ac
                                                \label{eq:(3.4)}
\end{equation}
where $F= \hbox{Sign} \ D$ and  we assume that the dimension $p$ of our space is 
finite,
which means that $(D+i)^{-1}$ is of order $1/p$, also $n\geq p$ is an odd integer. 
There are similar formulas involving the grading $\gamma$ in the even case, and it 
is quite satisfactory (\cite{theta} \cite{JLO}) that both cyclic cohomology and 
the 
chern Character formula adapt to the infinite dimensional case
 in which the only hypothesis is that $\hbox{exp}(-D^{2})$ is a trace class operator.

\smallskip

 \noindent It is difficult to compute the cocycle $\tau$ in general because the formula
(\ref{eq:(3.4)}) involves the ordinary trace instead of the local trace
${\int \!\!\!\!\! -}$ and it is crucial to obtain a local form of the above 
cocycle.

\smallskip

 \noindent This problem is solved by a general formula \cite{[C-M2]} which we now 
describe.

\noindent Let us make the following regularity hypothesis
on  $(\Ac ,\Hc ,D)$
\begin{equation}
a \ \hbox{and } \ [D,a] \ \in \ \cap \Dom \d^k , \ \fl \, a\in
\Ac                                                     \label{eq:(3.1)}
\end{equation}
where $\d$ is the derivation $\d(T) = [\vert D \vert,T]$ for any operator $T$.

\smallskip

\noindent We let $\Bc$ denote the algebra generated by  $\d^k (a)$,
$\d^k ([D,a])$. 
The usual notion of {\it dimension}  of a space is replaced by the
{\it dimension spectrum } which is a subset of $\Cb$.
\noindent The precise definition of the dimension spectrum is the subset
$\Si \sbs \Cb$ of singularities of the analytic functions
\begin{equation}
\z_b (z) = \Trace \, (b \vert D \vert^{-z}) \qquad \Re z
> p \ , \ b\in \Bc \, .
                                                        \label{eq:(3.2)}
\end{equation}
The dimension spectrum of a manifold $M$ is the set $\{
0,1,\ldots ,n\}$, $n=\dim M$; it is simple.  Multiplicities appear for
singular manifolds. Cantor sets provide examples of  complex points $z
\notin \Rb$ in the dimension  spectrum.

\noindent We assume  that $\Si$ is discrete and simple, i.e. that
$\z_b$ can be extended to  $\Cb / \Si$  with simple poles in  $\Si$.

\smallskip

\noindent We refer to  \cite{[C-M2]} for the case of a
spectrum with multiplicities.  Let $(\Ac ,\Hc ,D)$ be a spectral triple satisfying 
the hypothesis
(\ref{eq:(3.1)}) and (\ref{eq:(3.2)}). The local index theorem is the following,  
\cite{[C-M2]}:

\bigskip

\noindent {\bf Theorem 7.} 
\begin{enumerate}
\item {\it The equality
\[
{\int \!\!\!\!\!\! -} P = \Res_{z=0} \, \Trace \, (P\vert
D\vert^{-z})
\]
defines a  trace on the algebra generated by
$\Ac$, $[D,\Ac]$ and $\vert D \vert^z$, where}
$z\in \Cb$.

\item
{\it There is  only  a finite number of non--zero terms in the following formula which
defines the odd components $(\vp_n)_{n=1,3,\ldots}$
of a cocycle in the bicomplex $(b,B)$ of $\Ac$,
 \[
\vp_n (a^0 ,\ldots ,a^n) = \sum_k c_{n,k}
{\int \!\!\!\!\!\! -} a^0 [D,a^1]^{(k_1)} \ldots
[D,a^n]^{(k_n)} \, \vert D \vert^{-n-2\vert k\vert}
\qquad \fl \, a^j \in \Ac
 \]
where the following notations are used:
$T^{(k)} = \nb^k (T)$ and $\nb (T) = D^2 T
- TD^2$, $k$ is a multi-index}, $\vert k \vert = k_1 +\ldots + k_n$,
\[
c_{n,k} =
(-1)^{\vert k \vert} \, \sqrt{2i} (k_1 ! \ldots k_n
!)^{-1} \, ((k_1 +1) \ldots (k_1 + k_2 + \ldots + k_n
+n))^{-1} \, \G \left( \vert k \vert + {n\over 2}
\right).
\]

\item {\it The pairing of the cyclic cohomology class $(\vp_n) \in HC^* (\Ac)$
with $K_1 (\Ac)$ gives the Fredholm index of $D$ with coefficients in}
$K_1(\Ac)$.
\end{enumerate}

\bigskip

\noindent For the normalization of the pairing
between $HC^*$ and  $K(\Ac)$ see \cite{[Co]}.
In the even case, i.e. when $\Hc$ is $\Zb /2$ graded by $\g$,
 \[
\g = \g^*, \ \ \g^2 =1, \ \ \g a = a \g \quad \fl \, a \in \Ac, \ \g
D = -D\g,
\]
there is an analogous formula for a cocycle
$(\vp_n)$, $n$ even, which gives the Fredholm index of $D$
with  coefficients in  $K_0$. However,
$\vp_0$ is not expressed  in terms of the residue
${\int \!\!\!\!\! -}$ because it is not local for a
finite dimensional  $\Hc$.

\bigskip
\section{Diffeomorphism invariant Geometry}

\noindent The power of the above operator theoretic local trace formula lies in its 
generality and in the existence of really new geometric examples to which it applies.

\noindent In this section we shall explain how the transverse structure of 
foliations is described by a spectral triple $(\Ac 
,\Hc,D)$ with simple dimension spectrum. This allows moreover to give a precise meaning to 
diffeomorphism invariant 
geometry on a manifold M, by the construction of a spectral triple $(\Ac ,\Hc ,D)$ 
where the 
algebra $\Ac$ is the crossed product of the algebra of smooth functions on the 
finite 
dimensional bundle $P$ of metrics on M by the natural action of the diffeomorphism 
group of 
M. While ordinary geometric constructions are "covariant" with respect to diffeomorphisms,
our construction (\cite{CM}) is "invariant" inasmuch as the algebra now incorporates the full
group of diffeomorphisms and the metrics involved are canonical.

\noindent The operator $D$ is an hypoelliptic operator (\cite{[H-S]}) which is directly associated to 
the 
reduction 
of the structure group of the manifold $P$ to a group of triangular matrices whose 
diagonal blocks are orthogonal. By construction the fiber of $P 
\stackrel{\pi}{\ra} 
M$ is the 
quotient $F^+/SO(n)$ of the $GL^+(n)$--principal bundle $F^+$ of oriented frames on $M$ by the 
action of the
orthogonal group $SO(n) \sbs GL^+(n)$. The space $P$ admits a canonical foliation: 
the 
vertical 
foliation $V \sbs TP$, $V=\Ker \pi_*$ and   on the fibers $V$ and on  $N = (TP)/V$ 
the 
following Euclidean structures. A choice of $GL^+(n)$--invariant Riemannian metric
on  $GL^+(n)/SO(n)$ determines a metric on  $V$. The metric on $N$ is defined 
tautologically: 
for every $p\in P$ one has a metric on  $T_{\pi (p)} (M)$ which is isomorphic to  
$N_p$  by 
$\pi_*$.

\noindent We first consider the hypoelliptic 
signature operator  $Q$  on $F^+$. It is not a scalar operator but it acts in the 
tensor product
\begin{equation}
{\Hc}_0= L^2 (F^+ , v) \ot E 
     \label{eq:(12.2)}
\end{equation}
where $E$ is a finite dimensional representation of $SO(n)$ specifically given 
by
\begin{equation}
E = \wedge \, P_n \ot  \wedge \, \Rb^n \ , \ P_n = S^2 \, \Rb^n \, . 
  \label{eq:(12.3)}
\end{equation}
The operator $Q$ is the graded sum,
\begin{equation}
Q = (d_V^* \, d_V - d_V \, d_V^*) \op (d_H + d_H^*) 
 \label{eq:(12.4)}
\end{equation}
where the horizontal (resp. vertical) differentiation $d_H$ (resp. $d_V$) is a 
matrix in the horizontal and vertical vector fields ${\bf X}_i$ and  ${\bf Y}_{\ell}^k$  as well as their adjoints 
(which 
also involve scalars). When $n$ is equal to 1 or 2 modulo 4 one has to replace
$F^+$ by its product by $ S^1$ so that the dimension of the vertical fiber is even (it is then $1+{n(n+1) \over 2}$ )
and the vertical signature operator makes sense.  The longitudinal part is not elliptic but only 
transversally elliptic with respect to the action of $SO(n)$.  
 Thus to get an 
hypoelliptic operator one restricts $Q$ to the Hilbert space,
\begin{equation}
{\Hc} = (L^2 (F^+ , v) \ot E)^{SO(n)} 
\label{eq:(12.5)}
\end{equation}
and one takes the following algebra $\Ac$,
\begin{equation}
{\Ac}= C_c^{\ify} (P) \semi \hbox{Diff}^+ \ , \ P = F^+ / SO(n) \, . 
\label{eq:(13.6)}
\end{equation}
Let us note that the operator $Q$ is in fact the image under the right regular
 representation of the affine group $G_{affine}$ of a (matrix valued) hypoelliptic symmetric element in the
envelopping algebra ${\Uc} (G_{affine})$.
By an easy adaptation of a theorem of Nelson and Stinespring, it then follows that $Q$ is essentially selfadjoint (with core any dense
 $G_{affine}$-invariant subspace of the space of $C^{\ify}$ vectors of the right regular representation of $G_{affine}$).

\medskip
\noindent {\bf Theorem 8.} \cite{CM}
 {\it 
Let  $\Ac$ be
the crossed product $C_c^{\ify} (P) \semi \hbox{Diff}^+$ acting in $\Hc$ as above.
\begin{enumerate}

\item The equality $D\vert D\vert =Q$ defines a spectral triple $(\Ac ,\Hc ,D)$ 
 which satisfies the hypotheses
of theorem 7; its dimension spectrum is simple and given by $\Si =
\{ 0,1,...,2n+\frac{n(n+1)}{2}\}$.

\item The cocycle $\vp_j$ given by the local index formula
 (theorem 7) is  the image by the characteristic map
of a universal Gelfand-Fuchs cohomology class.
\end{enumerate}
}

\noindent The equality $D\vert D\vert =Q$ defining $D$ while $Q$ is a differential operator
of second order, is characteristic of "quartic" geometries.

\noindent The computation of the local index formula for diffeomorphism invariant 
geometry 
\cite{CM} was quite complicated even in the case of codimension 1
foliations: there  were innumerable terms to be computed; this could be done by
hand, by 3 weeks of eight hours per day tedious computations, but it was of course 
hopeless to proceed by direct computations in the general case. Henri and I 
finally 
found how to get the answer for the general
case after discovering that the computation generated a Hopf algebra $\Hc(n)$ 
which 
only depends on n= codimension of
the foliation, and which allows to organize the computation provided cyclic 
cohomology is suitably adapted to Hopf algebras as in the next section.

\noindent The Hopf algebra $\Hc(n)$ only depends upon the integer $n$ and is neither
 commutative nor cocommutative. We proved in \cite{CM} that it is isomorphic to the 
bicrossed product Hopf algebra  
(\cite{[Kac]}, \cite{[Baaj-S]}, \cite{[Majid]}) associated to the following pair of 
subgroups of $G=\hbox{Diff}(\Rb^n)$. 

\noindent We let $G_1 \sbs G$ be the subgroup of affine diffeomorphisms,
\begin{equation}
k(x) = Ax + b \qquad \fl \, x \in \Rb^n 
\end{equation}
and we let $G_2 \sbs G$ be the subgroup,
\begin{equation}
\vp \in G \, , \ \vp (0) = 0 \, , \ \vp' (0) = 1 \, . 
\end{equation}
Given $\vp \in G$ it has a unique decomposition $\vp = k \, \psi$ where $k 
\in G_1$, $\psi \in G_2$ which allows to perform the bicrossed product construction.

\bigskip
\section{Characteristic classes for actions of Hopf algebras}

\noindent Hopf algebras arise very naturally from their actions on
noncommutative algebras \cite{[M]}. Given an algebra $A$, an 
action of the Hopf algebra $\Hc$ on $A$ is given by a linear map,
\begin{equation}
\Hc \ot A \ra A, \quad h \ot a \ra h(a) 
\end{equation}
satisfying $h_1 (h_2  a) = (h_1  h_2) (a)$, $\fl  h_i \in {\Hc}$,
$a \in A$ and
\begin{equation}\label{ba}
h(ab) = \sum  h_{(1)}  (a)  h_{(2)}  (b)  \qquad \fl  a,b  \in A,
h \in {\Hc}.
\end{equation}
where the coproduct of $h$ is,
\begin{equation}\label{bb}
\D(h)=  \sum  h_{(1)}  \ot  h_{(2)} 
\end{equation}
In concrete examples, the algebra $A$ appears first, together with
linear maps $A \ra A$ satisfying a relation of the form (\ref{ba}) which dictates 
the Hopf algebra structure. This is exactly what occured in the above example (see 
\cite{CM} for the description of $\Hc(n)$ and its relation with Diff($\Rb^n$)).

\noindent The theory of characteristic classes for actions of $\Hc$
extends the construction \cite{C3} of cyclic cocycles from 
a Lie algebra of derivations of a $C^*$ algebra $A$, 
together with an \textit{invariant trace} $\tau$ on $A$.
 
\noindent This theory was developped in \cite{CM} in order to solve the above 
computational 
problem for diffeomorphism invariant geometry but it was shown in \cite{CM2} that 
the correct framework 
for the cyclic cohomology of Hopf algebras is that of modular pairs in involution. 
It is 
quite satisfactory that exactly the same structure
 emerged from the analysis of locally compact quantum groups.
The resulting cyclic cohomology appears to
be the natural candidate for the analogue of Lie algebra cohomology in
the context of Hopf algebras. We fix a group-like element $\sigma$ and a
character $\delta$ of $\Hc$ with $\delta(\sigma)=1$. They 
will play the role of the module of locally compact groups.

\noindent We then introduce the twisted antipode,
\begin{equation}\label{bc}
\wt S (y) = \sum  \delta (y_{(1)})  S (y_{(2)}) \ , \ y \in {\Hc}  , \
\D y = \sum  y_{(1)} \ot y_{(2)}.
\end{equation}
 \noindent We shall say that the modular
pair
($\sigma$, $\delta$) is in involution if the ($\sigma$, $\delta$)-twisted antipode 
is an involution,
\begin{equation}\label{ch} 
(\sigma^{-1}\wt{S})^2 = I.
\end{equation}
  We associate a cyclic complex (in fact a 
$\Lambda $-module, where $ \Lambda$ is the cyclic category), 
to any Hopf algebra together with a modular pair in involution.
More precisely the following graded vector space $\Hc_{(\delta,\sigma)}^{\natural} = \{ \Hc^{\ot n} \}_{n \geq 1}$ equipped 
with the operators given by the following formulas (\ref{ce})--(\ref{cg}) defines 
a 
module 
over the cyclic category $\Lambda$.
First, by transposing the standard simplicial operators underlying the 
Hochschild homology complex of an algebra,
one associates to $\Hc$, viewed only as a coalgebra, the 
natural cosimplicial  module $\{ \Hc^{\ot n} \}_{n \geq 1}$, with face operators 
$\delta_i: 
\Hc^{\ot n-1} \ra \Hc^{\ot n}$,
\begin{eqnarray}\label{ce}
&&\delta_0 (h^1 \ot \ldots \ot h^{n-1}) = 1 \ot h^1 
\ot \ldots \ot h^{n-1} \nonumber \\
&& \nonumber \\
&&\delta_j (h^1 \ot \ldots \ot h^{n-1}) = h^1 \ot \ldots \ot \D h^j \ot 
\ldots \ot h^n,\  \fl  1 \leq j \leq n-1, \nonumber\\
&&  \\
&&\delta_n (h^1 \ot \ldots \ot h^{n-1}) = h^1 \ot \ldots \ot h^{n-1}
\ot \sigma \nonumber
\end{eqnarray}
and degeneracy operators $\s_i : \Hc^{\ot n+1} \ra \Hc^{\ot n}$,
\begin{equation}\label{cf}
\s_i (h^1 \ot \ldots \ot h^{n+1}) = h^1 \ot \ldots \ot \ve (h^{i+1}) 
\ot \ldots \ot h^{n+1} \ , \ 0 \leq i \leq n.
\end{equation}

\noindent The remaining two essential features of a Hopf algebra 
--\textit{product} and \textit{antipode} -- are brought into play, to
define the \textit{cyclic operators} $\tau_n : \Hc^{\ot n} \ra \Hc^{\ot n}$,
\begin{equation}\label{cg} 
\tau_n (h^1 \ot \ldots \ot h^n) = (\D^{n-1}  \wt S (h^1)) \cdot h^2 \ot \ldots \ot 
h^n \ot 
\sigma.
\end{equation}
The theory of characteristic classes applies to actions of the Hopf algebra
on an algebra endowed with a $\delta$-invariant $\sigma$-trace.
A linear form $\tau$ on $A$ is a $\sigma$-trace
under the action of $\Hc$ iff one has,
\begin{equation}
\tau (ab) = \tau (b \sigma (a)) \qquad \fl  a,b  \in A  .  \label{mod}
\end{equation}
A $\sigma$-trace $\tau$ on $A$ is $\delta$-invariant
under the action of $\Hc$ iff 
\begin{equation}
\tau (h(a)b) = \tau (a  \wt S (h)(b)) \qquad \fl  a,b  \in A  , \ h
\in 
{\Hc}. 
\end{equation}
\noindent  Note that equation (\ref{mod}) is an excellent guide in order to construct 
Hopf algebra actions, since by the modular theory any positive linear 
functional $\tau$ on an algebra $A$ gives rise to an (unbounded) automorphism
$\sigma$ of its weak closure fulfilling equation (\ref{mod}).

\noindent The theory of characteristic classes for actions of Hopf algebras
is governed by the following 
general result:

\medskip
\noindent {\bf Theorem 9.}
 (\cite{CM2}) {\it 
Let $\Hc$ be a Hopf algebra endowed
with a modular pair in involution
Then $\Hc_{\delta, \sigma}^{\natural} = \{ \Hc^{\ot n} \}_{n \geq 1}$ equipped 
with the operators given by (\ref{ce})--(\ref{cg}) defines a module over the
cyclic category $\Lambda$.
Let $\Hc$ act on an algebra $A$ endowed with a
$\delta$-invariant $\sigma$-trace $\tau$
, then the following defines a canonical map from $HC^*_{\delta, \sigma}
({\Hc})$ to 
$HC^* (A)$,
$$
\matrix{
\g (h^1 \ot \ldots \ot h^n) \in C^n (A)  , \ \g (h^1 \ot \ldots \ot
h^n) 
(x^0 , \ldots , x^n) = \cr
\cr
\tau (x^0  h^1 (x^1) \ldots h^n (x^n)). \cr
}
$$
}
\noindent We refer to \cite{CM2} for the discussion of the remarkable agreement
of this theory with the standard theory of quantum groups and their locally compact
versions.

\bigskip
\section{Hopf algebras, Renormalization and the 
Riemann-Hilbert problem}

 \noindent We describe in this section our joint work with Dirk Kreimer.
Perturbative renormalization is by far the most successful
technique for computing physical quantities in quantum field
theory. It is well known for instance that it accurately predicts
the first ten decimal places of the anomalous magnetic moment of
the electron.

\smallskip

\noindent The physical motivation behind the renormalization technique
is quite clear and goes back to the concept of effective mass  in
nineteen century hydrodynamics. To appreciate it, one should dive
under water with a ping-pong ball and start applying
Newton's law,
\begin{equation}
F = m \, a \label{ball}
\end{equation}
to compute the initial acceleration of the ball B when we let it loose
(at zero speed relative to the water). If one naively applies \ref{ball},
one finds (see the QFT course by Sidney Coleman) an unrealistic initial
acceleration of about 20g! In fact as explained in loc. cit. due to the 
interaction of B with the surrounding field of water, the inertial mass $m$ 
involved in \ref{ball} is
not the bare mass $m_0$ of B but is modified to
\begin{equation}
m = m_0 + {\textstyle \frac{1}{2}} \, M 
\end{equation}
where $M$ is the mass of the water occupied by B.

\smallskip

\noindent It follows for instance that the initial acceleration $a$ of
B is given, using the Archimedean law, by
\begin{equation}
-(M-m_0) g = \left( m_0 + {\textstyle \frac{1}{2}} \, M \right) a 
\end{equation}
and is always of magnitude less than $2g$.

\smallskip

\noindent The additional inertial mass $\d \, m = m - m_0$ is due to
the interaction of B with the surounding field of water and if
this interaction could not be turned off (which is the case
if we deal with an electron instead of a ping-pong ball) there would be no way to
measure the bare mass $m_0$.

\smallskip

\noindent The analogy between hydrodynamics and electromagnetism led
(through the work of Thomson, Lorentz, Kramers$\ldots$ \cite{MD}) to the
crucial distinction between the bare parameters, such as $m_0$, which
enter the field theoretic equations, and the observed parameters, such
as the inertial mass $m$.

\smallskip

\smallskip

\noindent A quantum field theory in $D=4$ dimensions, is given by
a classical action functional, 
\begin{equation}
S \, (A) =  \int \Lc \, (A) \,
d^4 x   \label{action}
\end{equation}
 where $A$ is a classical field and the
Lagrangian is of the form, 
\begin{equation}
 \Lc \, (A) =  (\part A)^2 / 2 -
\frac{m^2}{2} \, A^2 + \Lc_{\rm int} (A) \label{lag}
\end{equation}
 where
$\Lc_{\rm int} (A)$ is usually a polynomial in $A$ and possibly
its derivatives.

\smallskip

\noindent One way to describe the quantum fields $\phi (x)$, is by
means of the time ordered Green's functions,
\begin{equation}
G_N (x_1 , \ldots , x_N) = \lgl \, 0 \, \vert T \, \phi (x_1) \ldots
\phi (x_N) \vert \, 0 \, \rgl \label{16-3}
\end{equation}
where the time ordering symbol $T$ means that the $\phi (x_j)$'s are
written in order of increasing time from right to left.

\smallskip

\noindent The probability amplitude of a classical field configuration
$A$ is given by,
\begin{equation}
e^{i \, \frac{S(A)}{\hbar}} \label{16-4}
\end{equation}
and if one could ignore the renormalization problem, the Green's
functions would be computed as,
\begin{equation}
G_N (x_1 , \ldots , x_N) = \Nc \int e^{i \, \frac{S(A)}{\hbar}} \ A (x_1)
\ldots A (x_N) \, [dA] \label{16-5}
\end{equation}
where $\Nc$ is a normalization factor required to ensure the
normalization of the vacuum state,
\begin{equation}
\lgl \, 0 \mid 0 \, \rgl = 1 \, . \label{16-6}
\end{equation}

\smallskip

\noindent If one could ignore renormalization, the functional integral
\ref{16-5} would be easy to compute in perturbation theory, i.e. by
treating the term $\Lc_{\rm int}$ in \ref{lag} as a perturbation of
\begin{equation}
\Lc_0 (A) = (\part A)^2 / 2 - \frac{m^2}{2} \, A^2 \, .
\label{16-7}
\end{equation}
With obvious notations the action functional splits as
\begin{equation}
S (A) = S_0 (A) + S_{\rm int} (A) \label{16-8}
\end{equation}
where the free action $S_0$ generates a Gaussian measure $\exp \, (i
\, S_0 (A)) \, [d A] = d \mu$.

\smallskip

\noindent The series expansion of the Green's functions is then given
in terms of Gaussian integrals of polynomials as,
\begin{eqnarray}
G_N (x_1 , \ldots , x_N) = \left( \sum_{n=0}^{\ify} \, i^n / n! \int
A (x_1) \ldots A (x_N) \, (S_{\rm int} (A))^n \, d \mu
\right) \nonumber \\
\quad \left( \sum_{n=0}^{\ify} \, i^n / n! \int S_{\rm int} (A)^n \, d
\mu \right)^{-1} \, \nonumber 
\end{eqnarray}
The various terms of this expansion are computed using integration by
parts under the Gaussian measure $\mu$. This generates a large number
of terms $U(\G)$, each being labelled by a Feynman graph $\G$, and
having a numerical value $U(\G)$ obtained as a multiple integral in a
finite number of space-time variables. As a rule the unrenormalized values
$U(\G)$ are given by nonsensical divergent integrals.

\smallskip

\noindent The conceptually really nasty divergences are called
ultraviolet and are associated to
the presence of arbitrarily large frequencies or equivalently to
the unboundedness of momentum space on which integration has to be
carried out. Equivalently, when one attempts to integrate in
coordinate space, one confronts divergences along diagonals,
reflecting the fact that products of field operators are defined
only on the configuration space of distinct spacetime points.

\smallskip

\noindent The physics resolution of this problem is obtained by
first introducing a cut-off in momentum space (or any suitable
regularization procedure) and then by cleverly making use of the
unobservability of the bare parameters, such as the bare mass
$m_0$. By adjusting, term by term of the perturbative expansion,
the dependence of the bare parameters on the cut-off parameter, it
is possible for a large class of theories, called renormalizable,
to eliminate the unwanted ultraviolet divergences.

\smallskip

\noindent The main calculational complication of this subtraction procedure
occurs for diagrams which possess non-trivial subdivergences,
i.e.~subdiagrams which are already divergent. In that situation
the procedure becomes very involved since it is no longer a simple
subtraction, and this for two obvious reasons: i) the divergences
are no longer given by local terms, and ii) the previous
corrections (those for the subdivergences) have to be taken into
account.

\noindent To have an example for the combinatorial burden imposed by these
difficulties  consider the problem below of the
renormalization of a two-loop four-point function in massless
scalar $\phi^4$ theory in four dimensions, given by the following
Feynman graph:

$$
\hbox{
\psfig{figure=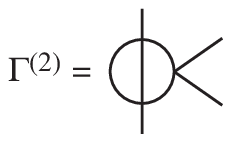}
}
$$

\noindent It contains a divergent subgraph:

$$
\hbox{
\psfig{figure=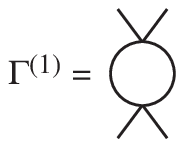}
}
$$

\noindent We work in the Euclidean framework and introduce a cut-off
$\lambda$ which we assume to be always much bigger than the square
of any external momentum $p_i$. Analytic expressions for these
Feynman graphs are obtained by utilizing a map $\Gamma_\lambda$
which assigns integrals to them according to the Feynman rules and
employs the cut-off $\lambda$ to the momentum integrations. Then
$\Gamma_\lambda[\Gamma^{[1,2]}]$ are  given by
\[
\Gamma_\lambda[\Gamma^{[1]}](p_i)=\int d^4k
\frac{\Theta(\lambda^2-k^2)}{k^2}\frac{1}{(k+p_1+p_2)^2},
\]
and \[ \Gamma_\lambda[\Gamma^{[2]}](p_i)=\int d^4l
\frac{\Theta(\lambda^2-l^2)}{l^2(l+p_1+p_2)^2}
\Gamma_\lambda[\Gamma^{[1]}(p_1,l,p_2,l)].
\]
It is easy to see that $\Gamma_\lambda[\Gamma^{[1]}]$
decomposes into the form $b
\log\lambda$ (where $b$ is a real number)
plus terms which remain finite for
 $\lambda\to\infty$, and hence will produce a
divergence which is a non-local function of external momenta
\[
\sim\log\lambda\int d^4l
\frac{\Theta(\lambda^2-l^2)}{l^2(l+p_1+p_2)^2}\sim\log\lambda\;\log(p_1+p_2)^2.
\]
Fortunately, the counterterm ${\cal L}_{\Gamma^{[1]}}$
$\sim\log\lambda$ generated to subtract the divergence in
$\Gamma_\lambda[\Gamma^{[1]}]$ will precisely cancel this
non-local divergence in $\Gamma^{[2]}$.

\noindent That this type of cancellation occurs at any order of perturbation theory,
i. e. that the two diseases above actually cure each other in general is a
very non-trivial fact that took decades to prove \cite{BPHZ}.

\noindent The detailed combinatorics is governed by the $\bar{R}$ operation
of Bogoliubov and Parasiuk (for a
1PI graph $\G$)
\begin{equation}
 \bar{R}(\G)= U (\G) + \sum_{\g \sbsneq \G} C (\g) \, U (\G / \g) \label{BP1}
\end{equation}
which prepares a given graph with unrenormalized value $U (\G)$ by adding the 
counterterms $ C (\g)$. The latter are constructed by induction using
\begin{equation}
C \, (\G) = -T \left( U (\G) + \sum_{\g \sbsneq \G} C (\g) \, U (\G / \g)
\right) \label{BP2}
\end{equation}
where, using dimensional regularization $T$ is just the extraction of the 
pole part in $D=4- \epsilon$.
The renormalized graph is then given by
\begin{equation}
R \, (\G) = U (\G) + C (\G) + \sum_{\g \sbsneq \G} C (\g) \, U (\G / \g) \,.
\label{BP3}
\end{equation}
 For
a mathematician the intricacies of the detailed combinatorics and the lack of any
obvious mathematical structure underlying it make it totally
inaccessible, in spite of the existence of a satisfactory formal
approach to the problem \cite{EG}.

\noindent This situation was drastically changed by the discovery by Dirk Kreimer (\cite{dhopf})
who understood that the formula for the $\bar{R}$ operation
in fact dictates a Hopf algebra coproduct on the free commutative
algebra ${\cal H}_K$ generated by the 1PI graphs $\G$,
\begin{equation}
\D \, \G= \G \ot 1 + 1 \ot \G+ \sum_{\g \sbsneq \G} \, \g \ot \G / \g \ , 
\label{36}
\end{equation}
(In fact he first formulated it in terms of rooted trees,
but the graph formulation is easier to explain).

\noindent This Hopf algebra is commutative as an algebra and we showed in \cite{CK1}, \cite{CK} that it 
is the dual Hopf algebra of the
envelopping algebra of a Lie algebra $\underline G$ whose basis is labelled by
the one particle irreducible Feynman graphs. The Lie bracket of two such
graphs is computed from insertions of one graph in the other and vice
versa. The corresponding Lie group $G$ is the group of characters of
$\Hc$.

\smallskip
\noindent The next breakthrough came from our joint discovery  \cite{CK} that identical
formulas to equations (12, 13, 14) occur in the solution of the Riemann Hilbert problem for an arbitrary pronilpotent Lie group $G$!

\noindent This really unveils the true nature of this seemingly complicated combinatorics and shows
 that it is a special case of a general extraction of
finite values based on the Riemann-Hilbert problem.

\noindent  The Riemann-Hilbert problem comes
from Hilbert's $21^{\rm st}$ problem which he formulated as
follows:

\begin{itemize}
\item[] ``Prove that there always exists a Fuchsian linear
differential equation with given singularities and given monodromy.''
\end{itemize}

\noindent In this form it admits a positive answer due to Plemelj
and Birkhoff (cf.~\cite{beauv} for a careful exposition). When
formulated in terms of linear systems of the form, 
\begin{equation}
y'(z) = A(z)
\, y(z) \ , \ A(z) = \sum_{\a \in S} \ \frac{A_{\a} }{z - \a} \,
, \label{13} 
\end{equation}
(where $S$ is the given finite set of
singularities, $\ify \not\in S$, the $A_{\a}$ are complex matrices
such that 
\begin{equation}
\sum \ A_{\a} = 0 \label{14} 
\end{equation}
to avoid
singularities at $\ify$), the answer is not always positive
\cite{boli}, but the solution exists when the monodromy matrices
$M_{\a}$ are sufficiently close to 1. It can then
be explicitly written as a series of polylogarithms \cite{beauv}.

\noindent Another formulation of the Riemann-Hilbert problem,
intimately tied up to the classification of holomorphic vector
bundles on the Riemann sphere $P_1 (\Cb)$, is in terms of the
Birkhoff decomposition 
\begin{equation}
\g \, (z) = \g_- (z)^{-1} \, \g_+ (z)
\qquad z \in C \label{15} 
\end{equation}
where we let $C \sbs P_1 (\Cb)$ be a
smooth simple curve, $C_-$ the component of the complement of $C$
containing $\ify \not\in C$ and $C_+$ the other component. Both
$\g$ and $\g_{\pm}$ are loops with values in ${\rm GL}_n (\Cb)$,
\begin{equation}
\g \, (z) \in G = {\rm GL}_n (\Cb) \qquad \fl \, z \in \Cb
\label{16} 
\end{equation}
and $\g_{\pm}$ are boundary values of holomorphic
maps (still denoted by the same symbol)
\begin{equation}
\g_{\pm} : C_{\pm} \ra {\rm GL}_n (\Cb) \, .
\label{17}
\end{equation}
The normalization condition $\g_- (\ify) = 1$
ensures that, if it exists, the decomposition (18) is unique
(under suitable regularity conditions).

\smallskip

\noindent The existence of the Birkhoff decomposition (18) is
equivalent to the vanishing,
\begin{equation}
c_1 \, (L_j) = 0 \label{18}
\end{equation}
of the Chern numbers $n_j = c_1 \, (L_j)$ of the holomorphic line
bundles of the Birkhoff-Grothendieck decomposition,
\begin{equation}
E = \op \, L_j \label{19}
\end{equation}
where $E$ is the holomorphic vector bundle on $P_1 (\Cb)$ associated to
$\g$, i.e. with total space:
\begin{equation}
(C_+ \ts \Cb^n)\cup_{\g} (C_- \ts \Cb^n) \, . \label{20}
\end{equation}
The above discussion for $G = {\rm GL}_n (\Cb)$ extends to arbitrary
complex Lie groups.

\smallskip

\noindent When $G$ is a simply connected nilpotent complex Lie group the
existence (and uniqueness) of the Birkhoff decomposition (18)
is valid for any $\g$. When the loop $\g : C \ra G$ extends to a
holomorphic loop: $C_+ \ra G$, the Birkhoff decomposition is given by
$\g_+ = \g$, $\g_- = 1$. In general, for $z_0 \in C_+$ the evaluation,
\begin{equation}
\g \ra \g_+ (z_0) \in G \label{21}
\end{equation}
is a natural principle to extract a finite value from the singular
expression $\g (z_0)$. This extraction of finite values is a multiplicative 
removal of the pole part for a meromorphic loop $\g$ when we let
 $C$ be an infinitesimal circle centered at $z_0$.

\noindent We are now ready to apply this procedure in Quantum Field Theory.
First, using dimensional regularization, the
bare (unrenormalized) theory gives rise to a meromorphic loop,
\begin{equation}
\g (z) \in G \ , \qquad z \in \Cb
\end{equation}
 Our main result \cite{ck, CK}  is that the renormalized theory is just
the evaluation at the integer
dimension $z_0 = D$ of space-time of the holomorphic part $\g_+$ of the
Birkhoff decomposition of $\g$.

\smallskip

\noindent In fact, the original loop $d \ra \g (d)$ not only depends upon the
parameters of the theory but also on the additional ``unit of
mass'' $\mu$ required by dimensional analysis. We showed in \cite {cknew}
 that the mathematical concepts developped in our earlier papers
provide very powerful tools to lift the usual concepts of the
$\b$-function and renormalization group from the space of coupling
constants of the theory to the complex Lie group $G$.

\noindent  We first observed that even though the loop $\g (d)$ does
depend on the additional parameter $\mu$, 
\begin{equation}
\mu \ra \g_{\mu} (d) \, , \label{1} 
\end{equation}
the negative part $\g_{\mu^-}$ in the Birkhoff
decomposition, 
\begin{equation}
\g_{\mu} (d) = \g_{\mu^-} (d)^{-1} \, \g_{\mu^+}
(d) \label{2} 
\end{equation}
is actually independent of $\mu$, 
\begin{equation}
\frac{\partial}{\partial \mu} \, \g_{\mu^-} (d) = 0 \, . \label{3} 
\end{equation}
This is a restatement of a well known fact and follows
immediately from dimensional analysis. Moreover, by construction,
the Lie group $G$ turns out to be graded, with grading, 
\begin{equation}
\t_t \in {\rm Aut} \, G \ , \quad t \in \Rb \, , \label{4} 
\end{equation}
inherited from the grading of the Hopf algebra $\Hc$ of Feynman
graphs given by the loop number, 
\begin{equation}
L (\G) = \hbox{loop number of}
\ \G \label{5} 
\end{equation}
for any 1PI graph $\G$.

\noindent  The straightforward equality, 
\begin{equation}
\g_{e^t \mu} (d) = \t_{t \ve}
(\g_{\mu} (d)) \qquad \fl \, t \in \Rb \, , \ \ve = D-d \label{6}
\end{equation}
shows that the loops $\g_{\mu}$ associated to the
unrenormalized theory satisfy the striking property that the
negative part of their Birkhoff decomposition is unaltered by the
operation, 
\begin{equation}
\g (\ve) \ra \t_{t\ve} (\g (\ve)) \, , \label{7} 
\end{equation}
In other words, if we replace $\g (\ve) $ by $\t_{t\ve} (\g
(\ve))$ we dont change the negative part of its Birkhoff
decomposition. We settled now for the variable, 
\begin{equation}
\ve = D - d \in
\Cb \backslash \{ 0 \} \, . \label{8} 
\end{equation}
We give in \cite {cknew}
 a complete characterization of the loops $\g (\ve)
\in G$ fulfilling the above striking invariance. This
characterization only involves the negative part $\g_- (\ve)$ of
their Birkhoff decomposition which by hypothesis fulfills, 
\begin{equation}
\g_-
(\ve) \, \t_{t \ve} (\g_- (\ve)^{-1}) \ \hbox{is convergent for} \
\ve \ra 0 \, . \label{9} 
\end{equation}
It is easy to see that the limit of
(34) for $\ve \ra 0$ defines a one parameter subgroup, 
\begin{equation}
F_t \in G
\, , \ t \in \Rb \label{10} 
\end{equation}
and that the generator $\b =
\left( \frac{\partial}{\partial t} \, F_t \right)_{t=0}$ of this
one parameter group is related to the {\it residue} of $\g$ 
\begin{equation}
\build{\rm Res}_{\ve = 0}^{} \g = - \left( \frac{\partial}{\partial
u} \, \g_- \left( \frac{1}{u} \right) \right)_{u=0} \label{11} 
\end{equation}
by the simple equation, 
\begin{equation}
\b = Y \, {\rm Res} \, \g \, , \label{12} 
\end{equation}
where $Y = \left( \frac{\partial}{\partial t} \, \t_t
\right)_{t=0}$ is the grading.

\noindent  This is straightforward but our result is the following formula
(39) which gives $\g_- (\ve)$ in closed form as a function of
$\b$. We shall for convenience introduce an additional generator
in the Lie algebra of $G$ (i.e.~primitive elements of $\Hc^*$)
such that, 
\begin{equation}
[Z_0 , X] = Y(X) \qquad \fl \, X \in \hbox{Lie} \ G
\, . \label{13} 
\end{equation} 
The scattering formula for $\g_- (\ve)$ is
then, 
\begin{equation}
\g_- (\ve) = \lim_{t \ra \ify} e^{-t \left(
\frac{\b}{\ve} + Z_0 \right)} \, e^{t Z_0} \, . \label{14} 
\end{equation} 
Both factors in the right hand side belong to the semi direct
product, 
\begin{equation}
\wt G = G \, \semi_{\t} \, \Rb \label{15} 
\end{equation} 
of the
group $G$ by the grading, but of course the ratio (39) belongs to
the group $G$.

\noindent  This shows (\cite {cknew}) that the higher pole structure of the
divergences is uniquely determined by the residue and gives a
strong form of the t'Hooft relations, which will come as an
immediate corollary.

\noindent  The main new result of \cite {cknew}, specializing to the massless case
and taking $\vp_6^3$ as an illustrative example to
fix ideas and notations, is that the
formula for the bare coupling constant, 
\begin{equation}
g_0 = g \, Z_1 \,
Z_3^{-3/2} \label{16} 
\end{equation} 
where both $g \, Z_1 =g + \d g$ and the
field strength renormalization constant $Z_3$ are thought of as
power series (in $g$) of elements of the Hopf algebra $\Hc$, does
define a Hopf algebra homomorphism, 
\begin{equation}
\Hc_{CM} \build
\longra_{}^{g_0} \Hc_K \, , \label{17} 
\end{equation} 
from the Hopf algebra
$\Hc_{CM}$ of coordinates on the group of formal diffeomorphisms
of $\Cb$ such that, 
\begin{equation}
\vp (0) = 0 \, , \ \vp' (0) = {\rm id}
\label{18} 
\end{equation} 
to the Hopf algebra $\Hc_K$ of the massless theory.
We had already constructed in \cite{CK} a Hopf algebra
homomorphism from $\Hc_{\rm CM}$ to the Hopf algebra of rooted
trees, but the physical significance of this construction was
unclear.

\noindent  The homomorphism (42) is quite different in that for instance the
transposed group homomorphism, 
\begin{equation}
G \build \longra_{}^{\rho} {\rm
Diff} (\Cb) \label{19} 
\end{equation} 
lands in the subgroup of {\it odd}
diffeomorphisms, 
\begin{equation}
\vp (-z) = -\vp (z) \qquad \fl \, z \, .
\label{20}
\end{equation} 
Moreover its physical significance is 
transparent. In particular the image by $\rho$
of $\b = Y \, {\rm Res} \, \g$ is the usual $\b$-function of the
coupling constant $g$.

\noindent  We discovered the homomorphism (42) by lengthy concrete
computations which were an excellent test for the explicit ways of handling the coproduct,
coassociativity, symmetry factors$\ldots$ that underly the theory.

\noindent  As a corollary of the construction of $\rho$ we get an {\it
action} by (formal) diffeomorphisms of the group $G$ on the space
$X$ of (dimensionless) coupling constants of the theory. We can
then in particular formulate the Birkhoff decomposition {\it
directly} in the group, 
\begin{equation}
{\rm Diff} \, (X) \label{21} 
\end{equation} 
of
formal diffeomorphisms of the space of coupling constants.

\medskip
\noindent {\bf Theorem 10.} (\cite{cknew})
 {\it Let the unrenormalized effective 
coupling constant $g_{\rm eff} (\ve)$ be viewed as a formal power series 
in $g$ and let $g_{\rm eff}(\ve) =  g_{{\rm eff}_+} (\ve)\, (g_{{\rm eff}_-}(\ve))^{-1}$ 
be its (opposite) Birkhoff decomposition in the group of formal 
diffeomorphisms. Then the loop $g_{{\rm eff}_-} (\ve)$ is the bare 
coupling constant and $g_{{\rm eff}_+} (0)$ is the renormalized 
effective coupling.}
\medskip

 \noindent  This result is now, in its statement, no longer depending upon
our group $G$ or the Hopf algebra $\Hc$. But of course the proof
makes heavy use of the above ingredients.
\noindent  It is a  challenge to physicists to find a direct proof of this result.

\noindent  Finally the Birkhoff decomposition of a loop, 
\begin{equation}
\d (\ve) \in {\rm
Diff} \, (X) \, , \label{26} 
\end{equation} 
admits a beautiful geometric
interpretation. If we let $X$ be a complex manifold and pass from
formal diffeomorphisms to actual ones, the data (47) is the initial
data to perform, by the clutching operation, the construction of a
complex bundle, 
\begin{equation}
P = (S^+ \ts X) \, \cup_{\d} (S^- \ts X) \label{27} 
\end{equation} 
over the sphere $S = P_1 (\Cb) = S^+ \cup S^-$, and with
fiber $X$, 
\begin{equation}
X \longra P \build \longra_{}^{\pi} S \, . \label{28}
\end{equation}
The meaning of the Birkhoff decomposition, 
\begin{equation}
\d (\ve) = \d_- (\ve)^{-1} \, \d_+ (\ve) 
\end{equation}
is then exactly
captured by an isomorphism of the bundle $P$ with the trivial
bundle, 
\begin{equation}
S \ts X \, . \label{29} 
\end{equation}

\bigskip

    \section{Number theory}  

\noindent I shall conclude these notes by giving a brief glimpse at the connection  
between noncommutative geometry and number theory. 
 There are two points of 
contact of the two subjects, the first gives a spectral interpretation 
of zeros of zeta and $L$-functions in terms of a construction 
involving adeles, more specifically the noncommutative space of adele classes.
The second has to do with the missing Galois theory 
at Archimedean places. For 
the specialists of quantum chaos looking for a spectral realization of the non-trivial 
zeros of the Riemann zeta function from the quantization of 
classical mechanical systems, the adeles might look rather exotic 
at first sight and we first need to explain briefly (for non specialists)
why Ideles and Adeles are natural and important in number theory.

\smallskip

\noindent Let us start with the reciprocity law (Gauss 1801)
\begin{equation}
\left( \frac{\ell}{p} \right) = (-1)^{\ve (p)\ve (\ell)} \, \left( 
\frac{p}{\ell} \right) \ , \qquad \ve (p) = \frac{p-1}{2} \ ({\rm 
mod} \, 2)
\end{equation}
where $\ell$ and $p$ are odd primes and $\left( \frac{\ell}{p} 
\right)$ is the Legendre symbol whose value is $+1$ if the equation
\begin{equation}
x^2 = \ell \qquad ({\rm mod} \, p)
\end{equation}
admits a solution, and is $-1$ if it does not.

\smallskip

\noindent For instance, with $\ell = 5$ we see that whether the equation $x^2 
= 5 ({\rm mod} \, p)$ admits a solution only depends upon $p({\rm mod} \, 5)$, i.e. only on 
the last digit of $p$. Thus the answer is the same for $p=7$ and 
$p=1997$ or for $p=19$ and $p=1999$. It follows that the primes $p$ 
thus fall into {\it classes}. The language of Adeles and Ideles extends this simple notion of {\it classes} of 
primes to those of {\it ideal classes} and then of {\it Idele 
classes}.

\smallskip

\noindent To the proof of Dirichlet of the existence of infinitely many 
primes in an arithmetic progression corresponds the construction of 
an $L$-function associated to a character $\chi$ modulo $m$,

\begin{equation}
L(s,\chi) = \prod \frac{1}{1 - \chi (p) \, p^{-s}} \, .
\end{equation}
More generally a Hecke $L$-function is associated to a character of 
the ideal class group modulo $m$ and in fact also to a 
Gr\"ossencharakter which is a character of the Idele class group of 
a number field $k$.

\smallskip

\noindent The quickest way to introduce the Idele class group of a number field 
$k$ is to understand  (Cf. 
Iwasawa {\it Ann. of Math.} {\bf 57} (1953)) that such a field sits uniquely as a discrete 
cocompact subfield of a unique locally compact ( semi-simple and non discrete)
 ring $A$
\begin{equation}
k \subset A \ , \qquad k \, \, \, \hbox{cocompact}
\end{equation}
called the ring of {\it Adeles} of $k$. One then has,
\begin{equation}
\hbox{Idele class group of} \ k = {\rm GL}_1 (k) \backslash {\rm 
GL}_1 (A) \, ,
\end{equation}
and a Gr\"ossencharakter is a character of this locally compact 
group. Iwasawa and Tate showed how to use analysis on adeles to 
prove the basic properties of the Hecke $L$-functions which were 
then extended to $L$-functions associated to automorphic forms 
which appear in the action of ${\rm GL}_n (A)$ on the Hilbert space
\begin{equation}
L^2 ({\rm GL}_n (k) \backslash {\rm GL}_n (A)) \, .
\end{equation}
To understand the other language involved in the basic dictionary 
which underlies the Langlands program let us go back to the 
equation (2) say with $\ell = 5$ and simply adjoin $\sqrt{5}$ to the 
field $\Qb$ of rational numbers which gives an algebraic extension
$K= \Qb 
(\sqrt{5})$ of $k=\Qb$. The Galois group ${\rm Gal} (\Qb 
(\sqrt{5}) : \Qb) = {\rm Gal} (K:k)$ is of course $\Zb / 2$ in this 
case and admits an obvious non trivial one dimensional 
representation $\pi$. In general, the {\it Artin} $L$-function 
associated to a representation
\begin{equation}
{\rm Gal} (K:k) \ra {\rm GL} (n,\Cb)
\end{equation}
of the Galois group of a finite Galois extension $K$ of $k$, is
\begin{equation}
L(s,\pi) = \prod_p L_p (s,\pi)
\end{equation}
where $p$ runs through the prime ideals of $k$ and the local $L$ 
factor is given at unramified $p$ by,
\begin{equation}
L_p (s,\pi) = \det (1 - \pi (\s) \, N(p)^{-s})
\end{equation}
where $\s$ is the Frobenius automorphism of $p$.

\smallskip

\noindent When $K/k$ is an abelian extension and $\pi$ a one dimensional 
representation it follows from class field theory that $\pi$ defines 
a character modulo the conductor of $K/k$ and that the Artin 
$L$-function equals the Hecke $L$-function. This Artin reciprocity 
law is a far reaching extension of the Gauss reciprocity law (1).

\smallskip

\noindent The Langlands program extends Hecke's theory of Euler products 
associated to automorphic forms on ${\rm GL} (2)$ to arbitrary 
reductive groups $G$ and gives a correspondance, extending Artin's 
reciprocity law to the non-Abelian case, between automorphic 
representations of $G$ and representations,
\begin{equation}
{\rm Gal} (K:k) \ra \ ^L{G}
\end{equation}
in the Langlands dual $^L{G}$ of $G$.

\smallskip

\noindent A basic tool of the theory is the trace formula \cite{Arthur} which 
extends to the adelic context the Selberg trace formula. The trace 
formula is the equality obtained by computing in two different ways 
the trace of operators of the form, 
\begin{equation}
{\rm Trace} \, (C_{\Lb} \, \pi (f))
\end{equation}
where (for $G = {\rm GL}_n$), $\pi$ is the natural representation of ${\rm GL}_n (A)$ in the 
Hilbert space 
\begin{equation}
L_{\chi}^2 ({\rm GL}_n (k) \backslash {\rm GL}_n (A))
\end{equation}
where $\chi$ is a Gr\"ossencharakter, and $C_{\Lb}$ is a cutoff.
The Gr\"ossencharakter $\chi$ allows one to restrict to vectors with a 
fixed behaviour relative to ${\rm GL}_1$.

\smallskip

\noindent The spectral side of the trace formula is obtained from the harmonic 
analysis of the representation $\pi$. The geometric side expresses the 
trace as a sum of orbital integrals.

\smallskip

\noindent The restriction imposed in (12) by the Gr\"ossencharakter $\chi$ shows 
that the case $n=1$ becomes trivial and concentrates essentially on 
the ${\rm SL}_n$ aspect for $n \geq 2$. So far the zeros of 
$L$-functions do not appear in this language.

\smallskip

\noindent Our contribution to this subject is to show that both the zeros of 
$L$-functions and the Riemann-Weil explicit formulas appear directly 
in a refinement of the trace formula obtained as follows. Instead of 
restricting the Hilbert space,
\begin{equation}
L^2 ({\rm GL}_n (k) \backslash {\rm GL}_n (A))
\end{equation}
by the choice of Gr\"ossencharakter $\chi$ as above, one introduces on 
the full Hilbert space (13) a finer cutoff operator $Q_{\Lb}$ taking 
care of the ``${\rm GL}_1$'' behaviour of vectors. 

\smallskip

\noindent To understand in which way the corresponding trace formula refines the 
Arthur trace formula, it is simplest to restrict to the case of ${\rm 
GL}_1$. In order to simplify even further we shall replace the number 
field $k$ by a function field of positive characteristic. This allows 
for a straightforward definition of the cutoff operators $Q_{\Lb}$ as 
the orthogonal projection on the subspace,
\begin{equation}
Q_{\Lb} \subset L^2 ({\rm GL}_1 (k) \backslash {\rm GL}_1 (A))
\end{equation}
spanned by the vectors $\xi \in {\cal 
S} (A)$ (averaged on ${\rm GL}_1 (k)$) which vanish as well as their Fourier transform for $\vert x 
\vert > \Lb$. Note that we use Fourier transform on the {\it additive} 
group of adeles so that the space ${\rm GL}_1 (k) \backslash A$ of 
Adele classes is implicit in this construction. To define this Fourier transformation
 we needed to choose a  basic character $\a =\prod \a_v$ of the additive group $A$ 
for which the lattice $k$ is selfdual. 

\smallskip

\noindent The spectral computation of the trace of $Q_{\Lb} \, \pi (f)$ involves 
all the nontrivial zeros of Hecke $L$-functions and is given by the 
following formula (\cite{CoRH}),
\begin{eqnarray}
{\rm Trace} (Q_{\Lb} \, \pi (f)) = 2 \left( \sum_{{\rm GL}_1} f(k) 
\right) \log' \Lb \nonumber \\
+ \, \wh f (0) + \wh f (1) - \sum_{{L \left( \chi , \frac{1}{2} + \rho 
\right) = 0 \atop \rho \in B/N^{\perp}}} N \left( \chi , \frac{1}{2} + 
\rho \right) \int_{i \Rb} \wh f (\chi , z) \, d\mu_{\rho} (z) + o(1)
\end{eqnarray}
where $B$ is the open strip $B = \left\{ \rho \in \Cb \, ; \vert {\rm 
Re} \, \rho \vert < \frac{1}{2} \right\}$, $N \left( \chi , 
\frac{1}{2} + \rho \right)$ is the multiplicity of $\frac{1}{2} + 
\rho$ as a zero of the $L$ function $L(\chi , s)$, $\chi$ varying 
through Gr\"ossencha\-rakters (modulo principal ones), $N$ being the 
module,
\begin{equation}
N = {\rm Mod} (k) \, .
\end{equation}
Also $2 \log' \L = \int_{ \vert \lb \vert \in
[\L^{-1}, \L]} d^* \lb$, and the measure $d\mu_{\rho} (z)$ is the harmonic measure of $\rho \in \Cb$ with 
respect to the line $i\Rb$. In particular if the zero $\frac{1}{2} + 
\rho$ is on the critical line $d\mu_{\rho} (z)$ is just the Dirac mass 
at $z = \rho$. Finally the Fourier transform of $f$ is given by,
\begin{equation}
\wh f (\chi , z) = \int_{{\rm GL}_1 (A)} f(u^{-1}) \, \chi (u) \, \vert u 
\vert^z \, d^* u \, .
\end{equation}
The geometric side of the trace formula has so far only be fully 
justified in the simplified situation where only finitely many places 
are used. It is then given by the following formula (\cite{CoRH})
\begin{equation}
{\rm Trace} (Q_{\Lb} \, \pi (f)) =  2 \left( \sum_{{\rm GL}_1} f(k) 
\right) \log' \Lb  + \sum_{v,k } \int'_{k^*_v} 
{ f(ku) \over \vert 1-u \vert} \, d^* u + o(1) \, ;
\end{equation}
where each $k^*_v$ is embedded in $({\rm GL}_1 (k) \backslash {\rm GL}_1 (A))$ 
by the map $u \rightarrow (1,1,...,u,...,1)$ and the principal
value $\int'$ is uniquely determined by the pairing with the 
unique distribution on $k_v$ which agrees with ${du \over \vert 1-u 
\vert}$ for $u \not= 1$ and whose Fourier transform relative to $\a_v$ vanishes 
at $1$.

\noindent By proving that it entails the
positivity of the Weil distribution, we showed in \cite{CoRH}  that the validity of the geometric
 side, i.e., the global trace formula, is equivalent to the Riemann 
Hypothesis for all $L$-functions with Gr\"ossencharakter. 

\medskip
\noindent {\bf Theorem 11.} {\it The following two conditions are equivalent: 

\smallskip

a) When $\L \ra \ify$, one has, for all $f \in \Sc ({\rm GL}_1 (k) \backslash {\rm GL}_1 (A))$ with
compact support,
\begin{equation}
{\rm Trace} (Q_{\Lb} \, \pi (f)) =  2 \left( \sum_{{\rm GL}_1} f(k) 
\right) \log' \Lb  + \sum_{v,k } \int'_{k^*_v} 
{ f(ku) \over \vert 1-u \vert} \, d^* u + o(1) \, ;
\end{equation}  }
\smallskip

{\it b) All $L$ functions with Gr\"ossencharakter on $k$ satisfy the Riemann 
Hypothesis.}

\medskip

\noindent We have thus obtained a spectral interpretation of the zeros of zeta and L-functions
 as 
an absorption spectrum, i.e., as missing spectral lines. All zeros do play a 
role in the spectral side of the trace formula, but while the critical zeros 
do appear perse, the noncritical ones appear as resonances and 
enter in the trace formula through their harmonic potential with respect to the 
critical line. The spectral side is entirely canonical, and its validity is justified
in the global case \cite{CoRH}. It is quite important to understand 
why a crucial negative sign in the analysis of the statistical fluctuations of 
the zeros of zeta indicated from the start that the spectral interpretation should be as an 
absorption spectrum, or equivalently should be of a cohomological nature.

$$
\hbox{
\psfig{figure=Fig.10.eps}
}
$$

\noindent The number of
zeros of zeta whose imaginary part is less than $E > 0$,
\begin{equation}
N(E) = \# \ \hbox{of zeros} \ \rho \ , \ 0 < {\rm Im} \,
\rho < E \label{1}
\end{equation}
has an asymptotic expression (\cite{[R]}) given by
\begin{equation}
N(E) = {E \over 2\pi} \ \left( \log \left( {E \over 2\pi}
\right) -1 \right) + {7 \over 8} + o(1) + N_{\rm osc} (E)
\label{2}
\end{equation}
where the oscillatory part of this step function is
\begin{equation}
N_{\rm osc} (E) = {1\over \pi} \ {\rm Im} \, \log \, \z \,
\left( {1\over 2} + iE \right) \label{3}
\end{equation}

\noindent  which is of the order of $ \hbox{Log}(E)$ (We assume that $E$ is not the imaginary part of a zero and
take for the logarithm the branch which is $0$ at $+\ify$).
The Euler
product formula for the zeta function yields (cf. \cite{Berry}) a
heuristic asymptotic formula for $N_{\rm osc} (E)$,
\begin{equation}
N_{\rm osc} (E) \sm {-1 \over \pi} \sum_p \sum_{m=1}^{\ify}
{1\over m} \, {1 \over p^{m/2}} \, \sin \, (m \, E \, \log \,
p) \, . \label{7}
\end{equation}

\medskip

\noindent One can compare this formula with what appears in the direct atempt  \cite{Berry} to construct 
a spectral realization of zeros of zeta from quantization of a classical dynamical system. 
 In this theory the quantization of the classical dynamical system 
given by the phase space $X$ and hamiltonian $h$ gives rise to a Hilbert
space $\Hc$ and a selfadjoint operator $H$ whose spectrum is the essential 
physical observable of the system. For complicated systems the only useful 
information about this spectrum is that, while the average part of the counting
function, 
\begin{equation}
N(E) = \# \ \hbox{eigenvalues of $H$ in} \ [0,E] \label{17.1}
\end{equation}
is computed by a semiclassical approximation mainly as a
volume in phase space, the oscillatory part,
\begin{equation}
N_{\rm osc} (E) = N(E) - \lgl N(E) \rgl \label{17.2}
\end{equation}
is the same as for a random matrix, governed by the statistic
dictated by the symmetries of the system.

\smallskip
\noindent One can then (\cite{Berry}) write down an asymptotic semiclassical approximation to the 
oscillatory function $N_{\rm osc} (E)$
\begin{equation}
N_{\rm osc} (E) = {1\over \pi} \ {\rm Im} \int_0^{\ify} {\rm
Trace} (H-(E+i\eta))^{-1} \, id\eta \label{17.5}
\end{equation}
using the stationary phase approximation of the
corresponding functional integral. For a system whose
configuration space is 2-dimensional, this gives (\cite{Berry} (15)),
\begin{equation}
N_{\rm osc} (E) \sm {1\over \pi} \sum_{\g_p} \sum_{m=1}^{\ify}
{1\over m} \, {1\over 2{\rm sh} \left( {m\lb_p \over
2}\right)} \, \sin (S_{\rm pm} (E)) \label{17.6}
\end{equation}
where the $\g_p$ are the primitive periodic orbits, the label
$m$ corresponds to the number of traversals of this orbit,
while the corresponding instability exponents are $\pm
\lb_p$. The phase $S_{\rm pm} (E)$ is up to a constant equal
to $m \, E \, T_{\g}^{\#}$ where $T_{\g}^{\#}$ is the period
of the primitive orbit $\g_p$.

\smallskip

\noindent Comparing the formulas  one sees a  fundamental mismatch  (cf.\cite{Berry}) which is
the overall {\it minus sign} in front of formula (23) as opposed to the plus sign of (27). This problem is resolved in our
spectral interpretation by the minus sign present in the spectral side of the trace formula (15).
The point is that the spectral analysis of the action of the Idele class group on the Adele class space
shows (\cite{CoRH}) white light with dark absorption lines labelled by the zeros of zeta and L-functions.
This also provides the correct explanation for the asymptotic form of
 the formula for the average number of zeros
\begin{equation}
\lgl N(E) \rgl \sim (E/2\pi) (\log (E/2\pi) -1) + 7/8 + o(1) 
\end{equation}
from a semiclassical computation for the number of quantum 
mechanical states in one degree of freedom which fulfill the conditions
\begin{equation}
\vert q \vert \leq \Lambda , \vert p \vert \leq \Lambda , \vert H \vert 
\leq E \, , 
\end{equation}
where $H = 2\pi q p$ is the Hamiltonian which generates the group
 involved in the action of the Idele class group namely the scaling 
transformations (see (\cite{CoRH}) for precise normalization).
 We are thus computing the area of,
\begin{equation}
  D = \{ (p,q) ; pq \geq 0 , \vert q \vert \leq \Lambda , \vert p \vert \leq \Lambda , 
\vert pq \vert \leq E/2\pi \} . 
\end{equation}
(since we deal with zeta alone we restrict ourselves to 
even functions so that we exclude the region $pq \leq 0$ of the semiclassical 
$(p,q)$ plane).
The computation yields
\begin{equation}
1/2 \int_{D}dp dq = {2E \over 2\pi} \, \log \Lambda - {E \over 2\pi} 
\left( \log \, {E \over 2\pi} - 1 \right) \, .
\end{equation}
In this formula we see in the right hand side the overall term $\langle N(E) \rangle$  which
appears with a {\it minus} sign. This shows that the number of quantum 
mechanical states is equal to ${4E \over 2\pi} \, 
\log \Lambda$ minus the first approximation to the number of zeros of zeta 
whose imaginary part is less than $E$ in absolute value (one just multiplies by 
$2$ the equality (31)). Now ${1 \over 
2\pi} \, (2E) (2\log \Lambda)$ is the number of quantum states in the 
Hilbert space $L^2 (\Rb_+^* , d^* x)$ which are localized in $\Rb_+^*$ 
between $\Lambda^{-1}$ and $\Lambda$ and are localized in the dual group 
$\Rb$ (for the pairing $\langle \lambda , t \rangle = \lambda^{it}$) 
between $-E$ and $E$. Thus we see clearly that the first approximation to 
$N(E)$ appears as the lack of surjectivity of the map which associates to 
quantum states $\xi$ with support in $D$ the function on $\Rb_+^*$,
\begin{equation}
E(\xi) (x) = \vert x \vert^{1/2} \sum_{n\in \Zb} \xi (nx) \label{49}
\end{equation}
where we assume the additional conditions $\xi (0) = \int \xi (x) dx = 0$.

\smallskip

\noindent A finer analysis, which is just what the trace formula is doing,
would yield the additional terms $ 7/8 + o(1) + N_{osc} (E)$.
The above discussion yields an explicit construction of a large
matrix whose spectrum approaches the zeros of zeta as $ \L \ra \ify$.

\noindent While the above discussion clearly gives the sought for spectral interpretation of zeros
of zeta it is unclear that one can expect to justify the (geometric side of) trace formula without a deeper
understanding of the symmetries of the situation, which might well involve quantum groups.

\noindent As we mentionned earlier, the second point of contact between noncommutative geometry
and number theory has to do with the missing Galois theory at Archimedean places.

\noindent   Let $k$ be a {\it global} field, when the characteristic of $k$ is $p>1$ so that $k$ is a function field 
over $\Fb_q$, one has
$$
k \subset k_{\rm un} \subset k_{\rm ab} \subset k_{\rm sep} \subset 
\overline k \, , 
$$
where $\overline k$ is an algebraic closure of $k$, $k_{\rm 
sep}$ the separable algebraic closure, $k_{\rm ab}$ the maximal abelian 
extension and $k_{\rm un}$ is obtained by adjoining to $k$ all roots of 
unity of order prime to $p$.

\noindent One defines the Weil group $W_k$ as the subgroup of ${\rm Gal} 
(k_{\rm ab} : k)$ of those automorphisms which induce on $k_{\rm un}$ an  
integral power of the Frobenius automorphism $\sigma$,
$$
\sigma (\mu) = \mu^q \qquad \forall \, \mu \ \hbox{root of 1 of order 
prime to} \ p \, . 
$$
The main theorem of global class field theory asserts the existence of a  
canonical isomorphism,
$$
W_k \simeq C_k = GL_1 (A) / GL_1 (k) \, ,
$$
of locally compact groups.

\noindent When $k$ is of characteristic 0, i.e. is a number field, one has a 
canonical isomorphism,
$$
{\rm Gal} (k_{\rm ab} : k) \simeq C_k / D_k \, , 
$$
where $D_k$ is the connected component of identity in the Idele class 
group $C_k $, but because of the Archimedian places 
of $k$ there is no interpretation of $C_k$ analogous to the Galois group 
interpretation for function fields. According to A.~Weil \cite{W4}, ``La 
recherche d'une interpr\'{e}tation pour $C_k$ si $k$ est un corps de nombres, 
analogue en quelque mani\`{e}re \`{a} l'interpr\'{e}tation par un groupe de Galois 
quand $k$ est un corps de fonctions, me semble constituer l'un des 
probl\`{e}mes fondamentaux de la th\'{e}orie des nombres \`{a} l'heure actuelle~; il 
se peut qu'une telle interpr\'{e}tation renferme la clef de l'hypoth\`{e}se de 
Riemann~$\ldots$''.
\smallskip
 
\noindent Galois groups are by construction 
projective limits of the finite groups attached to finite extensions. To  
get connected groups one clearly needs to relax this finiteness condition 
which is the same as the finite dimensionality of the central simple 
algebras of the Brauer theory. 
Since Archimedian places of $k$ are responsible for the non triviality of $D_k$ 
it is natural to ask the following preliminary question,

\smallskip

\noindent ``Is there a non trivial Brauer theory of central simple algebras over $\Cb$.''

\smallskip

\noindent We showed in \cite{ccc} that the {\it approximately finite 
dimensional} simple central algebras over $\Cb$ (called factors) provide a satisfactory answer to this question.
They are classified by their 
module,
$$
{\rm Mod} (M) \mathop{\subset}_{\sim} \ \Rb_+^* \, ,
$$
which is a virtual closed subgroup of $\Rb_+^*$.

 \noindent One can in fact go much further and understand that the renormalization group, 
once properly formulated mathematically as we did in section XVI, really appears as 
a perfect ambiguity group between solutions to a (physics) problem. It hence plays a 
role very similar to that of the Galois group of an algebraic equation and is an ideal
candidate for the missing Galois group at the Archimedian place. One can explore this 
idea further by making
use of the relation between the (conjectural) Hopf algebra
of Euler-Zagier numbers (\cite{sasha},  \cite{zagier}) and the Kreimer Hopf algebra.

\section{Appendix, the cyclic category}
At the conceptual level, cyclic cohomology is a way to embed the nonadditive 
category of 
algebras and
algebra homomorphisms in an additive category of modules. The latter is the 
additive 
category of
$\Lambda$-modules where $\Lambda$ is the cyclic category. Cyclic cohomology is 
then 
obtained as an $Ext$
functor (\cite{C2}).

\noindent The cyclic category is a small category which can be defined by
generators and relations. It has the same objects as the small
category $\D$ of totally ordered finite sets and increasing maps which plays a key 
role in simplicial topology. Let
us recall that $\D$ has one object $[n]$ for each
integer $n$, and is generated by faces $\delta_i, [n-1] \ra [n]$ (the
injection that misses $i$), and degeneracies $\s_j,[n+1] \ra [n] $
(the surjection which identifies $j$ with $j+1$), with the relations,
\begin{equation}\label{ad}
\delta_j  \delta_i = \delta_i  \delta_{j-1}
\ \hbox{for} \ i < j  , \ \s_j  \s_i = 
\s_i  \s_{j+1} \qquad i \leq j 
\end{equation}
\begin{equation}\label{ax}
\s_j  \delta_i = \left\{ \matrix{
\delta_i  \s_{j-1} \hfill &i < j \hfill \cr
1_n \hfill &\hbox{if} \ i=j \ \hbox{or} \ i = j+1 \cr
\delta_{i-1}  \s_j \hfill &i > j+1  . \hfill \cr
} \right.
\end{equation}
To obtain the cyclic category $\Lambda$ one adds for each $n$ a new morphism $\tau_n, [n]
\ra [n]$ such that,
\begin{equation}\label{ae}
\matrix{
\tau_n  \delta_i = \delta_{i-1}  
\tau_{n-1} &1 \leq i \leq n , &\tau_n  \delta_0 = 
\delta_n \hfill \cr
\cr
\tau_n  \s_i = \s_{i-1} 
\tau_{n+1} &1 \leq i \leq n , &\tau_n  \s_0 = 
\s_n  \tau_{n+1}^2 \cr
\cr
\tau_n^{n+1} = 1_n  . \hfill \cr
} 
\end{equation}

\noindent The original definition of $\Lambda$ (cf.~\cite{C2})
used homotopy classes of non decreasing maps from $S^1$ to $S^1$ of
degree~1, mapping $\Zb / n$ to $\Zb / m$ and is trivially equivalent
to the above. 

\noindent Given an algebra $A$ one obtains a module over the small category
$\Lambda$ by assigning to each integer $n \geq 0$ the vector space
$C^n$ of $n+1$-linear forms $\vp (x^0 , \ldots , x^n)$ on $A$, while
the basic operations are given by
\begin{equation}\label{ag}
\matrix{
(\delta_i  \vp) (x^0 , \ldots , x^n) &=& \vp (x^0 , \ldots , x^i
x^{i+1} , 
\ldots , x^n), \quad i=0,1,\ldots , n-1  \cr
\cr
(\delta_n  \vp) (x^0 , \ldots , x^n) &=& \vp (x^n  x^0 , x^1 , \ldots
, 
x^{n-1}) \hfill \cr
\cr
(\s_j  \vp) (x^0 , \ldots , x^n) &=& \vp (x^0 , \ldots , x^j , 1 ,
x^{j+1} 
, \ldots , x^n), \quad j=0,1,\ldots , n  \cr
\cr
(\tau_n  \vp) (x^0 , \ldots , x^n) &=& \vp (x^n , x^0 , \ldots ,
x^{n-1}) 
 . \hfill \cr
}
\end{equation}
\noindent These operations satisfy the relations (\ref{ad}) (\ref{ax}) and (\ref{ae}). This
shows that any algebra $A$ gives rise canonically to a $\Lambda$-module and
allows \cite{C2,L} to interpret the cyclic cohomology groups $HC^n(A)$
as $Ext^n$ functors. All of the general properties of cyclic
cohomology such as the long exact sequence relating it to Hochschild
cohomology are shared by Ext of general $\Lambda$-modules and can be
attributed to the equality of the classifying space $B\Lambda$ of the
small category $\Lambda$ with the classifying space $BS^1$ of the
compact one-dimensional Lie group $S^1$.
One has
\begin{equation}\label{af}
B\Lambda = BS^1 =P_{\infty}(\Cb)
\end{equation}

\end{document}

\end{thebibliography}
\end{document}